\documentclass[12pt]{article}

\newcommand{\real}{{\bf R}}
\newcommand{\complex}{{\bf C}}
\newcommand{\intplus}{{\bf N}}
\newcommand{\allint}{{\bf Z}}
\renewcommand{\div}{\mathop{\rm div}}
\newcommand{\rot}{\mathop{\rm rot}}
\newcommand{\Lip}{{\rm Lip}}
\renewcommand{\span}{{\rm span}}

\newcommand{\e}{{\rm e}}		
\renewcommand{\d}{\,{\rm d}}		
\newcommand{\D}{{\rm d}}		
\def\I{{\rm i}}				
\newcommand{\tr}{{\rm T}}
\newcommand{\mes}{{\rm mes}}

\newcommand{\cL}{{\cal L}}
\newcommand{\cN}{{\cal N}}
\newcommand{\cM}{{\cal M}}
\newcommand{\cO}{{\cal O}}
\newcommand{\cR}{{\cal R}}
\newcommand{\cS}{{\cal S}}
\newcommand{\cF}{{\cal F}}
\newcommand{\cG}{{\cal G}}

\newcommand{\pp}{{\bf p}}
\newcommand{\uu}{{\bf u}}
\newcommand{\vv}{{\bf v}}
\newcommand{\xx}{{\bf x}}
\newcommand{\yy}{{\bf y}}
\renewcommand{\AA}{{\bf A}}
\newcommand{\BB}{{\bf B}}
\newcommand{\FF}{{\bf F}}
\newcommand{\bbeta}{{\mbox{\boldmath$\beta$}}}
\newcommand{\ggamma}{{\mbox{\boldmath$\gamma$}}}
\newcommand{\GGamma}{{\mbox{\boldmath$\Gamma$}}}
\newcommand{\eeta}{{\mbox{\boldmath$\eta$}}}
\newcommand{\xxi}{{\mbox{\boldmath$\xi$}}}
\newcommand{\oone}{{\bf 1}}

\newtheorem{theorem}{Theorem}[section]
\newtheorem{lemma}[theorem]{Lemma}

\newtheorem{proposition}[theorem]{Proposition}
\newtheorem{corollary}[theorem]{Corollary}
\newtheorem{remark}[theorem]{Remark}

\newcommand{\reff}[1]{(\ref{#1})}

\newcommand{\proof}{{\noindent \bf Proof:\ }}
\newcommand{\half}{{\frac{1}{2}}}
\newcommand{\chiro}{{\chi_{r_0}}}
\newcommand{\phit}{{\Phi_{\tau}^{r_0,m}}}
\newcommand{\phirm}{{\Phi_1^{r_0,m}}}
\newcommand{\phione}{{\Phi_1^{r_0,m}}}
\newcommand{\wsloc}{{W_s^{loc}}}
\newcommand{\wcloc}{{W_c^{loc}}}
\newcommand{\wca}{{W_c^{\alpha^*}}}
\newcommand{\Eca}{{E_c^{\alpha^*}}}
\newcommand{\tf}{{\tilde{f}}}
\newcommand{\tg}{{\tilde{g}}}
\newcommand{\tFF}{{\tilde{{\bf F}}}}
\newcommand{\app}{{\rm app}}

\def\eqdef{\buildrel\hbox{\small{def}}\over =}
\def\build#1_#2^#3{\mathrel{
  \mathop{\kern 0pt#1}\limits_{#2}^{#3}}}
\def\QED{\mbox{}\hfill QED}

\newdimen\texpscorrection
\newdimen\figcenter
\def\figurewithtex #1 #2 #3 #4 #5\cr{\null
  {\goodbreak\figcenter=\hsize\relax
  \advance\figcenter by -#4truecm
  \divide\figcenter by 2
  \begin{figure}[hbt]
  \vskip #3truecm\noindent\hskip\figcenter
  \includegraphics{#1}{\hskip\texpscorrection\input #2 }
  \vskip 0.8truecm\noindent \vbox{\noindent {\footnotesize #5}}
  \end{figure}}}
\def\point#1 #2 #3 {\rlap{\kern #1 truecm
\raise #2 truecm \hbox{#3}}}

\usepackage{amsfonts}

\begin{document}

\title{Invariant manifolds and the long-time asymptotics of the
Navier-Stokes and vorticity equations on $\real^2$}

\author{Thierry Gallay \\ 
Universit\'e de Paris-Sud \\ 
Math\'ematiques \\
B\^atiment 425 \\
F-91405 Orsay \\
France
\and
C. Eugene Wayne \\
Department of Mathematics \\
and  Center for BioDynamics \\
Boston University \\
111 Cummington Street\\
Boston, MA 02215, USA}
\date{February 23, 2001}

\maketitle
\begin{abstract}
We construct finite-dimensional invariant manifolds in the phase
space of the Navier-Stokes equation on ${\real^2}$ and show
that these manifolds control the long-time behavior of the solutions.
This gives geometric insight into the existing results on the asymptotics 
of such solutions and also allows one to extend those results in a number 
of ways.
\end{abstract}


\section{Introduction} \label{intro}

In the last decade and a half, starting with the 
work of M. Schonbek, the long-time behavior of solutions of the
Navier-Stokes equation (and the related vorticity equation) on unbounded 
spatial domains has been extensively studied. (See \cite{schonbek:1985}, 
\cite{giga:1988b}, \cite{schonbek:1991} \cite{carpio:1994}, \cite{carpio:1996},
\cite{fujigaki:2000}, \cite{schonbek:2000}, \cite{marcel:2000} and 
\cite{miyakawa:2000} for a small sampling of this literature.)  
This prior work used a variety of techniques including energy estimates, 
the Fourier splitting method, and a detailed analysis of the semigroup of 
the linear part of the equation.  In the present paper we introduce another 
approach, based on ideas from the theory of dynamical systems, to compute these
asymptotics.  We prove that there exist finite-dimensional invariant 
manifolds in the phase space of these equations, and that all solutions 
in a neighborhood of the origin approach one of these manifolds with a rate 
that can be easily computed. Thus, computing the asymptotics of solutions 
up to that order is reduced to the simpler task of determining the 
asymptotics of the system of {\em ordinary} differential equations that
result when the original partial differential equation is restricted to the 
invariant manifold.  Although
it is technically quite different from our work, Foias and
Saut, \cite{foias:1984}, have also used invariant manifold theory to study the
long-time behavior of solutions of the Navier-Stokes equation
in a {\em bounded} domain.
 
As we shall see, if one specifies some given order of decay in advance, 
say $\cO(t^{-k})$, one can find a (finite-dimensional) invariant manifold 
such that all small solutions approach the manifold with at least that 
rate. Therefore, one can in principle compute the asymptotics of solutions 
to any order with this method. In Subsection~\ref{oseen}, for instance, 
we illustrate how these ideas can be used to extend known results about 
the stability of the Oseen vortex. In particular, we show that as 
solutions approach the vortex solution their velocity fields have a 
universal profile.

If one adopts a dynamical systems point of view towards the question of 
the long-time behavior of solutions of the Navier-Stokes equation, 
many phenomena which have been investigated in finite-dimensional 
dynamics immediately suggest possible effects in the Navier-Stokes 
equation.  As an example, one knows that in ordinary differential equations 
the presence of resonances between the eigenvalues of the linearized equations
may produce dramatic changes in the asymptotic behavior of solutions.  
Exploiting this observation in Subsection~\ref{logarithm} we construct 
solutions of the Navier-Stokes and vorticity equations whose asymptotic 
expansion contains logarithmic terms in time.

Finally, we feel that the geometric insights that the invariant
manifold method provides are very valuable.  As an example of
the sorts of insights this method provides, we reexamine the results of
Miyakawa and Schonbek \cite{miyakawa:2000} on optimal decay rates 
in Subsection~\ref{optimal}. We find that the moment conditions derived 
in \cite{miyakawa:2000} have a very simple geometric interpretation -- 
they are the analytic expression of the requirement that the solutions 
lie on certain invariant manifolds.

To explain our results and methods in somewhat more detail, we
recall that the Navier-Stokes equation in $\real^2$ has the form
\begin{equation}\label{NS2}
\frac{\partial \uu}{\partial t} +
  (\uu \cdot \nabla)\uu = \Delta \uu - \nabla p\ , \quad
  \nabla \cdot \uu = 0\ ,
\end{equation}
where $\uu = \uu(x,t) \in \real^2$ is the velocity field, $p = p(x,t) 
\in \real$ is the pressure field, and $x \in \real^2$, $t \ge 0$. 
For simplicity, the kinematic viscosity has been rescaled to $1$ 
in Eq.\reff{NS2}. 

Our results are based on two ideas. The first observation is that it is 
easier to derive the asymptotics of \reff{NS2} by working with the 
vorticity formulation of the problem. That is, we set $\omega = \rot \uu
= \partial_1 u_2 - \partial_2 u_1$ and study the equation \reff{V2} for 
$\omega$ rather than \reff{NS2} itself. One can then recover the solution 
of the Navier-Stokes equation via the Biot-Savart law \reff{BS2} (see
Lemma~\ref{HLS2} and Appendix~\ref{velocity2d} for estimates of various 
$L^p$ norms of the velocity field in terms of the vorticity.) While there 
have been some studies of the asymptotics
of solutions of the vorticity equation (notably \cite{carpio:1994}
and \cite{giga:1988b}), the relationship between the vorticity and 
velocity does not seem to have been systematically exploited to study 
the asymptotics of the Navier-Stokes equations. We feel that approaching 
the problem through the vorticity equation has a significant advantage.  
It has been realized for some time that the decay rate in time of solutions 
of \reff{NS2} is affected by the decay rate in space of the initial data.
For instance, general solutions of the Navier-Stokes equation in 
$L^2(\real^2)$ with initial data in $L^2 \cap L^1$ have $L^2$ norm decaying 
like $C t^{-1/2}$ as $t \to +\infty$. However, it follows from 
Wiegner's result \cite{wiegner:1987} that, if the initial data satisfy
$\int_{\real^2}(1+|x|)|\uu_0(x)| \d x < \infty$, this decay rate
can be improved to $C t^{-1}$. However, the semiflow generated by the 
Navier-Stokes equation does not preserve such a condition. For example, in 
Subsection~\ref{optimal}, we prove that there exist solutions of \reff{NS2}  
which satisfy $\int_{\real^2} (1+|x|)|\uu(x,t)| \d x < \infty$ 
when $t=0$, but not for later times. The vorticity equation does not
suffer from this drawback -- as we demonstrate in Section~\ref{scaling}, 
solutions of the vorticity equation in $\real^2$ which lie in weighted 
Sobolev spaces at time zero remain in those spaces for all time. 
A similar result holds for small solutions of the vorticity
equation in $\real^3$, see \cite{gallay:2001}. As we will see in later 
sections, these decay conditions are crucial for determining the asymptotic 
behavior of the solutions, and the fact that the vorticity equation 
preserves them makes it natural to work in the vorticity formulation.
It should be noted, however, that this approach requires the vorticity
to decay sufficiently rapidly as $|x| \to \infty$ and this is not assumed in
\cite{wiegner:1987}, for instance.  On the other hand, 
assuming that the vorticity decreases rapidly at infinity is
very reasonable from a physical point of view. This is the case, for instance,
if the initial data are created by stirring the fluid with
a (finite size) tool.  Note that in contrast, even very localized
stirring may not result in a velocity field that decays rapidly
as $x$ tends to infinity.  As an example, in Subsection~\ref{oseen}
we discuss the Oseen vortex which is a solution of \reff{NS2} whose
vorticity is Gaussian, but whose velocity field is not even in 
$L^2$.

The second idea that helps to understand
the long-time asymptotics of \reff{NS2} is the introduction of scaling
variables (see \reff{omega-w} and \reff{u-v} below).  It has, of course, been 
realized for a long-time that in studying the long-time asymptotics of 
parabolic equations it is natural to work with rescaled spatial variables. 
However, these variables do not seem to have been  used very much
in the context of the Navier-Stokes equation. (Though they are exploited 
in a slightly different context in \cite{cannone:1996}.)  For our work 
they offer a special advantage. If one linearizes the Navier-Stokes equation 
or vorticity equation around the zero solution, the resulting linear equation
has continuous spectrum that extends all the way from minus infinity to 
the origin.  If one wishes to construct finite-dimensional invariant 
manifolds in the phase space of these equations which control the 
asymptotics of solutions, it is not obvious how to do this even for the 
linearized equation -- let alone for the full nonlinear problem.  
However, building on ideas of \cite{wayne:1997} we show that, when 
rewritten in terms of the scaling variables, the linearized operator has
an infinite set of eigenvalues with explicitly computable eigenfunctions,
and the continuous spectrum can be pushed arbitrarily far into the left 
half plane by choosing the weighted Sobolev spaces in which we work 
appropriately. (See Appendix \ref{spectrum} for more details.) We are 
then able to construct invariant manifolds tangent at the origin to 
the eigenspaces of the point spectrum of this linearized equation and 
exploit the ideas that have been developed in finite-dimensional dynamics 
to analyze the asymptotic behavior of solutions of these partial 
differential equations.

We conclude this introduction with a short survey of the remainder
of the paper.  Although our ideas are applicable in all dimensions
we consider here the case of fluids in two dimensions.
Because the vorticity is a scalar in this case many calculations
are simplified, and we feel that the central ideas of our method
can be better seen without being obscured by technical details.
Thus, in Section~\ref{2Dvorticity} we begin by surveying the
(well developed) existence and uniqueness theory for the 
two-dimensional vorticity equation.  In Section~\ref{scaling}
we introduce scaling variables and the weighted Sobolev spaces
that we use in our analysis.  We then show (see Theorem~\ref{invman}) 
that in this formulation there exist families of finite-dimensional 
invariant manifolds in the phase space of the problem, and that all 
solutions near the origin either lie in, or approach these manifolds 
at a computable rate. Furthermore, we obtain a geometrical 
characterization of solutions that approach the origin ``faster than 
expected'', see Theorem~\ref{strongman}. In Section~\ref{applications} 
we apply the invariant manifolds constructed in the previous section
to derive the above mentioned results about the long-time behavior of 
the vorticity and Navier-Stokes equations. In a companion paper 
\cite{gallay:2001}, we obtain similar results for the small solutions 
of the Navier-Stokes equation in three dimensions by a slightly 
different method. Finally, there are three appendices 
which derive a number of facts we need in the previous sections. 
In Appendix~\ref{spectrum}, we study the spectrum of the linear operator
$\cL$ appearing in the rescaled vorticity equation \reff{SV2}, and we
obtain sharp estimates on the semigroup it generates. In
Appendix~\ref{velocity2d}, we study in detail the relationship between
the velocity field $\uu$ and the associated vorticity $\omega$. 
In particular, we show how the spatial decay of $\uu$ is related 
to the moments of the vorticity $\omega$. Appendix~\ref{Gapp} 
derives a technical estimate on the invariant manifold constructed 
in Subsection~\ref{logarithm}.
\medskip

\noindent
{\bf Notation:} Throughout the paper we use boldface letters for
vector-valued functions, such as $\uu(x,t)$.  However, to avoid
a proliferation of boldface symbols we use standard italic
characters for points in $\real^2$, such as $x=(x_1,x_2)$.
In both cases, $|\cdot |$ denotes the Euclidean norm in
$\real^2$.  For any $p\in [1,\infty]$ we denote by $| f |_p$
the norm of a function in the Lebesgue space $L^p(\real^2)$.
If ${\bf f}\in (L^p(\real^2))^2$, we set $| f |_p
= |~|f|~|_p$.  Weighted norms play an important role in this
paper. We introduce the weight function $b:\real^2 \to \real$
defined by $b(x) = (1+|x|^2)^{\half}$.  For
any $m\ge 0$, we set $\| f \|_m = |b^m f|_2$, and denote the
resulting Hilbert space by $L^2(m)$.  If $f \in
C^0([0,T],L^p(\real^2))$, we often write $f(\cdot,t)$ or
simply $f(t)$ to denote the map $x \to f(x,t)$.  Finally,
we denote by $C$ a generic positive constant, which may differ from
place to place, even in the same chain of inequalities.
\medskip

\noindent
{\bf Acknowledgements.} Part of this work was done when C.E.W.
visited the University of Paris-Sud and Th.G. the Department of 
Mathematics and Center for BioDynamics of Boston University.  The
hospitality of both institutions is
gratefully acknowledged.
We also thank I.~Gallagher, A.~Mielke, G.~Raugel, J.-C.~Saut, 
M.~Vishik, and P.~Wittwer for stimulating discussions.  We are
especially indebted to A.~Mielke for bringing to our attention
the work of \cite{miyakawa:2000}, which triggered our
interest in this problem, and to M.~Vishik who first suggested to
one of us that the ideas of \cite{wayne:1997} might be
useful in the context of the Navier-Stokes equation.  The research
of C.E.W. is supported in part by the NSF under grant 
DMS-9803164.


\section{The Cauchy problem for the vorticity equation}
\label{2Dvorticity}

In this section we describe existence and uniqueness results for
solutions of the vorticity equation.  As stressed in the
introduction, our approach is to study in detail the behavior of
solutions of the vorticity equation, and then to derive information 
about the solutions of the Navier-Stokes equation as a corollary.
The results in this section are not new and are reproduced here for easy 
reference.

In two dimensions, the vorticity equation is
\begin{equation}\label{V2}
  \omega_t + (\uu \cdot \nabla) \omega = \Delta \omega\ , 
\end{equation}
where $\omega = \omega(x,t) \in \real$, $x = (x_1,x_2) \in \real^2$, 
$t \ge 0$. The velocity field $\uu$ is defined in terms of the vorticity 
via the Biot-Savart law
\begin{equation}\label{BS2}
  \uu(x) = \frac{1}{2\pi} \int_{\real^2} 
  \frac{(\xx - \yy)^{\perp}}{|x -y|^2} \omega(y)\d y\ ,
  \quad x \in \real^2\ .
\end{equation}
Here and in the sequel, if $x=(x_1,x_2) \in \real^2$, we denote
$\xx = (x_1,x_2)^\tr$ and $\xx^{\perp} = (-x_2,x_1)^\tr$. 

\medskip
The following lemma collects useful estimates for the velocity $\uu$ in 
terms of $\omega$. 

\begin{lemma}\label{HLS2} Let $\uu$ be the velocity field obtained from
$\omega$ via the Biot-Savart law \reff{BS2}.\\
{\bf (a)} Assume that $1 < p < 2 < q < \infty$ and $\frac{1}{q} = 
\frac{1}{p} - \half$. If $\omega \in L^p(\real^2)$, then $\uu \in 
L^q(\real^2)^2$, and there exists $C>0$ such that
\begin{equation}\label{HLS}
  |\uu|_q \le C |\omega|_p\ .
\end{equation}
{\bf (b)} Assume that $1 \le p < 2 < q \le \infty$, and define $\alpha \in 
(0,1)$ by the relation $\half = \frac{\alpha}{p} + \frac{1-\alpha}{q}$. 
If $\omega \in L^p(\real^2) \cap L^q(\real^2)$, then $\uu \in 
L^\infty(\real^2)^2$, and there exists $C>0$ such that
\begin{equation}\label{interpol}
  |\uu|_\infty \le C |\omega|_p^\alpha |\omega|_q^{1-\alpha}\ .
\end{equation}
{\bf (c)} Assume that $1 < p < \infty$. If $\omega \in L^p(\real^2)$, then 
$\nabla \uu \in L^p(\real^2)^4$ and there exists $C>0$ such that
\begin{equation}\label{Calderon}
  |\nabla\uu|_p \le C |\omega|_p\ .
\end{equation}
In addition, $\div \uu = 0$ and $\rot \uu \equiv \partial_1 u_2 - 
\partial_2 u_1 = \omega$. 
\end{lemma}

\proof {\bf (a)} follows from the Hardy-Littlewood-Sobolev inequality, 
see for instance Stein \cite{stein:1970}, Chapter~V, Theorem~1. 
To prove {\bf (b)}, we remark that, for all $R > 0$, 
\begin{eqnarray*}
  |\uu(x)| &\le& \frac{1}{2\pi} \int_{|y|\le R} |\omega(x-y)|
  \frac{1}{|y|}\d y + \frac{1}{2\pi} \int_{|y|\ge R} |\omega(x-y)|
  \frac{1}{|y|}\d y \\
  &\le& C |\omega|_q R^{1-\frac{2}{q}} + C|\omega|_p 
   \frac{1}{R^{\frac{2}{p}-1}}\ ,
\end{eqnarray*}
by H\"older's inequality. Choosing $R = (|\omega|_p/|\omega|_q)^\beta$, 
where $\beta = \frac{\alpha}{1-2/q} = \frac{1-\alpha}{2/p-1}$, we obtain 
\reff{interpol}. Finally, \reff{Calderon} holds because
$\nabla \uu$ is obtained from $\omega$ via a singular integral kernel
of Calder\'on-Zygmund type, see \cite{stein:1970}, Chapter~II, Theorem~3.
\QED

\medskip
The next result shows that the Cauchy problem for \reff{V2} is globally 
well-posed in the space $L^1(\real^2)$.

\begin{theorem}\label{L1V2} For all initial data $\omega_0 \in L^1(\real^2)$, 
equation \reff{V2} has a unique global solution $\omega \in 
C^0([0,\infty),L^1(\real^2)) \cap C^0((0,\infty),L^\infty(\real^2))$
such that $\omega(0) = \omega_0$. Moreover, for all $p \in [1,+\infty]$, 
there exists $C_p > 0$ such that
\begin{equation}\label{lpvorticity}
 |\omega(t)|_p \le \frac{C_p |\omega_0|_1}{t^{1-\frac{1}{p}}}\ ,\quad t>0\ .
\end{equation}
Finally, the total mass of $\omega$ is preserved under the evolution:
$$
  \int_{\real^2} \omega(x,t) \d x = \int_{\real^2} \omega_0(x) \d x\ ,
  \quad t \ge 0\ .
$$
\end{theorem}

\proof The strategy of the proof is to rewrite \reff{V2} as an 
integral equation (see \reff{V2int} below), and then to solve this equation
using a fixed point argument in some appropriate function space. We refer
to Ben-Artzi \cite{ben-artzi:1994} and Brezis \cite{brezis:1994} for
details. \QED

\begin{remark} $L^1(\real^2)$ is not the largest space in which the
Cauchy problem for equation \reff{V2} can be solved. For instance, it
is shown in \cite{giga:1988}, \cite{giga:1988b} that \reff{V2} has a 
global solution for any $\omega_0 \in \cM(\real^2)$, the set of all finite 
measures on $\real^2$. However, uniqueness of this solution is known only 
if the atomic part of $\omega_0$ is sufficiently small. Moreover, 
such a solution is smooth and belongs to $L^1(\real^2)$ for any 
$t > 0$. Thus, since we are interested in the long-time behavior 
of the solutions, there is no loss of generality in assuming that
$\omega_0 \in L^1(\real^2)$. 
\end{remark}

If $\omega(t)$ is the solution of \reff{V2} given by Theorem~\ref{L1V2}, 
it follows from Lemma~\ref{HLS2} that the velocity field $\uu(t)$ 
constructed from $\omega(t)$ satisfies $\uu \in C^0((0,\infty),
L^q(\real^2)^2)$ for all $q \in (2,\infty]$, and that there exist 
constants $C_q > 0$ such that
\begin{equation}\label{lqvelocity}
  |\uu(t)|_q \le \frac{C_q |\omega_0|_1}{t^{\half - \frac{1}{q}}}\ ,
  \quad t>0\ .
\end{equation}
Moreover, one can show that $\uu(t)$ is a solution of the integral 
equation associated to \reff{NS2}. 

In \reff{lqvelocity}, the fact that we cannot bound the $L^2$ norm of the
velocity field is not a technical restriction.  As we will see,
even if $\omega(x)$ is smooth and rapidly decreasing, the velocity field 
$\uu(x)$ given by \reff{BS2} may not be in $L^2(\real^2)^2$. A typical 
example is the so-called {\sl Oseen vortex}, which will be studied in
Section~\ref{oseen} below. In fact, it is not difficult to verify that, 
if $\uu \in L^2(\real^2)^2$ and $\omega = \rot \uu \in L^1(\real^2)$, then
necessarily $\int_{\real^2} \omega(x)\d x = 0$. In this situation, 
the decay estimates \reff{lpvorticity}, \reff{lqvelocity} can be 
improved as follows:

\begin{theorem}\label{zeromass} Assume that $\uu_0 \in L^2(\real^2)^2$
and that $\omega_0 = \rot \uu_0 \in L^1(\real^2)$. Let $\omega(t)$ 
be the solution of \reff{V2} given by Theorem~\ref{L1V2}. Then
\begin{equation}\label{omegadecay}
  \lim_{t\to\infty} t^{1-\frac{1}{p}} |\omega(t)|_p = 0\ , \quad
  1 \le p \le \infty\ .
\end{equation}
If $\uu(t)$ is the velocity field obtained from $\omega(t)$ via the 
Biot-Savart law \reff{BS2}, then
\begin{equation}\label{uudecay}
  \lim_{t\to\infty} t^{\frac{1}{2}-\frac{1}{q}} |\uu(t)|_q = 0\ , \quad
  2 \le q \le \infty\ .
\end{equation}
\end{theorem}

\proof The decay estimate \reff{omegadecay} is a particular case of 
Theorem~1.2 in \cite{carpio:1994}. However, Theorem~\ref{zeromass} can 
also be proved by the following simple argument. Since $\uu(t)$ is
a solution of the Navier-Stokes equation in $L^2(\real^2)^2$, it is 
well-known that $\nabla \uu \in L^2((0,+\infty),L^2(\real^2)^4)$. 
In particular, $\int_0^\infty |\omega(t)|_2^2\d t < \infty$. On the other 
hand, since $\frac{\D}{\D t} |\omega(t)|_2^2 = -2|\nabla \omega(t)|_2^2$, 
the function $t \mapsto |\omega(t)|_2$ is non-increasing. Therefore, 
$$
   t|\omega(t)|_2^2 \le 2 \int_{\frac{t}{2}}^t |\omega(s)|_2^2\d s \eqdef
   \epsilon(t)^2 \,\build\hbox to 10mm{\rightarrowfill}_{t\to +\infty}^{}\,
   0\ ,
$$
which proves \reff{omegadecay} for $p = 2$. 

We now use the integral equation satisfied by $\omega(t)$, namely
\begin{equation}\label{V2int}
  \omega(t) = \e^{t\Delta}\omega_0 -\int_0^t \nabla\cdot \e^{(t-s)\Delta}
  \bigl(\uu(s)\omega(s)\bigr)\d s = \omega_1(t) + \omega_2(t)\ .
\end{equation}
Since $\omega_0 \in L^1(\real^2)$ and $\int_{\real^2}\omega_0(x)\d x = 0$,
a direct calculation shows that $|\omega_1(t)|_1 \to 0$ as $t \to +\infty$.
On the other hand, 
\begin{eqnarray*}
  |\omega_2(t)|_1 &\le& \int_0^t \frac{C}{\sqrt{t-s}}|\uu(s)\omega(s)|_1
   \d s \le \int_0^t \frac{C}{\sqrt{t-s}}|\uu(s)|_2 |\omega(s)|_2 \d s \\
  &\le& \int_0^t \frac{C|\uu_0|_2}{\sqrt{t-s}}\frac{\epsilon(s)}{\sqrt{s}}
  \d s = C |\uu_0|_2 \int_0^1 \frac{\epsilon(tx)}{\sqrt{x(1-x)}} \d x 
  \,\build\hbox to 10mm{\rightarrowfill}_{t\to +\infty}^{}\,0\ ,
\end{eqnarray*}
by Lebesgue's dominated convergence theorem. Thus, \reff{omegadecay} holds 
for $p = 1$, hence for $1 \le p \le 2$ by interpolation. Next, using 
the integral equation again and proceeding as above, it is straighforward 
to verify that \reff{omegadecay} holds for all $p \in [1,+\infty]$. From
Lemma~\ref{HLS2}, we then obtain \reff{uudecay} for $2 < q \le \infty$.
Finally, the fact that $|\uu(t)|_2 \to 0$ as $t \to +\infty$ is established
in \cite{kato:1984}. \QED


\section{Scaling Variables and Invariant Manifolds}\label{scaling}

Our analysis of the long-time asymptotics of \reff{V2} depends on 
rewriting the equation in terms of ``scaling variables'' or 
``similarity variables'':
$$
  \xi = \frac{x}{\sqrt{1+t}}\ , \qquad \tau = \log(1+t)\ .
$$
If $\omega(x,t)$ is a solution of \reff{V2} and if $\uu(t)$ is the 
corresponding velocity field, we introduce new functions $w(\xi,\tau)$, 
$\vv(\xi,\tau)$ by
\begin{eqnarray}\label{omega-w}
  \omega(x,t) &=& \frac{1}{1+t} \,w\Bigl(\frac{x}{\sqrt{1+t}},
  \log(1+t)\Bigr)\ , \\ \label{u-v}
  \uu(x,t) &=& \frac{1}{\sqrt{1+t}} \,\vv\Bigl(\frac{x}{\sqrt{1+t}},
  \log(1+t)\Bigr)\ .
\end{eqnarray}
Then $w(\xi,\tau)$ satisfies the equation
\begin{equation}\label{SV2}
  \partial_\tau  w = \cL w - (\vv \cdot \nabla_\xi) w\ ,
\end{equation}
where
\begin{equation}\label{Ldef}
  \cL w = \Delta_\xi w + \half (\xxi \cdot \nabla_\xi)w +w
\end{equation}
and 
\begin{equation}\label{SBS2}
  \vv(\xi,\tau) = \frac{1}{2\pi} \int_{\real^2}
  \frac{(\xxi - \eeta)^{\perp}}{|\xi - \eta |^2} w(\eta,\tau) 
  \d\eta\ .
\end{equation}
Scaling variables have been previously used to study the evolution
of the vorticity by Giga and Kambe \cite{giga:1988b} and 
Carpio \cite{carpio:1994}.

The results of the previous section already provide us with some 
information about solutions of \reff{SV2}.  For example, for all $w_0 \in
L^1(\real^2)$, there exists a unique solution $w \in C^0([0,\infty),
L^1(\real^2)) \cap C^0((0,\infty),L^\infty(\real^2))$ of \reff{SV2} 
such that $w(0) = w_0$. Translating \reff{lpvorticity} and \reff{lqvelocity}
into rescaled variables, we see that, for all $\tau > 0$, 
\begin{eqnarray}\label{wlpest}
  |w(\tau)|_p &\le& \frac{C_p |w_0|_1}{a(\tau)^{(1-\frac{1}{p})}}\ , 
    \quad 1 \le p \le \infty\ , \\ \label{vlqest}
  |\vv(\tau)|_q &\le& \frac{C_q |w_0|_1}{a(\tau)^{(\half - \frac{1}{q})}}\ ,
    \quad 2 < q \le \infty\ ,
\end{eqnarray}
where $a(\tau) = 1-\e^{-\tau}$. Moreover, if $\vv_0 = \vv(\cdot,0) \in 
L^2(\real^2)^2$, then $\int_{\real^2} w(\xi,\tau)\d\xi = 0$ for all 
$\tau \ge 0$, and it follows from \reff{omegadecay} that $|w(\tau)|_p \to 0$ 
as $\tau \to +\infty$ for all $p \in [1,+\infty]$. 

However, in order to derive the long-time asymptotics of these
solutions we need to work not only in ``ordinary'' $L^p$ spaces, but
also in {\sl weighted} $L^2$ spaces. For any $m \ge 0$, we define the 
Hilbert space $L^2(m)$ by
\begin{eqnarray*}
  L^2(m) &=& \bigl\{ f \in L^2(\real^2) \,|\, \|f\|_m < \infty \bigr\}\ ,
  \quad \hbox{where}\\ \nonumber
  \|f\|_m &=& \left(\int_{\real^2}(1+|\xi|^2)^m |f(\xi)|^2 \d\xi \right)^{1/2}
  = |b^m f|_2\ ,
\end{eqnarray*}
where $b(\xi) = (1+|\xi|^2)^{1/2}$. If $m > 1$, then $L^2(m) \hookrightarrow 
L^1(\real^2)$. In this case, we denote by $L^2_0(m)$ the closed subspace 
of $L^2(m)$ given by
$$
  L^2_0(m) = \Bigl\{w \in L^2(m)\,\Big|\, \int_{\real^2} w(\xi)\d\xi = 0
  \Bigr\}\ .
$$
Weighted Sobolev spaces can be defined in a similar way. For instance, 
$$
  H^1(m) = \bigl\{w \in L^2(m) \,|\, \partial_i f \in L^2(m) \hbox{ for }
  i = 1,2 \bigr\}\ .
$$

As is well known, the operator $\cL$ is the generator of a strongly continuous 
semigroup $\e^{\tau\cL}$ in $L^2(m)$, see Appendix~\ref{spectrum}. 
Note that, unlike $\e^{t\Delta}$, the semigroup $\e^{\tau\cL}$ does not 
commute with space derivatives. Indeed, since $\partial_i \cL = 
(\cL+\frac{1}{2})\partial_i$ for $i = 1,2$, we have $\partial_i \e^{\tau\cL}
= \e^\frac{\tau}{2} \e^{\tau\cL}\partial_i$ for all $\tau \ge 0$. Using 
this remark and the fact that $\nabla\cdot\vv = 0$, we can rewrite 
equation \reff{SV2} in integral form:
\begin{equation}\label{winteq}
  w(\tau) = \e^{\tau\cL}w_0 - \int_0^\tau \e^{-\frac{1}{2}(\tau-s)}
  \nabla\cdot \e^{(\tau-s)\cL} (\vv(s)w(s))\d s\ ,
\end{equation}
where $w_0 = w(0) \in L^2(m)$. The following lemma shows that the 
quadratic nonlinear term in \reff{winteq} is bounded (hence smooth) 
in $L^2(m)$ if $m > 0$. 

\begin{lemma}\label{Restim} Fix $m > 0$ and $T > 0$. Given $w_1, w_2 \in 
C^0([0,T],L^2(m))$, define
$$
  R(\tau) = \int_0^\tau \nabla\cdot \e^{(\tau-s)\cL} (\vv_1(s)w_2(s))\d s\ ,
   \quad 0 \le \tau \le T\ ,
$$
where $\vv_1$ is obtained from $w_1$ via the Biot-Savart law \reff{SBS2}.
Then $R \in C^0([0,T],L^2(m))$, and there exists $C_0 = C_0(m,T) > 0$ 
such that
$$
  \sup_{0 \le \tau \le T}\|R(\tau)\|_m \le C_0
  \Bigl(\,\sup_{0 \le \tau \le T}\|w_1(\tau)\|_m\Bigr)
  \Bigl(\,\sup_{0 \le \tau \le T}\|w_2(\tau)\|_m\Bigr)\ .
$$
Moreover, $C_0(m,T) \to 0$ as $T \to 0$. 
\end{lemma}

\proof Choose $q \in (1,2)$ such that $q > \frac{2}{m+1}$. Using 
Proposition~\ref{LpLqestim} (with $N = p = 2$), we obtain
$$
  |b^m R(\tau)|_2 \le C \int_0^\tau \frac{1}{a(\tau{-}s)^{(\frac{1}{q}
  -\frac{1}{2})+\frac{1}{2}}}\,|b^m \vv_1(s) w_2(s)|_q \d s\ ,
$$
where $a(\tau) = 1-\e^{-\tau}$ and $b(\xi) = (1+|\xi|^2)^{1/2}$. Moreover, 
using H\"older's inequality and Lemma~\ref{HLS2}, we have
$$
  |b^m \vv_1 w_2|_q \le |b^m w_2|_2 |\vv_1|_\frac{2q}{2-q} 
  \le C \|w_2\|_m |w_1|_q \le C\|w_2\|_m \|w_1\|_m\ ,
$$
since $L^2(m) \hookrightarrow L^q(\real^2)$. Therefore, for all $\tau \in
[0,T]$, we find
$$
  \|R(\tau)\|_m \le C \Bigl(\int_0^\tau \frac{1}{a(s)^\frac{1}{q}}
  \d s\Bigr) \Bigl(\,\sup_{0 \le s \le T}\|w_1(s)\|_m\Bigr)
  \Bigl(\,\sup_{0 \le s \le T}\|w_2(s)\|_m\Bigr)\ .
$$
This concludes the proof of Lemma~\ref{Restim}. \QED

\medskip
We now show that \reff{SV2} has global solutions in $L^2(m)$ if $m > 1$. 

\begin{theorem} \label{weightedL2} Suppose that $w_0 \in L^2(m)$ for some
$m > 1$.  Then \reff{SV2} has a unique global solution $w \in 
C^0([0,\infty),L^2(m))$ with $w(0) = w_0$, and there exists
$C_1 = C_1(\|w_0\|_m) > 0$ such that
\begin{equation}\label{wglobound}
  \|w(\tau)\|_m \le C_1\ , \quad \tau \ge 0\ .
\end{equation}
Moreover, $C_1(\|w_0\|_m) \to 0$ as $\|w_0\|_m \to 0$. Finally, if 
$w_0 \in L^2_0(m)$, then $\int_{\real^2} w(\xi,\tau)\d\xi = 0$ for all 
$\tau \ge 0$, and $\displaystyle{\lim_{\tau\to\infty} \|w(\tau)\|_m = 0}$. 
\end{theorem}

\proof Using Lemma~\ref{Restim} and a fixed point argument, it is 
easy to to show that, for any $K > 0$, there exists $\tilde T = \tilde 
T(K) > 0$ such that, for all $w_0 \in L^2(m)$ with $\|w_0\|_m \le K$, equation
\reff{winteq} has a unique local solution $w \in C^0([0,\tilde T],L^2(m))$. 
This solution $w(\tau)$ depends continuously on the initial data $w_0$, 
uniformly in $\tau \in [0,\tilde T]$. Moreover, $\tilde T$ can be chosen
so that $\|w(\tau)\|_m \le 2\|w_0\|_m$ for all $\tau \in [0,\tilde T]$.
Thus, to prove global existence, it is sufficient to show that any solution 
$w \in C^0([0,T],L^2(m))$ of \reff{winteq} satisfies the bound 
\reff{wglobound} for some $C_1 > 0$ (independent of $T$.)

Let $w_0 \in L^2(m)$, $T > 0$, and assume that $w \in C^0([0,T],L^2(m))$ 
is a solution of \reff{winteq}. Without loss of generality, we suppose
that $T \ge \tilde T \equiv \tilde T(\|w_0\|_m)$. Since $L^2(m) 
\hookrightarrow L^p(\real^2)$ for all $p \in [1,2]$, there exists 
$C > 0$ such that $|w(\tau)|_p \le C \|w(\tau)\|_m \le 2C \|w_0\|_m$
for all $\tau \in [0,\tilde T]$. By \reff{wlpest}, we also have 
$|w(\tau)|_p \le C_p a(\tilde T)^{-1+1/p}|w_0|_1 \le C\|w_0\|_m$
for all $\tau \in [\tilde T,T]$. Thus, there exists $C_2 = C_2(\|w_0\|_m) 
> 0$ such that $|w(\tau)|_p \le C_2$ for all $p \in [1,2]$ and all 
$\tau \in [0,T]$. Moreover, $C_2(\|w_0\|_m) \to 0$ as $\|w_0\|_m \to 0$.
To bound $||\xi|^m w(\tau)|_2$, we compute
\begin{eqnarray}\nonumber
 &&\half \frac{\D}{\D\tau} \int_{\real^2} |\xi|^{2m} w(\xi,\tau)^2 d \xi
 = \int |\xi|^{2m} w (\partial_\tau  w) \d\xi \\ \label{dum1}
 &&\qquad\quad = \int |\xi|^{2m} \{ w \Delta w + \frac{w}{2} (\xxi \cdot 
\nabla) w + w^2 - w(\vv \cdot \nabla)w\} \d\xi\ .
\end{eqnarray}
Integrating by parts and using the fact that $\div \vv = 0$, we can 
rewrite
\begin{eqnarray*}
 \int |\xi|^{2m} w(\Delta w) \d\xi &=& -\int |\xi|^{2m} |\nabla w|^2 \d\xi 
 + 2 m^2 \int |\xi |^{2m-2} w^2 \d\xi\ ,\\
 \int |\xi |^{2m} \frac{w}{2} (\xxi \cdot \nabla) w \d\xi
 &=& -\frac{m+1}{2} \int|\xi |^{2m} w^2 \d\xi\ ,\\
 \int |\xi |^{2m} w (\vv \cdot \nabla) w \d\xi
 &=& \half \int |\xi |^{2m} (\vv \cdot \nabla) w^2 \d\xi
 = \half \int |\xi |^{2m} \nabla\cdot(\vv w^2) \d\xi\\
 &=& -m \int |\xi |^{2m-2} (\xxi \cdot \vv) w^2 \d\xi\ .
\end{eqnarray*}
Inserting these expressions into \reff{dum1}, we find:
\begin{eqnarray*}
 \half \frac{\D}{\D\tau} \int |\xi|^{2m} w^2 d \xi
 &=& -\int |\xi|^{2m} |\nabla w|^2 \d\xi -\frac{m-1}{2} \int|\xi |^{2m} 
 w^2 \d\xi\\
 &&+\ 2 m^2 \int |\xi |^{2m-2} w^2 \d\xi +m\int |\xi |^{2m-2} 
 (\xxi \cdot \vv) w^2 \d\xi\ .
\end{eqnarray*}
We next remark that, for all $\epsilon > 0$, there exists $C_\epsilon > 0$
such that
\begin{eqnarray*}
 \int|\xi |^{2m-2} w^2 \d\xi &\le& \epsilon \int 
 |\xi |^{2m} w^2 \d\xi + C_{\epsilon} \int w^2 \d\xi\ ,\\
 \Big| \int |\xi |^{2m-2} (\xxi \cdot \vv) w^2 \d\xi \Big|
 &\le& \epsilon \int |\xi |^{2m} w^2 \d\xi 
 + C_{\epsilon} |\vv |_\infty^{2m} \int  w^2 \d\xi\ .
\end{eqnarray*}
Note that, just as in \reff{interpol} of Lemma~\ref{HLS2}, 
$|\vv(\tau)|_{\infty} \le C |w|_p^{\alpha} |w|_q^{1-\alpha}$
where $1 < p < 2 < q <\infty$ and $\frac{\alpha}{p} +
\frac{1-\alpha}{q} = \half$. We already know that $|w(\tau)|_p \le C_2$
for all $\tau \in [0,T]$. Choosing $\alpha = 1-\frac{1}{8m}$, 
$q = 4 + \frac{1}{2m}$, and using \reff{wlpest} to bound $|w(\tau)|_q$,
we obtain
$$
  |\vv(\tau)|_\infty^{2m} \le C C_2^{2m\alpha} \Bigl(\frac{C_q 
  |w_0|_1}{a(\tau)^{(1-\frac{1}{q})}}\Bigr)^{2m(1-\alpha)}
  \le C_3 (1+\tau^{-3/16})\ , \quad 0 < \tau \le T\ ,
$$
where $C_3 = C_3(m,\|w_0\|_m)$. Thus, we see that for any $\delta \in 
(0,m-1)$, there exists $\epsilon > 0$ sufficiently small so that
\begin{eqnarray}\nonumber
  \frac{\D}{\D\tau} \int |\xi|^{2m} w^2 d \xi
  &\le& - 2\int |\xi |^{2m} |\nabla w|^2 \d\xi - \delta \int
  |\xi |^{2m} w^2 \d\xi\\ \label{diffineq}
  && \quad  +\ C_{\epsilon}(1+C_3)(1+\tau^{-3/16}) \int w^2 \d\xi\ .
\end{eqnarray}
Integrating this inequality, we obtain that $\int |\xi|^{2m} w^2 \d\xi 
\le C_4$ for all $\tau \in [0,T]$, where $C_4 = C_4(m,\|w_0\|_m) \to 0$ 
as $\|w_0\|_m \to 0$. This proves that all solutions of \reff{winteq} in 
$L^2(m)$ are global and satisfy \reff{wglobound}. Finally, if 
$w_0 \in L^2_0(m)$ and if $w_0$ is the velocity field obtained from 
$w_0$ via the Biot-Savart law \reff{SBS2}, then $\vv_0 \in L^2(\real^2)^2$
by Corollary~\ref{velint1}. In this case, we already observed that 
$|w(\tau)|_2$ converges to zero as $\tau \to +\infty$, and so does 
$||\xi|^m w(\tau)|_2$ by \reff{diffineq}. This concludes the proof of 
Theorem~\ref{weightedL2}. \QED

\medskip
We now explain the motivation for introducing the scaling variables
and discuss briefly how we will proceed.  The spectrum of the
operator $\cL$ acting on $L^2(m)$ consists of a sequence of
eigenvalues $\sigma_d = \{ -\frac{k}{2} ~|~ k=0,1,2,\dots\}$, plus 
continuous spectrum $\sigma_c = \{ \lambda \in \complex ~|~
\Re \lambda \le -(\frac{m-1}{2}) \}$, see Appendix~\ref{spectrum}. 
Ignoring the continuous spectrum for a moment (since we can ``push it
out of the way'' by choosing $m$ appropriately) we choose 
coordinates for $L^2(m)$ whose basis vectors are the eigenvectors
of $\cL$. These eigenvectors can be computed explicitly.
Expressed in these coordinates \reff{SV2} becomes an infinite
set of ordinary differential equations with a very simple linear
part and a quadratic nonlinear term.  In the study of dynamical
systems, powerful tools such as invariant manifold theory and
the theory of normal forms have been developed to study such
systems of equations.  We will apply these tools to study
\reff{SV2}.  In particular, we will show that given any
$\mu > 0$, the long time behavior of solutions of \reff{SV2}
up to terms of order $\cO (\e^{-\mu \tau})$ (which corresponds to the
behavior of solutions of \reff{V2} up to terms of order
$\cO(t^{-(\mu+1)}$ ) is given by a {\em finite} system of ordinary
differential equations which results from restricting
\reff{SV2} to a finite-dimensional, invariant manifold in $L^2(m)$.
This manifold is tangent at the origin to the eigenspace spanned by
the eigenvalues in $\sigma_d$ with eigenvalues bigger than $-\mu$.
Using these manifolds we will see that we can, at least in 
principle, calculate the long-time asymptotics to any order.
As an example we will show that in computing the asymptotics 
of the velocity field we obtain additional terms beyond those 
calculated in \cite{carpio:1996} and \cite{fujigaki:2000} due to the
fact that we work with the vorticity formulation. We also exhibit 
an example that shows that in general the asymptotic behavior of 
the solution of \reff{V2} contains terms logarithmic in $t$.

\figurewithtex 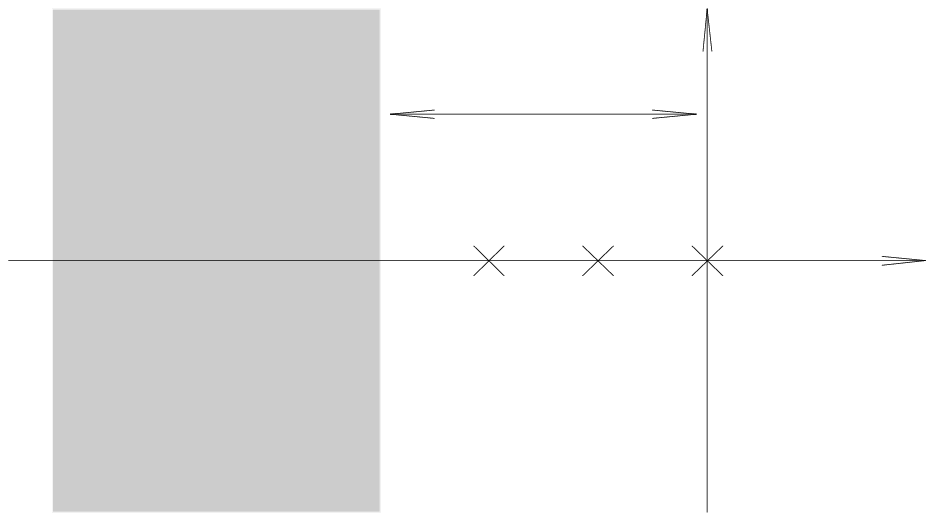 Fig1.tex 5.5 11 {\bf Fig. 1.} The spectrum
of the linear operator $\cL$ in $L^2(m)$, when $m = 4$.\cr

There are some technical difficulties which have to be overcome
in order to implement this picture.  The two most important ones
are the presence of the continuous spectrum of $\cL$, and
the fact that the nonlinear term in \reff{SV2} is not a 
smooth function from $L^2(m)$ to itself.  We circumvent these
problems by working not with the differential equation itself,
but rather with the semiflow defined by it.  For this semiflow
we can easily develop estimates which allow us to apply the
invariant manifold theorem of Chen, Hale and Tan \cite{chen:1997}.

In order to apply the invariant manifold theorem we consider the behavior
of solutions of \reff{SV2} in a neighborhood of the fixed point
$w = 0$.  Since we wish for the moment to concentrate only
on this neighborhood we cut-off the nonlinearity outside of
a neighborhood of size $r_0 > 0$.  More precisely, let $\chi : 
L^2(m) \to \real$ be a bounded $C^{\infty}$ function such that 
$\chi(w)=1$ if $\|w\|_m \le 1$ and $\chi(w)=0$ if $\|w\|_m \ge 2$.
For any $r_0 > 0$, we define $\chiro(w) = \chi(w/r_0)$ for all 
$w \in L^2(m)$. 

\begin{remark} Such a cut-off function $\chi$ may not exist for all Banach 
spaces. However, for a Hilbert space like $L^2(m)$ one can always find such
a function.
\end{remark} 

We now replace \reff{SV2} by
\begin{equation}\label{SV2c}
  \partial_\tau  w = \cL w - (\vv \cdot\nabla )(\chiro(w) w)\ .
\end{equation}
Note that \reff{SV2c} and \reff{SV2} are equivalent whenever
$\|w\|_m \le r_0$, but \reff{SV2c} becomes a linear equation
for $\|w\|_m \ge 2r_0$. Using Proposition~\ref{Sgestim} and the 
analogue of Lemma~\ref{Restim}, it is easy to show that all solutions
of \reff{SV2c} in $L^2(m)$ are global if $m > 0$, and stay bounded
if $m > 1$. In the sequel, we denote by $\phit$ the semiflow on $L^2(m)$ 
defined by \reff{SV2c}, and by $\Phi_\tau^m$ the semiflow associated
with \reff{SV2}. The following properties of $\phit$ will be useful:

\begin{proposition} \label{2Dsemiflow} Fix $m > 1$, $r_0 > 0$, 
and let $\phit$ be the semiflow on $L^2(m)$ defined by \reff{SV2c}.  
If $r_0 > 0$ is sufficiently small, then for $\tau=1$ we can decompose
\begin{equation}\label{phidef}
\phione = \Lambda + \cR\ ,
\end{equation}
where $\Lambda \equiv \e^{\cL}$ is a bounded linear operator on $L^2(m)$, 
and $\cR : L^2(m) \to L^2(m)$ is a $C^\infty$ map satisfying $\cR(0) = 0$, 
$D\cR(0) = 0$. Moreover, $\cR$ is globally Lipschitz, and there exists 
$L > 0$ (independent of $r_0$) such that $\Lip(\cR) \le Lr_0$. 
\end{proposition}

\proof Let $w_1,w_2 \in L^2(m)$, and define $w_i(\tau) = \phit w_i$ 
for $i = 1,2$ and $0 \le \tau \le 1$. Then
$$
   w_i(\tau) = \e^{\tau\cL}w_i - \int_0^\tau \e^{-\frac{1}{2}
   (\tau-s)}\nabla\cdot \e^{(\tau-s)\cL}\bigl(\vv_i(s)\chiro(w_i(s))
   w_i(s)\bigr)\d s\ , \quad i = 1,2\ ,
$$
where $\vv_i(\tau)$ is the velocity field obtained from $w_i(\tau)$ via 
the Biot-Savart law \reff{SBS2}. Proceeding as in the proof of 
Lemma~\ref{Restim}, we find that there exists $C \ge 1$ and $K > 0$
such that
\begin{equation}\label{hatestim}
  \sup_{0\le\tau\le 1} \|w_1(\tau) - w_2(\tau)\|_m \le 
  C\|w_1-w_2\|_m + K r_0 \sup_{0\le\tau\le 1} 
  \|w_1(\tau) - w_2(\tau)\|_m\ .
\end{equation}
Thus, assuming $Kr_0 \le 1/2$, we obtain $\|\phit w_1 -\phit w_2\|_m 
\le 2C\|w_1-w_2\|_m$ for all $\tau \in [0,1]$. We now define
$\cR = \phione -\Lambda$, so that
\begin{equation}\label{Rexpress}
  \cR(w_i) = -\int_0^1 \e^{-\frac{1}{2}(1-\tau)}\nabla\cdot
  \e^{(1-\tau)\cL}(\vv_i(\tau)\chiro(w_i(\tau))w_i(\tau))\d\tau\ , \quad
  i = 1,2\ .
\end{equation} 
In view of \reff{hatestim}, we have $\|\cR(w_1)-\cR(w_2)\|_m \le Lr_0
\|w_1-w_2\|_m$, where $L = 2CK$. Moreover, using \reff{Rexpress},
it is not difficult to show that $\cR : L^2(m) \to L^2(m)$ is smooth and 
satisfies $D\cR(0) = 0$. \QED

\medskip
In the rest of this section, we fix $k \in \intplus$ and we assume
that $m \ge k+2$. As is shown in Appendix~\ref{spectrum}, the spectrum 
of $\Lambda = \e^{\cL}$ in $L^2(m)$ is $\sigma(\Lambda) = \Sigma_c 
\cup \Sigma_d$, where $\Sigma_d = \{\e^{-n/2}\,|\, n=0,1,2\dots\}$ and 
$\Sigma_c = \{ \lambda \in \complex \,|\, |\lambda| \le 
\e^{-(\frac{m-1}{2})}\}$. Since $m-1 > k$, this means that $\Lambda$ has at 
least $k+1$ isolated eigenvalues $\lambda_0,\lambda_1,\dots,\lambda_k$, 
where $\lambda_j = \e^{-j/2}$. Let $P_k$ be the spectral projection onto
the (finite-dimensional) subspace of $L^2(m)$ spanned by the eigenvectors
of $\Lambda$ corresponding to the eigenvalues $\lambda_0,\dots,\lambda_k$, 
and let $Q_k = \oone - P_k$. Applying to the semiflow $\phit$ the invariant
manifold theorem as stated in Chen, Hale, and Tan \cite{chen:1997}, we 
obtain the following result:

\begin{theorem} \label{invman} Fix $k \in \intplus$, $m \ge k+2$, 
and choose $\mu_1, \mu_2 \in \real$ such that $\frac{k}{2} < \mu_1 < \mu_2
< \frac{k+1}{2}$. Let $P_k, Q_k$ be the spectral projections defined 
above, and let $E_c = P_k L^2(m)$, $E_s = Q_k L^2(m)$. Then, for $r_0 > 0$ 
sufficiently small, there exists a $C^1$ and globally Lipschitz map 
$g : E_c \to E_s$ with $g(0) = 0$, $Dg(0) = 0$, such that the submanifold
$$
  W_c = \{ w_c + g(w_c) ~|~ w_c \in E_c \}
$$
has the following properties:\\
{\bf 1.} {\rm (Invariance)} The restriction to $W_c$ of the 
semiflow $\phit$ can be extended to a Lipschitz flow on $W_c$.
In particular, $\phit(W_c) = W_c$ for all $\tau\ge 0$, and
for any $w_0 \in W_c$, there exists a unique negative
semi-orbit $\{w(\tau)\}_{\tau\le 0}$ in $W_c$ with $w(0) = w_0$.
If $\{w(\tau)\}_{\tau\le 0}$ is a negative semi-orbit contained in 
$W_c$, then
\begin{equation}\label{char1}
  \limsup_{\tau \to -\infty} \frac{1}{|\tau|} \ln \| w(\tau) \|_m 
  \le \mu_1\ .
\end{equation}
Conversely, if a negative semi-orbit of $\phit$ satisfies
\begin{equation}\label{char2}
  \limsup_{\tau\to -\infty} \frac{1}{|\tau|} \ln \| w(\tau) \|_m 
  \le \mu_2\ ,
\end{equation}
then it lies in $W_c$.\\
{\bf 2.} {\rm (Invariant Foliation)} There is a continuous map
$h: L^2(m) \times E_s \to E_c$ such that, for each $w \in W_c$, 
$h(w,Q_k(w)) = P_k(w)$ and the manifold $M_{w} = \{ h(w,w_s)+w_s ~|~ 
w_s \in E_s\}$ passing through $w$ satisfies $\phit (M_w) \subset 
M_{\phit(w)}$ and
\begin{equation}\label{whiskercontract}
  M_w = \Bigl\{\tilde{w} \in L^2(m) ~\Big|~ \limsup_{\tau\to +\infty}
  \frac{1}{\tau}\ln \|\phit(\tilde{w}) - \phit(w)\|_m \le -\mu_2\Bigr\}\ .
\end{equation}
Moreover, $h: L^2(m) \times E_s \to E_c$ is $C^1$ in the $E_s$ direction.\\
{\bf 3.} {\rm (Completeness)} For every $w \in L^2(m)$, $M_w \cap
W_c$ is exactly a single point.  In particular, $\{M_w\}_{w\in W_c}$ is a 
foliation of $L^2(m)$ over $W_c$.
\end{theorem}

\begin{remark}\label{smoothness}
Although the nonlinearity $\cR$ is $C^\infty$, the invariant manifold
$W_c$ is in general not smooth. However, for any $\rho > 1$ such that
$k\rho < k+1$, one can prove that the map $g : E_c \to E_s$ is of 
class $C^\rho$ if $r_0$ is sufficiently small. (See \cite{henry:1981},
Section~6.1.) Note also that the manifold $W_c$ is not unique, since it 
depends on the choice of the cut-off function $\chi$. This is the 
only source of non-uniqueness however -- once the cutoff function
has been fixed, the construction of \cite{chen:1997} produces
a unique global manifold. 
\end{remark}

\proof All we need to do is verify that hypotheses ({\bf H.1})--({\bf H.4}) 
of Theorem 1.1 in \cite{chen:1997} hold for the semiflow $\phit$. 
Assumption ({\bf H.1}) is satisfied because $w \mapsto \phit w$ 
is globally Lipschitz, uniformly in $\tau \in [0,1]$, see \reff{hatestim}.
Hypothesis ({\bf H.2}) is nothing but the decomposition $\phit = \Lambda + 
\cR$ obtained in Proposition~\ref{2Dsemiflow}. Assumption ({\bf H.4})
is a smallness condition on $\Lip(\cR)$, which is easily achieved by
taking $r_0 > 0$ sufficiently small. To verify ({\bf H.3}), we recall
that $L^2(m) = E_c \oplus E_s$, and we set $\Lambda_c = P_k \Lambda P_k$, 
$\Lambda_s = Q_k \Lambda Q_k$. Since $\sigma(\Lambda_c) = 
\{1,\e^{-\frac{1}{2}},\dots,\e^{-\frac{k}{2}}\}$ and all eigenvalues 
of $\Lambda_c$ are semisimple, it is clear that $\Lambda_c$ has a bounded 
inverse, and that there exists $C_c \ge 1$ such that 
$$
  \|\Lambda_c^{-j}w\|_m \le C_c \ \e^{j \frac{k}{2}}\|w\|_m\ , \quad
  \hbox{for all } j \in \intplus \hbox{ and all }w \in E_c\ .
$$
On the other hand, take $\epsilon = 0$ if $m > k+2$ or $\epsilon > 0$
arbitrarily small if $m = k+2$. By Proposition~\ref{Sgestim}, there 
exists $C_s \ge 1$ such that
$$
  \|\Lambda_s^j w\|_m \le C_s \ \e^{-j (\frac{k+1}{2}-\epsilon)}\|w\|_m\ , 
  \quad \hbox{for all } j \in \intplus \hbox{ and all }w \in E_s\ .
$$
These estimates are exactly what is required in ({\bf H.3}). \QED

\figurewithtex 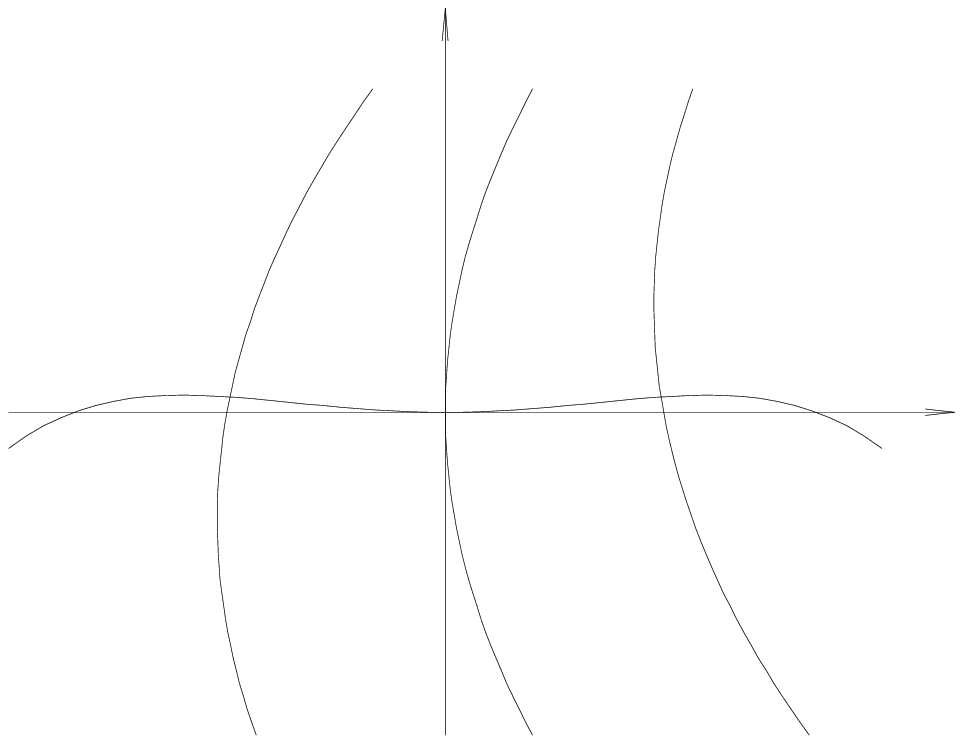 Fig2.tex 7.5 10.0 {\bf Fig. 2.} A picture 
illustrating the invariant manifold $W_c$ and the foliation 
$\{M_w\}_{w\in W_c}$ of $L^2(m)$ over $W_c$. The leaf $M_0$ 
of the foliation passing though the origin coincides with the
strong-stable manifold $W_s^{loc}$ in a neighborhood of the origin. 
Note that $\dim(E_c) < \infty$, whereas the subspace $E_s$ is 
infinite-dimensional.\cr

\medskip
By Theorem~\ref{weightedL2}, there exists $r_1 > 0$ such that 
all solutions of \reff{SV2} with $\|w(0)\|_m \le r_1$ satisfy 
$\|w(\tau)\|_m \le r_0$ for all $\tau \ge 0$. For such solutions the 
semiflow $\Phi_\tau^m$ defined by \reff{SV2} coincides with $\phit$,  
and thus the invariant manifolds in Theorem~\ref{invman} are also locally 
invariant with respect to $\Phi_\tau^m$. As a consequence, there exist a 
multitude of finite-dimensional invariant manifolds in the phase space of 
the vorticity equation \reff{SV2}, hence in the phase space of the 
Navier-Stokes equation.  Furthermore, points {\bf 2} and {\bf 3} in
Theorem~\ref{invman} give us an estimate of how rapidly solutions approach
these invariant manifolds. Since this result will be crucial in the 
applications, we state it here as a Corollary.

\begin{corollary}\label{locinvman}
Fix $k \in \intplus$, $m \ge k+2$, and let $W_c$ be the submanifold of 
$L^2(m)$ constructed in Theorem~\ref{invman}. Define
\begin{equation}\label{loccent}
  W_c^{loc} = W_c \cap \bigl\{w \in L^2(m)\,|\, \|w\|_m < r_0\bigr\}\ .
\end{equation}
Then $W_c^{loc}$ is locally invariant under the semiflow $\Phi_\tau^m$ 
defined by \reff{SV2}. If $\{w(\tau)\}_{\tau \le 0}$ is a negative 
semi-orbit of \reff{SV2} such that $\|w(\tau)\|_m < r_0$ for all 
$\tau \le 0$, then $w(\tau) \in W_c^{loc}$ for all $\tau \le 0$. 
Moreover, for any $\mu < \mu_2$ (where $\mu_2$ is as in Theorem~\ref{invman}),
there exist $r_2 > 0$ and $C > 0$ with the following property: for all 
$\tilde{w}_0 \in L^2(m)$ with $\|\tilde{w}_0\|_m \le r_2$, there exists 
a unique $w_0 \in W_c^{loc}$ such that $\Phi_\tau^m(w_0) \in W_c^{loc}$ 
for all $\tau \ge 0$ and
\begin{equation}\label{locattract}
  \|\Phi_\tau^m(\tilde{w}_0) - \Phi_\tau^m(w_0)\|_m \le C \e^{-\mu \tau}
  \ , \quad \tau \ge 0\ .
\end{equation}
\end{corollary}

\proof It follows immediately from point {\bf 1} in Theorem~\ref{invman} 
that $W_c^{loc}$ is locally invariant under the semiflow $\Phi_\tau^m$ 
and contains the negative semi-orbits that stay in a neighborhood of the 
origin. Choose $r_1 > 0$ so that $\|w_0\|_m \le r_1$ implies 
$\|\Phi_\tau^m(w_0)\|_m < r_0$ for all $\tau \ge 0$. By points {\bf 2} and 
{\bf 3} in Theorem~\ref{invman}, given $\tilde{w}_0 \in L^2(m)$, there 
exists a unique $w_0 \in M_{\tilde{w}_0} \cap W_c$. Setting $w_0 = w_c + 
g(w_c)$, we see from the definitions that $w_c$ solves the equation
\begin{equation}\label{fixedwc}
  w_c = h(\tilde{w}_0,g(w_c))\ .
\end{equation}
Since $h$ is continuous and $w_c \mapsto h(\tilde{w}_0,g(w_c))$ is 
Lipschitz (with a small Lipschitz constant), it is clear that \reff{fixedwc}
has a unique solution $w_c$ which depends continuously on $\tilde{w}_0$. 
Moreover, by \reff{whiskercontract}, $w_c = 0$ if $\tilde{w}_0 = 0$. 
Therefore, by continuity, there exists $r_2 \in (0,r_1]$ such that, if
$\|\tilde{w}_0\|_m \le r_2$, then $\|w_0\|_m = \|w_c + g(w_c)\|_m \le r_1$. 
In this case, $\max(\|\Phi_\tau^m(\tilde{w}_0)\|_m\,,\,\|\Phi_\tau^m(w_0)\|_m) 
< r_0$ for all $\tau \ge 0$, and \reff{locattract} follows from 
\reff{whiskercontract}. \QED

\medskip
In particular, given $\mu < \mu_2$, it follows from Corollary~\ref{locinvman}
that a solution $w(\tau)$ of \reff{SV2} on $W_c^{loc}$ cannot converge to 
zero faster than $\e^{-\mu\tau}$ as $\tau \to +\infty$, unless $w(\tau) 
\equiv 0$. For this reason, $W_c^{loc}$ is usually called the ``weak-stable'' 
manifold of the origin. In the applications, we will also be interested in 
solutions that approach the origin ``rapidly'', namely with a rate 
$\cO(\e^{-\mu_2 \tau})$ or faster. By Theorem~\ref{invman}, all such 
solutions lie in the leaf of the foliation $\{M_w\}_{w\in W_c}$ which passes 
through the origin. We refer to this leaf $M_0$ as the ``strong-stable''
manifold of the origin, and we denote by $W_s^{loc}$ its restriction to a 
neighborhood of zero. In contrast to $W_c^{loc}$, the local strong-stable 
manifold $W_s^{loc}$ is smooth and unique. In particular, it does not depend 
on the way in which the cut-off function $\chi$ is chosen.

\begin{theorem} \label{strongman}
Fix $k \in \intplus$, $m \ge k+2$, and let $E_c, E_s$ be as in 
Theorem~\ref{invman}. Then there exists $r_3 >0$ and a unique $C^\infty$ 
function $f : \{w_s \in E_s \,|\, \|w_s\|_m < r_3\} \to E_c$ with 
$f(0) = 0$, $Df(0) = 0$, such that the submanifold
$$
  \wsloc = \{w_s + f(w_s) ~|~ w_s \in E_s \ ,\ \| w_s \|_m < r_3\}
$$
satisfies, for any $\mu \in (\frac{k}{2},\frac{k+1}{2}]$,
\begin{equation}\label{charact}
  \wsloc = \Bigl\{w \in L^2(m) \,\Big|\, |||w|||_m < r_3\,,\ \  
  \limsup_{\tau \to +\infty} \tau^{-1} \ln\|\Phi_\tau^m w\|_m \le 
  -\mu\Bigr\}\ ,
\end{equation}
where $|||w|||_m = \max(\|P_k w\|_m,\|Q_k w\|_m)$. In particular, if 
$w_0 \in \wsloc$, there exists $T \ge 0$ such that $\Phi_\tau^m w_0 \in 
\wsloc$ for all $\tau \ge T$.
\end{theorem}

\proof Choose $r_3 > 0$ sufficiently small so that $\|\Phi_\tau^m w_0\|_m
\le r_0$ for all $\tau \ge 0$ whenever $\|w_0\|_m \le 2r_3$. Take the 
function $h(\cdot,\cdot)$ of point {\bf 2} in Theorem~\ref{invman}, and define
$$
  f(w_s) = h(0,w_s)\ , \quad \hbox{for all } w_s \in E_s \hbox{ with }
  \|w_s\|_m < r_3\ .
$$
Then $f$ is of class $C^1$, $f(0) = 0$, and the characterization 
\reff{charact} with $\mu = \mu_2$ follows immediately from 
\reff{whiskercontract} (with $w = 0$). In particular, $f$ is unique. 
Moreover, since any solution $w(\tau)$ on $\wsloc$ converges to zero
as $\tau \to +\infty$, and since $\mu_2 \in (\frac{k}{2},\frac{k+1}{2})$
was arbitrary, it is clear that \reff{charact} holds for any 
$\mu \in (\frac{k}{2},\frac{k+1}{2})$, hence for $\mu = \frac{k+1}{2}$ 
also. Finally, the smoothness of $f$ and the fact that $Df(0) = 0$ can be
proved using the integral equation satisfied by $f$. (See \cite{henry:1981},
Section~5.2.) \QED

\medskip
If $W_s^{loc}$ is as in Theorem~\ref{strongman}, we also define the global
strong-stable manifold of the origin by
\begin{eqnarray}\nonumber
  W_s &=& \bigl\{w_0 \in L^2(m) \,|\, \Phi_\tau^m(w_0) \in W_s^{loc}
   \hbox{ for some }\tau \ge 0\bigr\}\\ \label{globalss}
  &=& \bigl\{w_0 \in L^2(m) \,|\, \limsup_{\tau \to +\infty} \tau^{-1} 
   \ln\|w(\tau)\|_m \le -\mu\bigr\}\ ,
\end{eqnarray}
where $\mu \in (\frac{k}{2},\frac{k+1}{2}]$. Then $W_s$ is a smooth 
embedded submanifold of $L^2(m)$ of finite codimension. (See Henry 
\cite{henry:1981}, Sections~6.1 and 7.3.) The proof of this statement
is rather involved and uses two main ingredients. First, for any 
$\tau \ge 0$, the map $\Phi_\tau^m : L^2(m) \to L^2(m)$ is one-to-one. 
This property of the semiflow $\Phi_\tau^m$ is usually called 
{\it backwards uniqueness}. Next, since $\Phi_\tau^m$ preserves the
integral of $w$, it is clear that $W_s$ is contained in the subspace
$L^2_0(m) = \{w \in L^2(m)\,|\, \int_{\real^2}w(\xi)\d\xi = 0\}$. 
If $w(\tau)$ is any solution of \reff{SV2} in $L^2_0(m)$ and if 
$\vv(\tau)$ is the corresponding velocity field, it follows from 
Corollary~\ref{velint1} that $\vv(\tau) \in L^2(\real^2)$. 
Then, the classical energy equality
$$
   |\vv(\tau)|_2^2 = |\vv(0)|_2^2 -\int_0^\tau |\nabla \vv(s)|_2^2 \d s\ ,
   \quad \tau \ge 0\ ,
$$
shows that $|\vv(\tau)|_2$ is a strictly decreasing function of $\tau$
unless $\nabla \vv(\tau) \equiv 0$, which is equivalent to $w(\tau) 
\equiv 0$. In other words, the semiflow $\Phi_\tau^m$ is {\it strictly 
gradient} in the subspace $L^2_0(m)$. 


\section{The long-time asymptotics of solutions} \label{applications}

In this section we give three applications of the results of the 
preceding section.  In the first, we examine the long-time
asymptotics of small solutions of \reff{SV2} and show that
all such solutions with non-zero total vorticity asymptotically
approach the Oseen vortex, thereby recovering the results of
Giga and Kambe \cite{giga:1988b}.  The invariant manifold approach
yields additional information and as an example, we show that all
solutions approach the vortex in a ``universal'' way.

Prior investigations of the long-time asymptotics of the Navier-Stokes 
equations (\cite{carpio:1996}, \cite{fujigaki:2000}
\cite{miyakawa:2000}) have yielded expressions in which the terms
were proportional to inverse powers of $\sqrt{t}$.  As a second
application of the invariant manifold approach we show that
if one extends these calculations to higher order (and in principle
the invariant manifold approach allows one to extend the
asymptotics to any order) one must include terms in the asymptotics
proportional to $(\log(t))^{\alpha}/(\sqrt{t})^{\beta}$, where
$\alpha,\beta \in \intplus$. We also exhibit specific classes of solutions 
for which such logarithmic terms appear.

Finally, in a third application we extend some recent results
of Miyakawa and Schonbek \cite{miyakawa:2000} on solutions of
the Navier-Stokes equations that decay ``faster than expected''.
The invariant manifold approach  allows us both to
give a more complete characterization of the set of such
solutions and provides a natural geometrical interpretation of
the conditions in \cite{miyakawa:2000}.

\subsection{Stability of the Oseen Vortex}\label{oseen}

We begin by considering the behavior of small solutions of \reff{SV2}
in the space $L^2(m)$ with $m = 2$. Acting on $L^2(2)$, the operator
$\cL$ has a simple, isolated eigenvalue $\lambda_0 = 0$, with eigenfunction
$$
   G(\xi) \,=\, \frac{1}{4\pi}\ \e^{-|\xi|^2/4}\ , \quad \xi \in \real^2\ .
$$
(See Appendix~\ref{velocity2d}.) Let $\vv^G$ denote the corresponding
velocity field, satisfying $\rot \vv^G = G$. From the explicit expression
\reff{vgfield}, it is clear that $\xxi \perp \vv^G(\xi)$, hence
\begin{equation}\label{gident}
  \vv^G \cdot \nabla G = 0\ , \quad \xi \in \real^2\ .
\end{equation}
As a consequence, for any $\alpha \in \real$, the function $w(\xi) = 
\alpha G(\xi)$ is a stationary solution of \reff{SV2}, called 
{\it Oseen's vortex}. 

Let $E_c = \span\{G\}$, let $E_s$ be the spectral subspace of $\cL$ 
corresponding to the continuous spectrum $\sigma_c = \{ \lambda \in \complex 
~|~ \Re(\lambda) \le -\half \}$, and consider the local center manifold 
$\wcloc$ given by Corollary~\ref{locinvman} (with $k = 0$, $m = 2$). 
We claim that
\begin{equation}\label{loccent2}
  W_c^{loc} = \bigl\{ \alpha G \,|\, \alpha \in \real\ ,~|\alpha|\|G\|_2
  < r_0 \}\ .
\end{equation}
Indeed, if $|\alpha|\|G\|_2 < r_0$, then $w(\xi,\tau) = \alpha G(\xi)$ 
is a solution of \reff{SV2} such that $\|w(\tau)\|_2 < r_0$ for all 
$\tau \le 0$, hence $w(\tau) \in W_c^{loc}$ by Corollary~\ref{locinvman}. 
Since $W_c^{loc} \subset \{\alpha G + g(\alpha G)\,|\, \alpha \in \real\}$
for some $g : E_c \to E_s$, it follows that $g \equiv 0$ and that 
\reff{loccent2} holds. (Remark that, in this particular case, the local 
center manifold is {\it unique}.) Applying Corollary~\ref{locinvman} we 
conclude:

\begin{proposition} \label{oseen1} Fix $0 < \mu < \frac{1}{2}$. 
There exist positive constants $r_2$ and $C$ such that, for any initial 
data $w_0$ with $\|w_0\|_2 \le r_2$, the solution $w(\cdot,\tau)$ of 
\reff{SV2} satisfies:
\begin{equation}\label{vortexest}
  \|w(\cdot,\tau) - A G(\cdot)\|_2 \le C \e^{-\mu \tau}\ ,
  \quad \tau \ge 0\ ,
\end{equation}
where $A = \int_{\real^2} w_0(\xi) \d\xi$.
\end{proposition}

\proof If $r_2 > 0$ is sufficiently small, it follows from \reff{locattract} 
and \reff{loccent2} that \reff{vortexest} holds for some $A \in \real$. Now, 
an important property of \reff{SV2} is the {\it conservation of mass:} if 
$\alpha(\tau) = \int_{\real^2} w(\xi,\tau)\d\xi$, then
\begin{equation}\label{a-eqn}
  \dot{\alpha} = \int_{\real^2}(\cL w -\vv\cdot\nabla w)\d\xi = 
  \int_{\real^2}\nabla \cdot \bigl(\nabla w + \frac{1}{2}\xxi w -\vv w
  \bigr)\d\xi = 0\ .
\end{equation}
Since $\int_{\real^2}G(\xi)\d\xi = 1$, the conservation law \reff{a-eqn}
implies that $A = \int_{\real^2} w_0(\xi)\d\xi$ in \reff{vortexest}. \QED

\medskip
To facilitate the comparison of our results with those of \cite{giga:1988b} 
we revert to the unscaled variables $(x,t)$. Let 
$$ 
  \Omega(x,t) = \frac{1}{1+t}\,G\left(\frac{x}{\sqrt{1+t}}\right)\ , 
  \quad \uu^\Omega(x,t) = \frac{1}{\sqrt{1+t}}\,\vv^G\left(\frac{x}
  {\sqrt{1+t}}\right)\ . 
$$
Thus $\Omega$ is the solution of \reff{V2} corresponding, via the
change of variables \reff{omega-w}, to the solution $G$ of \reff{SV2}, 
and $\uu^\Omega$ is the associated velocity field. From 
Proposition~\ref{oseen1}, we obtain:

\begin{corollary} \label{oseen-gk} Fix $0 < \mu < \frac{1}{2}$. There
exists $r_2 > 0$ such that, for all initial data $\omega_0 \in L^2(2)$ 
with $\|\omega_0\|_2 \le r_2$, the solution $\omega(x,t)$ of \reff{V2}
satisfies
\begin{equation}\label{gk2}
  |\omega(\cdot,t) - A \Omega(\cdot,t)|_p \le 
  \frac{C_p}{(1+t)^{1+\mu -\frac{1}{p}}}\ ,\quad 1\le p \le 2\ ,
  \quad t \ge 0\ ,
\end{equation}
where $A = \int_{\real^2}\omega_0(x)\d x$. If $\uu(x,t)$ is the 
velocity field obtained from $\omega(x,t)$ via the Biot-Savart law
\reff{BS2}, then
\begin{equation}\label{gku}
  |\uu(\cdot,t) - A \vv^\Omega(\cdot,t)|_q \le 
  \frac{C_q}{(1+t)^{\frac{1}{2}+\mu -\frac{1}{q}}}\ , \quad 1 < q < \infty\ , 
  \quad t \ge 0\ .
\end{equation}
\end{corollary}

\proof Let $\omega(x,t)$ be the solution of \reff{V2} with $\omega(\cdot,0) 
= \omega_0$, and let $w(\xi,\tau)$ be the solution of \reff{SV2} with the 
same initial data. If $1 \le p \le 2$, then $L^2(2) \hookrightarrow
L^p(\real^2)$. Using \reff{omega-w} and \reff{vortexest}, we thus obtain 
\begin{eqnarray*}
 &&|\omega(\cdot,t) - A \Omega(\cdot,t)|_p = (1+t)^{-1+\frac{1}{p}}
  |w(\cdot,\log(1+t)) - A G(\cdot)|_p \\
 &&\hspace{1cm} \le C (1+t)^{-1+\frac{1}{p}} \|w(\cdot,\log(1+t)) - 
  A G(\cdot)\|_2 \le C (1+t)^{-1-\mu+\frac{1}{p}}\ .
\end{eqnarray*}
It then follows from Lemma~\ref{HLS2} that \reff{gku} holds for all 
$q \in (2,\infty)$. Finally, assume that $1 < q \le 2$, and fix $m \in 
(2/q,2)$. If $\tilde w(\tau) = w(\tau) -AG$ and if $\tilde \vv(\tau)$ 
denotes the corresponding velocity field, it follows from 
Proposition~\ref{velvort2} and H\"older's inequality that 
$$
   |\tilde \vv(\tau)|_q \le C |b^{m-\frac{1}{2}}\tilde \vv(\tau)|_4 
   \le C |b^m \tilde w(\tau)|_2 \le C \|\tilde w(\tau)\|_2 \le 
   C \e^{-\mu \tau}\ , \quad \tau \ge 0\ ,
$$
where $b(\xi) = (1+|\xi|^2)^{1/2}$. Using the change of variables 
\reff{u-v}, we thus obtain \reff{gku} for $1 < q \le 2$. \QED

\begin{remark} Once \reff{gk2} is known for $p \in [1,2]$, a bootstrap 
argument using the integral equation satisfied by $\omega(x,t)$ gives the 
same estimate for all $p \in [1,\infty]$ if $t \ge 1$ (see 
\cite{giga:1988b}, Proposition~5.3.) Similarly, \reff{gku} holds for 
$q \in (1,+\infty]$ if $t \ge 1$. However, as we shall see below, 
the difference $\uu(t) - A \vv^\Omega(t)$ is not in $L^1(\real^2)^2$ 
in general.
\end{remark} 

\begin{remark} In \cite{giga:1988b}, Giga and Kambe show that 
\reff{gk2} remains true if $\omega_0$ is any finite measure on 
$\real^2$ satisfying $\int_{\real^2}(1+|x|^2)|\omega_0|(\D x) < \infty$
and $\int_{\real^2}|\omega_0|(\D x) \ll 1$. Note that, in this
case, the solution $\omega(x,t)$ satisfies $\|\omega(\cdot,t)\|_2 \le
r_2$ for all $t \ge 1$, so that Corollary~\ref{oseen-gk} applies 
to $\omega(x,t)$ restricted to $t \ge 1$. In \cite{carpio:1994}, Carpio 
also shows that solutions of \reff{V2} with finite measures as initial 
data satisfy estimate \reff{gk2}. Instead of requiring $\omega_0$ to be
small in some norm, she assumes that the Cauchy problem for the vorticity
equation \reff{V2} with initial data $A\delta$ has a unique solution, 
where $A = \int_{\real^2}\omega_0(\D x)$ and $\delta$ is the Dirac mass
at the origin. However, to the best of our knowledge, the uniqueness of 
solutions for this problem has only been established if $|A|$ is 
sufficiently small. 
\end{remark}

We next note that not only does the invariant manifold approach yield
the stability of the Oseen vortex but it also allows us to systematically
compute the way in which the vortex is approached.  Somewhat surprisingly,
the approach to the vortex solution has a universal form if the norm
of the initial data is not too large.  More precisely, we prove:

\begin{theorem} \label{oseen2} Fix $\half < \mu < 1$.  There exist
$r_2 > 0$ and $C > 0$ such that, for all initial data $w_0 \in 
L^2(3)$ with $\|w_0\|_3 \le r_2$, the solution $w(\cdot,\tau)$ 
of \reff{SV2} satisfies
\begin{equation}\label{approach-vortex}
  \|w(\xi,\tau) - AG(\xi) + \half(B_1 \xi_1 + B_2 \xi_2) 
  G(\xi) \e^{-\tau/2}\|_3 \le C \e^{-\mu \tau}\ ,\quad \tau \ge 0\ ,
\end{equation}
where $A = \int w_0(\xi) \d\xi$, $B_1 = \int \xi_1 w_0(\xi) 
\d\xi$, and $B_2 = \int \xi_2 w_0(\xi) \d\xi$.
\end{theorem}

\proof Acting on $L^2(3)$, the operator $\cL$ has, in addition to the
simple eigenvalue $\lambda_0 = 0$, a double eigenvalue $\lambda_1 = 
-\frac{1}{2}$, with eigenfunctions $F_1(\xi) = -\frac{\xi_1}{2} G(\xi)$ 
and $F_2(\xi) = -\frac{\xi_2}{2} G(\xi)$. (See Appendix~\ref{velocity2d}.)
Let $E_c = \span\{G,F_1,F_2\}$, let $E_s$ be the spectral subspace of 
$\cL$ corresponding to the continuous spectrum $\sigma_c = \{\lambda \in 
\complex~|~ \Re(\lambda) \le 1\}$, and consider the local invariant 
manifold $W_c^{loc}$ given by Corollary~\ref{locinvman} (with $k = 1$, 
$m = 3$). 

By construction, a point $w \in W_c^{loc}$ can be written, in a unique way, 
as $w = \alpha G + \beta_1 F_1 +\beta_2 F_2 + g(\alpha,\bbeta)$, where 
$\bbeta = (\beta_1,\beta_2)$, and $g(\alpha, \bbeta) \in E_s$. The 
coefficients $\alpha, \bbeta$ are given by the formulas
\begin{equation}\label{alphabetcoeff}
  \alpha = \int_{\real^2}w(\xi)\d\xi\ ,\quad 
  \beta_i = -\int_{\real^2} \xi_i w(\xi)\d\xi\ ,\quad i = 1,2\ .
\end{equation}
Now, if $w(\cdot,\tau)$ evolves according to 
\reff{SV2}, we know from \reff{a-eqn} that $\dot{\alpha} = 0$, {\it i.e.}
$\alpha$ does not change with time. In an analogous way, we find
\begin{equation}\label{b-eqn}
  \dot{\beta_i} = -\int_{\real^2}\xi_i(\cL w -\vv\cdot\nabla w)
  \d\xi = -\frac{1}{2}\beta_i\ , \quad i = 1,2\ .
\end{equation}
Indeed, since $\nabla\cdot\vv = 0$ and $\rot \vv = \partial_1 v_2 
-\partial_2 v_1 = w$, we have the identities
\begin{eqnarray*}
  \xi_1 \cL w + \frac{1}{2}\xi_1 w &=& \partial_1\Bigl(\xi_1 \partial_1 w
   + \frac{1}{2}\xi_1^2 w - w\Bigr) + \partial_2\Bigl(\xi_1 \partial_2 w
   + \frac{1}{2}\xi_1 \xi_2 w\Bigr)\ ,\\
  \xi_1 \vv\cdot\nabla w &=& \partial_1\Bigl(\xi_1 v_1 w -v_1 v_2\Bigr)
  + \partial_2\Bigl(\xi_1 v_2 w + \frac{1}{2}(v_1^2-v_2^2)\Bigr)\ ,
\end{eqnarray*}
which prove \reff{b-eqn} for $i = 1$. The case $i = 2$ is similar.

\begin{remark} In the unscaled variables, \reff{b-eqn} means that 
equation \reff{V2} conserves the first moments of the vorticity
$\omega(x,t)$. 
\end{remark}

Thus, remarkably, the semiflow induced by \reff{SV2} on the 
three-dimensional invariant manifold is described by the {\it linear} 
equations \reff{a-eqn}, \reff{b-eqn}!
To complete the investigation of how the solutions evolve on the center
manifold, we must estimate the nonlinear term $g(\alpha,\bbeta)$.
Given $\frac{1}{2} < \mu < 1$, we claim that there exist $\epsilon > 0$
and $C_g > 0$ such that
\begin{equation}\label{g-bound}
  \|g(\alpha,\bbeta)\|_3 \le C_g |\bbeta |^{2\mu}\ ,
\end{equation}
for all $(\alpha,\bbeta) \in \real^3$ with $|\alpha| + |\bbeta| \le 
\epsilon$. Indeed, we know from Theorem~\ref{invman} (see also 
Remark~\ref{smoothness}) that $g : \real^3 \to E_s$ is of class 
$C^{2\mu}$, and that $g(0,0) = 0$, $Dg(0,0) = 0$. 
Moreover, since $\alpha G $ is a fixed point
for \reff{SV2}, in must lie in $ W_c^{loc}$ when $\alpha$ is sufficiently
small, so we have $g(\alpha,0) = 0$ for such values of $\alpha$.  
Thus, to prove \reff{g-bound}, it is sufficient to show that 
$\partial_\beta g(\alpha,0) = 0$ when $|\alpha|$ is small. 

To prove this, we linearize the map $\phirm$ in \reff{phidef} about the 
fixed point $\alpha^* G$ obtaining
\begin{equation}
  D_{\alpha^* G} \phirm = \Lambda +D_{\alpha^* G} \cR \ .
\end{equation}
Since $\cR$ is smooth, and $\cR = 0$, $D_0 \cR =0$, $D_{\alpha^* G} \cR$ 
can be made arbitrarily small in norm by taking $\alpha^*$ sufficiently small.
In particular, since we know the spectrum of $\Lambda$ explicitly, we can 
choose $\alpha^*$ sufficiently small so that the spectrum of $D_{\alpha^* G} 
\phirm$ is ``close'' to the spectrum of $\Lambda$. More precisely, if 
$0 < \delta < \min(\frac{1}{4},1-\mu)$, there exists $\alpha_0 > 0$ such 
that if $0 \le \alpha^* \le \alpha_0$, $D_{\alpha^* G} \phirm$ has a 
three-dimensional spectral subspace $\Eca$ with the spectrum of 
$D_{\alpha^* G} \phirm|_{\Eca}$ consisting of eigenvalues with absolute value 
greater than $\exp(-(\half +\delta))$ (these are the perturbations of the 
eigenvalues $1$ and $\e^{-1/2}$ of $\Lambda$) and the remainder of the 
spectrum of $D_{\alpha^* G} \phirm$ contained in a disk in the complex plane 
of radius $\exp(-(1-\delta))$. Applying the invariant manifold theorem of 
\cite{chen:1997} to the semiflow $\phit$ about the fixed point 
$\alpha^* G$ (instead of the origin), we find just as in Theorem~\ref{invman} 
that $\phit$ has a three-dimensional invariant manifold
$\wca$, which in is tangent to $\Eca$ at  $\alpha^* G$.  Furthermore, 
as we noted in Remark~\ref{smoothness}, once the cutoff function $\chi_{r_0}$ 
in \reff{SV2c} is fixed, the invariant manifold constructed in 
\cite{chen:1997} is unique. Therefore, it follows from the characterizations
\reff{char1}, \reff{char2} that $\wca$ actually coincides with $\wcloc$
in a neighborhood of $\alpha^* G$. We complete the argument by showing 
that $\Eca = \span\{G,F_1,F_2\}$. This implies that $\partial_\beta 
g(\alpha^*,0) =0$ for $\alpha^*$ sufficiently small, thus proving 
\reff{g-bound}.

To compute $\Eca$, note that if we linearize \reff{SV2} about
the fixed point $w=\alpha^* G$, we obtain
\begin{equation}
  \partial_{\tau} w = \cL^{\alpha^*} w =
  \Delta_\xi w + \half (\xxi \cdot \nabla_\xi)w +w
  - \alpha^*(\vv^G \cdot \nabla_\xi w + \vv \cdot \nabla_{\xi} G)\ ,
\end{equation}
where $\vv^G$ is the velocity field associated with the Oseen vortex $G$ 
and $\vv$ is the velocity field constructed from $w$ via \reff{SBS2}. 
Differentiating identity \reff{gident} with respect to $\xi_j$ ($j = 1,2$), 
we obtain
\begin{eqnarray} \nonumber
  0 &=& \partial_j (\vv^G \cdot \nabla_{\xi} G)
  = ( (\partial_j \vv^G) \cdot \nabla_{\xi} )G +
  (\vv^G \cdot \nabla_{\xi})(\partial_j G) \\ \label{fident}
  &=& (\vv^{F_j} \cdot \nabla_{\xi})G + \vv^G \cdot \nabla_{\xi} F_j \ .
\end{eqnarray}
Combining \reff{gident} and \reff{fident} we see immediately that
\begin{eqnarray*}
  \cL^{\alpha^*} G &=& \cL G - \alpha^*(\vv^G \cdot \nabla_{\xi} G
  + \vv^G \cdot \nabla_{\xi} G) = \cL G =0 \\
  \cL^{\alpha^*} F_j &=& \cL F_j - \alpha^*(\vv^G \cdot \nabla_{\xi} F_j
  + \vv^{F_j} \cdot \nabla_{\xi} G) = \cL F_j = -\half F_j \ , \quad
  j = 1,2\ .
\end{eqnarray*}
Exponentiating these results we find that $D_{\alpha^* G} \phirm
= \exp(\cL^{\alpha^*})$ has eigenvalues $1$ and $\e^{-\frac{1}{2}}$
with eigenspaces $\{G\}$ and $\{F_1,F_2\}$ respectively. This conludes
the proof of \reff{g-bound}. 

Assume now that $w(\cdot,\tau)$ is a solution of \reff{SV2} on $W_c^{loc}$,
and let $\alpha(\tau) = A$, $\beta_1(\tau) = B_1 \e^{-\tau/2}$, $\beta_2(\tau) 
= B_2 \e^{-\tau/2}$, where $A,B_1,B_2$ are as in Theorem~\ref{oseen2}. 
Then
$$
   w(\cdot,\tau) = \alpha(\tau) G + \beta_1(\tau) F_1 + \beta_2(\tau) F_2
   + g(\alpha(\tau),\bbeta(\tau))\ ,
$$
hence \reff{approach-vortex} follows directly from \reff{g-bound}. 

On the other hand, if $\|w_0\|_3 \le r_2$ and if $w(\tau)$ is the solution
of \reff{SV2} with initial data $w_0$, Corollary~\ref{locinvman} shows that
there exists a solution $\tilde{w}(\tau)$ on $W_c^{loc}$ such that 
$\|w(\tau)-\tilde{w}(\tau)\|_3 \le C\e^{-\mu\tau}$ for all $\tau \ge 0$. 
In view of the previous result, this means that \reff{approach-vortex}
holds for some $A, B_1, B_2 \in \real$. But equations \reff{a-eqn} and 
\reff{b-eqn}, which hold for any solution of \reff{SV2}, imply that 
$A = \int w_0(\xi) \d\xi$ and $B_i = \int \xi_i w_0(\xi) \d\xi$, 
$i = 1,2$.
This concludes the proof of Theorem~\ref{oseen2}. \QED

\begin{remark} Note that although the preceding computation
of the asymptotics of the solutions near the Oseen vortex
applies only if $\alpha^*$ is small, $0$ and $-\half$ are
eigenvalues of $\cL^{\alpha^*}$, with eigenfuctions $G$ and $F_1,F_2$,
for all values of $\alpha^*$.
\end{remark}

\medskip
One can also rewrite the result of Theorem~\ref{oseen2} in terms
of the unscaled variables as we did in Corollary~\ref{oseen-gk}.
Denote
\begin{eqnarray*}
  \omega_\app(x,t) &=& \frac{A}{1+t} \,G\Bigl(\frac{x}{\sqrt{1+t}}\Bigr) 
   + \sum_{i=1}^2 \frac{B_i}{(1+t)^{3/2}} \,F_i\Bigl(\frac{x}{\sqrt{1+t}}
   \Bigr)\ , \\ 
  \uu_\app(x,t) &=& \frac{A}{\sqrt{1+t}} \,\vv^G\Bigl(\frac{x}{\sqrt{1+t}}
   \Bigr) + \sum_{i=1}^2 \frac{B_i}{1+t} \,\vv^{F_i}\Bigl(\frac{x}{\sqrt{1+t}}
   \Bigr)\ . \\ 
\end{eqnarray*}

\begin{corollary}\label{oseen-asym}
Fix $\half < \mu < 1$.  There exists $r_2 > 0$ such that, for all 
initial data $\omega_0 \in L^2(3)$ with $\|w_0\|_3 \le r_2$, the
solution $\omega(\cdot,\tau)$ of \reff{SV2} satisfies
$$
  |\omega(\cdot,t) - \omega_\app(\cdot,t)|_p \le 
  \frac{C_p}{(1+t)^{1+\mu -\frac{1}{p}}}\ , \quad 
  1 \le p \le 2\ , \quad t \ge 0\ ,
$$
where $A = \int \omega_0(x)\d x$ and $B_i = \int x_i \omega_0(x) 
\d x$, $i = 1,2$. If $\uu(x,t)$ is the velocity field obtained from 
$\omega(x,t)$ via the Biot-Savart law \reff{BS2}, then
$$
  |\uu(\cdot,t) - \uu_\app(\cdot,t)|_q \le 
  \frac{C_q}{(1+t)^{\frac{1}{2}+\mu -\frac{1}{q}}}\ , \quad 
  1 \le q < \infty\ , \quad t \ge 0\ .
$$
\end{corollary}

\subsection{Secular terms in the asymptotics}\label{logarithm}

We next show that, in contrast to prior investigations of the
long-time asymptotics of solutions of \reff{V2} which yielded
expansions in inverse powers of the time, one will in general
encounter terms in the asymptotics with contain factors of $\log(t)$.
These terms arise from resonances between the eigenvalues of the linear 
operator $\cL$, and as we will see, computing them is straightforward 
using ideas from the theory of dynamical systems.

In the previous subsection we saw that, if $w_0 \in L^2(2)$ is sufficiently
small, the solution $w(\xi,\tau)$ of \reff{SV2} with initial data $w_0$ 
approaches the Oseen vortex $AG(\xi)$ as time goes to infinity, where
$A = \int w_0(\xi) \d\xi$. We now examine in more detail what happens if 
$A = 0$, namely $w_0 \in L^2_0(m)$. In this invariant
subspace, we know from Theorem~\ref{weightedL2} that all solutions 
converge to zero as $\tau \to +\infty$, without any restriction on the
size of the initial data. Assuming that $w_0 \in L^2(4)$, we will 
compute the long-time asymptotics of the solution up to terms of order 
$\cO(\e^{-\mu \tau})$, for any $\mu$ in the interval $(1,\frac{3}{2})$.

We first note that the operator $\cL$ acting on $L^2_0(4)$ has, in addition 
to $\lambda_1 = -\frac{1}{2}$, a triple eigenvalue $\lambda_1 = -1$, with 
eigenfunctions
\begin{equation}\label{Heigen}
  H_1(\xi) = \frac{1}{4}(|\xi|^2 - 4)G(\xi)\ , \quad
  H_2(\xi) = \frac{1}{4}(\xi_1^2 - \xi_2^2)G(\xi)\ , \quad
  H_3(\xi) = \frac{1}{4}\xi_1 \xi_2 G(\xi)\ .  
\end{equation}
Of course, these eigenfunctions are not unique, but the choices above are 
convenient ones for computation.

Let $E_c = \span\{F_1,F_2,H_1,H_2,H_3\}$, and let $E_s \subset L^2_0(4)$ be 
the spectral subspace of $\cL$ corresponding to the continuous spectrum 
$\sigma_c = \{\lambda \in \complex ~|~ \Re(\lambda) \le -\frac{3}{2}\}$.
Any function $w \in L^2_0(4)$ can be written as
\begin{equation}\label{w3decomp}
  w(\xi) = \sum_{i=1}^2 \beta_i F_i(\xi) + \sum_{j=1}^3 \gamma_j
  H_j(\xi) + R(\xi)\ ,
\end{equation}
where $R \in E_s$. The velocity field $\vv(\xi)$ associated to $w$ has a 
similar decomposition:
\begin{equation}\label{v3decomp}
  \vv(\xi) = \sum_{i=1}^2 \beta_i \vv^{F_i}(\xi) + \sum_{j=1}^3 \gamma_j 
  \vv^{H_j}(\xi) + \vv^R(\xi)\ ,
\end{equation}
see Appendix~\ref{velocity2d}. The coefficients $\beta_i$ are given by 
\reff{alphabetcoeff}, and the corresponding formulas for $\gamma_j$ read:
\begin{equation}\label{gammacoeff}
  \gamma_j = \int_{\real^2} p_j(\xi) w(\xi)\d\xi\ , \quad j = 1,2,3\ ,
\end{equation}
where $p_1(\xi) = \frac{1}{4}(|\xi|^2 - 4)$, $p_2(\xi) = \frac{1}{4}
(\xi_1^2-\xi_2^2)$ and $p_3(\xi) = \xi_1 \xi_2$. 

Assume now that $w(\xi,\tau)$ is a solution of \reff{SV2} in $L^2_0(4)$, 
and consider the evolution equations for the coefficients $\gamma_j$ 
and the remainder $R$ in \reff{w3decomp}. Proceeding as in \reff{b-eqn}, 
we find 
\begin{equation}\label{gam-eqn}
   \dot\gamma_j = -\gamma_j -\int_{\real^2} p_j(\xi)(\vv\cdot\nabla)
   w \d\xi\ , \quad j = 1,2,3\ .
\end{equation}
The following elementary result will be useful:

\begin{lemma}\label{proj_lem} Assume that $w \in L^2_0(2) \cap H^1(2)$, 
and let $\vv$ be the velocity field obtained from $w$ via the Biot-Savart 
law \reff{SBS2}. Then, for any quadratic polynomial $p(\xi)$, 
\begin{equation}\label{pquadident}
  \int_{\real^2} p(\xi) (\vv \cdot \nabla) w \d\xi
  = \int_{\real^2} \bigl(v_1 v_2 (\partial_1^2 p - \partial_2^2 p) 
  - (v_1^2 - v_2^2) \partial_1 \partial_2 p \bigr)\d\xi\ .
\end{equation}
\end{lemma}

\proof Since $\nabla\cdot\vv = 0$ and $w = \partial_1 v_2 - \partial_2 v_1$,
we have the identity
$$
  p(\vv\cdot\nabla)w = v_1 v_2 (\partial_1^2 p - \partial_2^2 p) - 
  (v_1^2 - v_2^2) \partial_1 \partial_2 p + \partial_1 E_1 + 
  \partial_2 E_2\ , 
$$
where $E_1 = pv_1 w - v_1 v_2 \partial_1 p + \frac{1}{2}(v_1^2-v_2^2)
\partial_2 p$ and $E_2 = pv_2 w + v_1 v_2 \partial_2 p + \frac{1}{2}
(v_1^2-v_2^2)\partial_1 p$. Applying Proposition~\ref{velvort2} with 
$m = 2$, we see from \reff{velv2} that $(1+|\xi|^2)\vv \in 
L^\infty(\real^2)^2$, hence $E_i \in L^2(\real^2)$ for $i = 1,2$. In 
addition, since $\nabla \vv \in L^2(\real^2)^4$ by Lemma~\ref{HLS2}, 
we have $\partial_i E_i \in L^1(\real^2)$ for $i = 1,2$, and 
\reff{pquadident} follows. \QED

\medskip
It follows in particular from Lemma~\ref{proj_lem} that $\int_{\real^2} 
p_1(\xi)(\vv \cdot\nabla)w \d\xi \equiv 0$. Thus, surprisingly enough, 
the equation for $\gamma_1$ is {\it linear:} $\dot \gamma_1 = -\gamma_1$.
In contrast, the equations for $\gamma_2$, $\gamma_3$, and $R$ contain 
nonlinear terms. For the purposes of this subsection, it will be useful 
to write out separately the terms that are quadratic in $\bbeta = 
(\beta_1,\beta_2)$. A direct calculation shows that 
\begin{equation}\label{betaquad}
  (\beta_1 \vv^{F_1} + \beta_2 \vv^{F_2})\cdot \nabla (\beta_1 F_1 +
  \beta_2 F_2) = \bigl( (\beta_1^2-\beta_2^2)\xi_1 \xi_2 - \beta_1 \beta_2
  (\xi_1^2-\xi_2^2)\bigr)\Phi(\xi)\ ,
\end{equation}
where $\Phi(\xi) = (8\pi^2 |\xi|^4)^{-1}\e^{-|\xi|^2/4}(\e^{-|\xi|^2/4}
-1 + |\xi|^2/4)$. Thus, defining
\begin{equation}\label{kappavalue}
  \kappa = \int_{\real^2} \xi_1^2 \xi_2^2 \Phi(\xi)\d\xi = 
  \frac{1}{4}  \int_{\real^2} (\xi_1^2 -\xi_2^2)^2 \Phi(\xi)\d\xi = 
  \frac{1}{32\pi}\ ,
\end{equation}
we see that the quadratic terms in the equations for $\gamma_2$, $\gamma_3$ 
are respectively $\kappa \beta_1 \beta_2$ and  $-\kappa (\beta_1^2 
-\beta_2^2)$. Therefore, the equations for $\bbeta,\ggamma,R$
(where $\ggamma = (\gamma_1,\gamma_2,\gamma_2)$) have the
following form:
\begin{eqnarray} \nonumber
  \dot{\beta}_1 &=& -\half \beta_1\ , \quad 
  \dot{\beta}_2 = -\half \beta_2\ , \quad
  \dot{\gamma}_1 = - \gamma_1\ ,\\ \label{log1}
  \dot{\gamma}_2 &=& - \gamma_2 +\kappa \beta_1 \beta_2
    + f_2(\bbeta,\ggamma,R)\ ,\\ \nonumber
  \dot{\gamma}_3 &=& - \gamma_3 - \kappa(\beta_1^2 - \beta_2^2)
    + f_3(\bbeta,\ggamma,R)\ ,\phantom{\frac{1}{2}} \\ \nonumber
  R_\tau &=& \cL R + \beta_1 \beta_2 \Psi_2 - (\beta_1^2-\beta_2^2)\Psi_3 
    + \nabla\cdot \FF(\bbeta,\ggamma,R) - \vv^R \cdot \nabla R \ ,
\end{eqnarray}
where $\Psi_2 = (\xi_1^2-\xi_2^2)\Phi -\kappa H_2$ and $\Psi_3 = \xi_1 
\xi_2 \Phi - \kappa H_3$. Moreover, since we have written out explicitly
the quadratic terms in $\bbeta$, the remainder terms $f_j$, $\FF $ in 
\reff{log1} satisfy the estimates
\begin{equation}\label{ffbounds}
  |f_j| + \|\FF\|_4 \le C (|\bbeta|+|\ggamma|)(|\ggamma|+\|R\|_4)
    + C\|R\|_4^2 \ .
\end{equation}

To compute the long-time asymptotics of the system \reff{log1}, a 
natural idea is first to study the behavior of the solutions on 
the five-dimensional invariant manifold $W_c^{loc}$ tangent to $E_c$ at 
the origin, and then to use the fact that all solutions of \reff{log1} 
approach $W_c^{loc}$ faster than $\cO(\e^{-\mu\tau})$ for any 
$\mu < \frac{3}{2}$, see Corollary~\ref{locinvman}. This is the approach 
we follow, but we begin with a simplifying observation.  Note that if one 
drops the terms $\nabla \cdot \FF-\vv^R \cdot \nabla R$ from the last 
equation in \reff{log1}, the terms $\beta_1 \beta_2 \Psi_2 - 
(\beta_1^2-\beta_2^2)\Psi_3$ act as a simple ``forcing term'' in the 
equation for $R$. As a result, one can find an explicit invariant manifold
for the simplified equations, namely:
\begin{equation}\label{expliman}
   R= g(\bbeta) = -\beta_1 \beta_2 (\cL+1)^{-1}\Psi_2 + 
   (\beta_1^2-\beta_2^2)(\cL+1)^{-1}\Psi_3\ .
\end{equation}
(Remark that $\Psi_2, \Psi_3 \in E_s$ and that the restriction of 
$\cL+1$ to $E_s$ is invertible.) As we will see, this manifold is a good 
approximation to an invariant manifold for the full equations \reff{log1}. 
We thus introduce the new variable
\begin{equation}\label{rhodef}
  \rho = R - g(\bbeta)\ ,
\end{equation}
and find that \reff{log1} can be rewritten as
\begin{eqnarray} \nonumber
  \dot{\beta}_1 &=& -\half \beta_1\ , \quad 
  \dot{\beta}_2 = -\half \beta_2\ , \quad
  \dot{\gamma}_1 = - \gamma_1\ ,\\ \label{log2}
  \dot{\gamma}_2 &=& - \gamma_2 +\kappa \beta_1 \beta_2
    + \tf_2(\bbeta,\ggamma,\rho)\ ,\\ \nonumber
  \dot{\gamma}_3 &=& - \gamma_3 - \kappa(\beta_1^2 - \beta_2^2)
    + \tf_3(\bbeta,\ggamma,\rho)\ ,\phantom{\frac{1}{2}} \\ \nonumber
  \rho_\tau &=& \cL \rho + \nabla\cdot \tFF(\bbeta,\ggamma,\rho) - 
\vv^\rho \cdot \nabla \rho \ ,
\end{eqnarray}
where $\tf_j$ and $\tFF$ obey the same estimates as $f_j$ and $\FF$, but 
with $\|\rho\|_4 + |\bbeta|^2$ substituted for $\|R\|_4$. Note further 
that, if $R \in E_s$, then $\rho \in E_s$ as well. 

If we now apply the results of \cite{chen:1997} to the system \reff{log2}, 
we find just as in Theorem~\ref{invman} that \reff{log2} has a 
five-dimensional invariant manifold tangent at the origin to $E_s$. The next 
proposition provides an estimate of the behavior of this manifold. Let 
$X$ be the Banach space $\real^2 \times \real^3 \times E_s$, equipped with 
the norm $\|(\bbeta,\ggamma,R)\|_X  = |\bbeta| + |\ggamma| + \|R\|_4$,
and set $E_c=\real^2 \times \real^3$, with coordinates $(\bbeta,\ggamma)$.

\begin{proposition}\label{invman2} Fix $\mu \in (1,\frac{3}{2})$ and
$\delta \in (0,\frac{3}{2}-\mu)$. There exists a $C^1$ map $\cG : E_c 
\to E_s$ satisfying $\cG(0) = 0$, $D\cG(0) = 0$, such that in a 
neighborhood of the origin the graph of $\cG$ is left invariant by 
the semiflow defined by \reff{log2}. Moreover, there exists $r_4 > 0$ and 
$C>0$ such that, for any solution $(\bbeta(\tau),\ggamma(\tau),\rho(\tau))$ 
of \reff{log2} with initial data $(\bbeta_0,\ggamma_0,\rho_0)$ satisfying 
$\|(\bbeta_0,\ggamma_0,\rho_0)\|_X \le r_4$, there exists a solution
$(\bbeta^c(\tau),\gamma^c(\tau),\cG(\bbeta^c(\tau),\gamma^c(\tau))$
on the invariant manifold with
\begin{equation}\label{bgrdec}
   |\bbeta(\tau) - \bbeta^c(\tau)| + |\ggamma(\tau) - \ggamma^c(\tau)| +
   \|\rho(\tau) - \cG(\bbeta^c(\tau),\ggamma^c(\tau))\|_4 \le 
   C \e^{-\mu \tau}\ .
\end{equation}
Finally, there exists $C_{\delta} > 0$ such that 
\begin{equation}\label{cGest}
  \|\cG(\bbeta,\ggamma) \|_4 \le C_{\delta} (|\bbeta|^{3-\delta}
  + |\ggamma|^{\frac{3}{2}-\delta} ) \ .
\end{equation}
\end{proposition}

With the exception of the last estimate, the proof of this
statement is almost identical to the proof of Theorem~\ref{invman} 
and Corollary~\ref{locinvman}. We prove \reff{cGest} in Appendix~\ref{Gapp} 
and concentrate here on showing how it can be used to determine the 
asymptotics of \reff{log2}.

We begin by studying the behavior of the solutions on the invariant 
manifold. Such solutions satisfy the system of ordinary differential 
equations:
\begin{eqnarray}
 \nonumber
  \dot{\beta}_1 &=& -\half \beta_1\ , \quad 
  \dot{\beta}_2 = -\half \beta_2\ , \quad
  \dot{\gamma}_1 = - \gamma_1\ ,\\ \label{log3}
  \dot{\gamma}_2 &=& - \gamma_2 +\kappa \beta_1 \beta_2
    + \tf_2(\bbeta,\ggamma ,\cG(\bbeta,\ggamma) )\\ \nonumber
  \dot{\gamma}_3 &=& - \gamma_3 - \kappa(\beta_1^2 - \beta_2^2)
    + \tf_3(\bbeta,\ggamma,\cG(\bbeta,\ggamma) ) \ .
\end{eqnarray}
There are various ways to analyze the asymptotics of solutions
of \reff{log3}, but perhaps the simplest is to make the change
of variables:
\begin{eqnarray}\label{normal-change}
 &&\Gamma_1 = \gamma_1\ ,\quad \Gamma_2 =
   \gamma_2 + \kappa \beta_1 \beta_2 |\log|\beta_1\beta_2||\ ,\\ \nonumber 
 &&\Gamma_3 = \gamma_3 - \kappa(\beta_1^2 |\log(\beta_1^2)| 
   - \beta_2^2 |\log(\beta_2^2)|)\ .
\end{eqnarray}
From a dynamical systems point of view, this is just a normal
form transformation that eliminates the resonant terms
in \reff{log3}.  Rewriting \reff{log3} in term of $\bbeta$
and $\GGamma=(\Gamma_1,\Gamma_2,\Gamma_3)$, we find
\begin{eqnarray}\nonumber
  \dot{\beta}_1 &=& -\half \beta_1\ , \quad 
  \dot{\beta}_2 = -\half \beta_2\ , \quad
  \dot{\Gamma}_1 = - \Gamma_1\ ,\\ \label{log4}
  \dot{\Gamma}_2 &=& - \Gamma_2 + \tg_2(\bbeta,\GGamma)\ , \quad
  \dot{\Gamma}_3 \ =\  - \Gamma_3 + \tg_3(\bbeta,\GGamma) \ ,
\end{eqnarray}
where $\tg_j(\bbeta,\GGamma)$ is just $\tf_j(\bbeta,\ggamma,\cG(\bbeta,
\ggamma))$ rewritten in terms of $(\bbeta, \GGamma)$ instead of $(\bbeta,
\ggamma)$. 

\begin{remark} The change of variables $(\bbeta,\ggamma) \to
(\bbeta,\GGamma)$ is not Lipshitz, so it might seem as if the standard
existence and uniqueness theory for solutions does not apply to \reff{log4}.
However, due to the very simple form of the first equations in \reff{log4},
we can first solve the equations for $\beta_1$ and $\beta_2$ explicitly 
and then insert these expressions into the equations for $\GGamma$. 
The equations for $\GGamma$ are then non-autonomous, but Lipshitz 
in $\Gamma$ and hence standard theorems imply that solutions of 
\reff{log4} exist and are unique.
\end{remark}

It is clear from \reff{log4} that $\beta_i(\tau) = b_i \e^{-\tau/2}$ 
for $i = 1,2$ and $\Gamma_1(\tau) = c_1 \e^{-\tau}$, where $b_i = \beta_i(0)$
and $c_1 = \Gamma_1(0)$. To determine the long-time behavior of $\Gamma_2$ 
and $\Gamma_3$, we note that $\e^{\tau} \Gamma_j(\tau) = \Gamma_j(0) + 
\int_0^\tau \e^s \tg_j(\bbeta(s),\GGamma(s))\d s$. Thus, defining
\begin{equation}\label{cjdef}
  c_j = \Gamma_j(0) + \int_0^{\infty} \e^s \tg_j(\bbeta(s),\GGamma(s))\d s\ ,
  \quad j = 2,3\ ,
\end{equation}
and using the estimates on $\tg_j$ which come from \reff{ffbounds}, we find
that for any $\mu \in (1,3/2)$ there exists $C > 0$ such that $|\Gamma_j(\tau) 
-\e^{-\tau} c_j| \le C \e^{-\mu \tau}$, for $j = 2,3$. Inverting the change 
of coordinates \reff{normal-change}, we immediately find:

\begin{lemma}\label{bgest} Fix $\mu \in (1,3/2)$. There exist $r_5 >0$ 
and $C > 0$ such that, for any solution of \reff{log3} with initial data
satisfying $|\bbeta(0)|+|\ggamma(0)| \le r_5$, there exist constants 
$b_1,b_2,c_1',c_2'$, and $c_3'$ such that, for all $\tau \ge 0$, 
\begin{eqnarray}\nonumber
  &&\beta_1(\tau) = b_1 \e^{-\tau/2}\ , \ \ 
  \beta_2(\tau)  = b_2 \e^{-\tau/2}\ , \ \ 
  \gamma_1(\tau) = c_1 \e^{-\tau}\ ,\\ \label{bgesteq}
  && |\gamma_2(\tau) - c_2 \e^{-\tau} + \kappa \tau \e^{-\tau} b_1 b_2 
  |\log|(b_1 b_2)|| \le  C \e^{-\mu \tau}\ , \\ \nonumber 
  && |\gamma_3(\tau) - c_3 \e^{-\tau} - \kappa \tau\e^{-\tau} (b_1^2 
  |\log(b_1^2)| -b_2^2 |\log(b_2^2)|)| \le C \e^{-\mu \tau}\ .
\end{eqnarray}
\end{lemma}

\begin{remark} The constants $c_j'$ ($j = 1,2,3$) in \reff{bgesteq} are 
related to $c_j$ by
$$
   c_1' = c_1\ ,\quad
   c_2' = c_2 -\kappa b_1 b_2 |\log|b_1 b_2||\ , \quad 
   c_3' = c_3 +\kappa \bigl(b_1^2 |\log(b_1^2)| -b_2^2 |\log(b_2^2)|\bigr)\ .
$$
\end{remark}

We now return to the full system of equations \reff{log2}.
We know that any solution $(\bbeta(\tau),\ggamma(\tau),\rho(\tau))$ 
of \reff{log2} with initial data in a sufficiently small neighborhood 
of the origin approaches a solution $(\bbeta^c(\tau),\ggamma^c(\tau),
\cG(\bbeta^c(\tau),\ggamma^c(\tau))$ on the invariant manifold with a rate 
$\cO(\e^{-\mu \tau})$, see \reff{bgrdec}. On the other hand, from 
Lemma~\ref{bgest}, we know that there exist constants $b_i, c_j'$ such that 
estimates \reff{bgesteq} hold with $\bbeta^c,\ggamma^c$ substituted for 
$\bbeta,\ggamma$. Finally, combining \reff{bgesteq} with \reff{cGest}, we 
obtain the bound $\|\cG(\bbeta^c(\tau),\gamma^c(\tau))\|_4 \le C 
\e^{-\mu \tau}$. Thus, the following estimate holds for all solutions 
of \reff{log2} in a neighborhood of the origin:
\begin{eqnarray}\nonumber
  &&|\beta_1(\tau) - b_1 \e^{-\tau/2}| +  
    |\beta_2(\tau) - b_2 \e^{-\tau/2}| +  
    |\gamma_1(\tau) - c_1' \e^{-\tau}| \\ \label{the-estim} 
  && +\ |\gamma_2(\tau) - c_2' \e^{-\tau} + \kappa \tau \e^{-\tau} b_1 b_2 
    |\log|(b_1 b_2)|| + \|\rho(\tau)\|_4 \\ \nonumber
  && +\ |\gamma_3(\tau) - c_3' \e^{-\tau} - \kappa \tau\e^{-\tau} (b_1^2 
  |\log(b_1^2)| -b_2^2 |\log(b_2^2)|)| \le C \e^{-\mu \tau}\ .
\end{eqnarray}

Finally, we undo the change of variables \reff{rhodef}. The result is:

\begin{theorem}\label{asymp2} Fix $1 < \mu < \frac{3}{2}$. If $w(\xi,\tau)$ 
is any solution of \reff{SV2} in $L^2_0(4)$, then there exist constants $b_i$, 
$c_j'$, and $C$ such that $\|w(\cdot,\tau) - w_\app(\cdot,\tau)\|_4 \le 
C\e^{-\mu \tau}$ for all $\tau \ge 0$, where
\begin{eqnarray}\nonumber
  w_\app(\xi,\tau) & = & \e^{-\frac{\tau}{2}} \bigl(b_1 F_1(\xi) 
   + b_2 F_2(\xi)\bigr) + \e^{-\tau} \bigl(c_1' H_1(\xi) + c_2' H_2(\xi) 
   + c_3' H_3(\xi) \bigr) \\ \label{wapp}
  && +~\kappa \tau \e^{-\tau} \bigl(-b_1 b_2 H_2(\xi) + (b_1^2-b_2^2)
   H_3(\xi)\bigr)\\ \nonumber
  && +~\e^{-\tau} \bigl(-b_1 b_2 
(\cL+1)^{-1}\Psi_2(\xi) + (b_1^2-b_2^2)(\cL+1)^{-1}
\Psi_3(\xi)\bigr)\ .
\end{eqnarray}
\end{theorem}

\proof Let $w(\xi,\tau)$ be a solution of \reff{SV2} in $L^2_0(4)$, 
and define $\bbeta(\tau), \gamma(\tau), R(\tau)$ by \reff{w3decomp} 
and $\rho(\tau)$ by \reff{rhodef}. Since $\|w(\xi,\tau)\|_4 \to 0$ as 
$\tau \to +\infty$, we can assume without loss of generality that
$\bbeta(0),\ggamma(0),\rho(0)$ satisfy the assumptions of 
Proposition~\ref{invman2} and Lemma~\ref{bgest}. In particular, 
\reff{the-estim} holds for some constants $b_i,c_j \in \real$. The 
estimate $\|w(\cdot,\tau) - w_\app(\cdot,\tau)\|_4 \le C\e^{-\mu \tau}$ 
is then a direct consequence of \reff{w3decomp}, \reff{expliman}, and 
\reff{the-estim}. \QED

\medskip
The terms in \reff{wapp} we wish to call particular attention to are those
proportional to  $\tau \e^{-\tau}$.  When we revert to the original,
unscaled, variables $(x,t)$, these yield terms of the form 
$(1+t)^{-1}\log(1+t)$, which should be contrasted with
previous asymptotic expansions which yielded only inverse powers 
of $t$. Proceeding as in the proof of Corollary~\ref{oseen-gk}, 
we obtain the following result in the original variables:

\begin{corollary}\label{asymp3} Under the assumptions of 
Theorem~\ref{asymp2}, the solution $\omega(x,t)$ of \reff{V2}
with initial data $\omega(x,0) = w_0(x)$  satisfies
$$
   |\omega(\cdot,t) - \omega_\app(\cdot,t)|_p \ \le\ 
   \frac{C_p}{(1+t)^{1+\mu -\frac{1}{p}}}\ , \quad 1\le p \le 2\ ,
   \quad t \ge 0\ ,
$$
where $\omega_\app(x,t) = \frac{1}{1+t}w_\app(\frac{x}{\sqrt{1+t}},
\log(1{+}t))$. Similarly, if $\uu(x,t), \uu_\app(x,t)$ are the velocity 
fields obtained from $\omega(x,t), \omega_\app(x,t)$ via the Biot-Savart
law \reff{BS2}, then
$$
   |\uu(\cdot,t) - \uu_\app(\cdot,t)|_q \ \le\ 
   \frac{C_q}{(1+t)^{\frac{1}{2}+\mu -\frac{1}{q}}}\ , \quad 
   1 \le q < \infty\ , \quad t \ge 0\ .
$$
\end{corollary}

\begin{remark} Two other references which examine the asymptotics
of solutions of \reff{NS2} to this order are \cite{carpio:1996}
and \cite{fujigaki:2000}.  We provide a detailed comparison of
our asymptotics with these two references in \cite{gallay:2001},
but just mention here that since both of these references require
that the initial velocity field satisfy $(1+|x|)\uu_0 \in
L^1(\real^2)^2$ in order to derive the asymptotics, the coefficients
$b_1$ and $b_2$ are both zero for all the solutions they
study. (This follows from Corollary \ref{velint2}.)  
Thus many solutions of finite energy ({\it i.e.} of finite
$L^2$ norm) are excluded from consideration by this
hypothesis, and in particular one does not observe the logarithmic
terms in the asymptotics.  This is a further reason that we feel
it is preferable to impose the decay conditions on the vorticity
field rather than on the velocity.
\end{remark}

\subsection{Optimal Decay Rates}\label{optimal}

It has been known for a long time that there is a relationship
between the spatial and temporal decay rates of solutions of the
Navier-Stokes equation in $\real^N$, $N \ge 2$. For instance, for any 
initial data $\uu_0 \in L^2(\real^N)^N \cap L^1(\real^N)^N$, there exists 
a global weak solution of the Navier-Stokes equation satisfying
\begin{equation}\label{gen_bound}
  |\uu(t)|_2 \le C(1+t)^{-N/4}\ , \quad t \ge 0\ ,
\end{equation}
see \cite{schonbek:1985}, \cite{kajikiya:1986}. If, in addition, 
$(1+|x|)\uu_0 \in L^1(\real^N)^N$, if follows from 
Wiegner's result \cite{wiegner:1987} that 
\begin{equation}\label{better-bound}
  |\uu(t)|_2 \le C(1+t)^{-(N+2)/4}\ , \quad t \ge 0\ .
\end{equation}
In \cite{miyakawa:2000}, T.~Miyakawa and M.E.~Schonbek investigate the 
optimality of the decay rate \reff{better-bound}. More specifically, they 
prove:

\begin{theorem}\label{MiSchon} {\bf \cite{miyakawa:2000}}
Assume that $\uu_0 \in L^2(\real^N)^N$, $\nabla\cdot\uu_0 = 0$, and
$(1+|x|)\uu_0 \in L^1(\real^N)^N$, $N \ge 2$. Let $\uu(t)$ be a global
weak solution of the Navier-Stokes equation with initial data $\uu_0$ 
satisfying \reff{better-bound}. For all $k,\ell \in \{1,\dots,N\}$, define
\begin{equation}\label{bcmom}
  b_{k \ell} = \int_{\real^N} x_{\ell} (\uu_0)_k(x)\d x\ , \quad
  c_{k \ell} = \int_0^{\infty} \int_{\real^N} u_k(x,t) u_{\ell}(x,t) 
  \d x \d t\ .
\end{equation}
Then 
\begin{equation}\label{optidecay}
  \lim_{t\to\infty} t^{\frac{N+2}{4}}|\uu(t)|_2 = 0
\end{equation}
if and only if there exists $c \ge 0$ such that
\begin{equation}\label{more-moments}
  b_{k \ell} = 0 \quad \mbox{and }\quad c_{k \ell} = c \delta_{k \ell}\ ,
  \quad k,\ell \in \{1,\dots,N\}\ .
\end{equation}  
\end{theorem}

The proofs in \cite{miyakawa:2000} are clear, but do not provide much
intuition as to the meaning of the conditions \reff{more-moments}.
As the authors themselves remark, ``We know nothing about the
characterization of solutions satisfying  $(c_{k \ell}) = 
(c \delta_{k \ell})$.'' As we will demonstrate below in the case 
$N = 2$ (see also \cite{gallay:2001} for the three-dimensional case), 
the invariant manifold approach provides a simple 
geometrical explanation of the meaning of both conditions in 
\reff{more-moments}. In addition, it will allow us to construct additional 
solutions satisfying \reff{optidecay}, but which do not fit into the setting
of Theorem~\ref{MiSchon}, because $(1+|x|)\uu_0 \notin L^1(\real^2)^2$.
Note however that our work does require that we consider solutions of 
\reff{SV2} in the space $L^2(4)$, and this decay requirement on the
vorticity does not appear in the work of Miyakawa and Schonbek.

As before, we work with the vorticity equation rather than the Navier-Stokes
equations themselves.  This seems to us a particularly natural choice
when examining the long-time behavior of solutions since, as
we have noted above, the long-time asymptotics are influenced by
assumptions about the spatial decay of the initial conditions, and such
hypotheses are preserved by the evolution of the vorticity
equation (see Theorem~\ref{weightedL2}), but not by the Navier-Stokes
evolution. Since our goal is to recover the results by Miyakawa 
and Schonbek, we restrict ourselves to velocity fields that belong 
to $L^1(\real^2)$. At the level of the vorticity, this condition is 
equivalent to
\begin{equation} \label{moment-cond}
  \int_{\real^2} w(\xi) \d\xi = 0 \ , \quad 
  \int_{\real^2} \xi_1 w(\xi) \d\xi = 0 \ , \quad  
  \int_{\real^2} \xi_2 w(\xi) \d\xi = 0 \ ,
\end{equation}
see Corollary~\ref{velint2}. Thus, we shall study the solutions of the
vorticity equation \reff{SV2} in the invariant subspace of $L^2(4)$ 
defined by \reff{moment-cond}. If $w(\xi,\tau)$ is such a solution and if 
$\vv(\xi,\tau)$ is the velocity field obtained from $w(\xi,\tau)$ 
via the Biot-Savart law \reff{BS2}, then $w$ and $\vv$ can be 
decomposed as follows
\begin{eqnarray}\label{wvvdecomp}
  w(\xi,\tau) &=& \gamma_1(\tau) H_1(\xi) + \gamma_2(\tau) H_2(\xi) + 
   \gamma_3(\tau) H_3(\xi) + R(\xi,\tau)\ ,\\ \nonumber
  \vv(\xi,\tau) &=& \gamma_1(\tau) \vv^{H_1}(\xi) + \gamma_2(\tau) 
   \vv^{H_2}(\xi) + \gamma_3(\tau) \vv^{H_3}(\xi) + \vv^R(\xi,\tau)\ , 
\end{eqnarray}
where the coefficients $\gamma_j$ ($j = 1,2,3$) are defined by 
\reff{gammacoeff}. 
Setting $\bbeta = 0$ in \reff{log1}, we obtain the evolution system
\begin{eqnarray} \nonumber
  \dot{\gamma}_1 &=& - \gamma_1\ ,\\ \label{opt1}
  \dot{\gamma}_2 &=& - \gamma_2 + f_2(0,\ggamma,R)\ ,\\ \nonumber
  \dot{\gamma}_3 &=& - \gamma_3 + f_3(0,\ggamma,R)\ ,\\ \nonumber
  R_\tau &=& \cL R + \nabla\cdot \FF(0,\ggamma,R) - \vv^R \cdot \nabla R\ .
\end{eqnarray}
Also, using \reff{gam-eqn} and Lemma~\ref{proj_lem}, we see that
\begin{equation}\label{f2f3}
  f_2(0,\ggamma,R) = -\int_{\real^2}v_1 v_2\d\xi\ , \quad
  f_3(0,\ggamma,R) = \int_{\real^2} (v_1^2-v_2^2)\d\xi\ ,
\end{equation}
where $\vv$ is given by \reff{wvvdecomp}. 

Now, let $E_c = \span\{H_1,H_2,H_3\}$, and let $E_s \subset L^2(4)$ be the 
spectral subspace of $\cL$ corresponding to the continuous spectrum 
$\sigma_c = \{ \lambda \in \complex ~|~ \Re(\lambda) \le -\frac{3}{2}\}$.
Given $1 < \mu < \frac{3}{2}$, Theorem~\ref{strongman} (with $k = 2$ and
$m = 4$) shows that, for sufficiently small $r_1 > 0$, the set $W_s^{loc}$
defined by \reff{charact} is a smooth (infinite-dimensional) manifold
which is tangent to $E_s$ at the origin. Since the equations for $\alpha$
and $\beta_i$ are linear, it is clear that $W_s^{loc}$ is contained in
the subspace of $L^2(4)$ defined by \reff{moment-cond}, and so is the
global strong-stable manifold $W_s$ defined by \reff{globalss}. 
The following result gives various characterizations of $W_s$. 

\begin{proposition}\label{caractws}
Fix $1 < \mu < \frac{3}{2}$, and assume that $w_0 \in L^2(4)$ satisfies 
\reff{moment-cond}. Let $w(\xi,\tau)$ be the solution of \reff{SV2} with 
initial data $w_0$, and let $\vv(\xi,\tau)$ be the velocity field obtained 
from $w(\xi,\tau)$ via the Biot-Savart law \reff{BS2}. Decompose $w(\xi,\tau)$ 
according to \reff{wvvdecomp}, and define the coefficients $c_{k\ell}$
by \reff{bcmom}, where
\begin{equation}\label{uvvscaling}
  \uu(x,t) = \frac{1}{\sqrt{1+t}}\,\vv\Bigl(\frac{x}{\sqrt{1+t}},
  \log(1+t)\Bigr)\ .
\end{equation}
Then the following statements are equivalent:\\
\null\quad 1) $w_0$ lies in the strong stable manifold $W_s$. \\
\null\quad 2) $\displaystyle{\lim_{\tau\to\infty} \e^\tau \|w(\cdot,\tau)\|_4 
= 0}$.\\
\null\quad 3) $\displaystyle{\lim_{t\to\infty} t|\uu(\cdot,t)|_2 = 0}$.\\
\null\quad 4) $\gamma_1(0) = 0$, $\gamma_2(0) = c_{12}$, and $\gamma_3(0) = 
c_{22}-c_{11}$.
\end{proposition}

\begin{remark}  Note that the global nature of the strong stable
manifold means that Proposition~\ref{caractws} is not limited to 
solutions of small norm, but applies to solutions of arbitrary
size.
\end{remark}

\proof Applying Theorem~\ref{asymp2} with $b_1 = b_2 = 0$, we obtain
\begin{equation}\label{cjestim}
  \|w(\cdot,\tau) - (c_1 H_1 + c_2 H_2 + c_3 H_3)\e^{-\mu\tau}\|_4 
  \le C \e^{-\tau}\ ,
\end{equation}
for some $c_1,c_2,c_3 \in \real$. We shall show that statements 1), 2), 
3), 4) in Proposition~\ref{caractws} are all equivalent to\\
\null\quad 5) $c_1 = c_2 = c_3 = 0$.\\
Indeed, since the functions $H_1, H_2, H_3$ are linearly independent, 
it is clear from \reff{charact} and \reff{cjestim} that 1) 
$\Leftrightarrow$ 2) $\Leftrightarrow$ 5). On the other hand, it 
follows from \reff{wvvdecomp} and \reff{opt1} that the coefficients
$c_j$ in \reff{cjestim} are given by the formulas
$$
  c_1 = \gamma_1(0)\ , \quad c_j = \gamma_j(0) 
  + \int_0^\infty \e^\tau f_j(0,\gamma(\tau),R(\tau))\d\tau\ , \quad j=2,3\ ,
$$
see also \reff{cjdef}. Furthermore, using \reff{f2f3} and the change of 
variables \reff{uvvscaling}, it is straightforward to verify that 
\begin{eqnarray*}
  c_2 &=& \gamma_2(0) - \int_0^\infty \e^\tau \int_{\real^2} v_1(\xi,\tau)
   v_2(\xi,\tau)\d\xi \d\tau \ \equiv\ \gamma_2(0) - c_{12}\ ,\\
  c_3 &=& \gamma_3(0) + \int_0^\infty \e^\tau \int_{\real^2} (v_1(\xi,\tau)^2
   -v_2(\xi,\tau)^2)\d\xi \d\tau \ \equiv\ \gamma_2(0) + c_{11}-c_{22}\ ,\\
\end{eqnarray*} 
where $c_{k\ell}$ is defined in \reff{bcmom}. Therefore, 4) 
$\Leftrightarrow$ 5). Finally, from Corollary~\ref{asymp3}, we have
$|\uu(\cdot,t) - \uu_\app(\cdot,t)|_2 \le C (1+t)^{-\mu}$, where
$$
  \uu_\app(x,t) = \frac{1}{(1+t)^{3/2}} \sum_{j=1}^3 c_j 
  \vv^{H_j}\Bigl(\frac{x}{\sqrt{1+t}}\Bigr)\ .
$$
Clearly, $|\uu_\app(\cdot,t)|_2 = \frac{K}{1+t}$, where $K = 0$ if and
only if $c_1 = c_2 = c_3 = 0$. Thus, 3) $\Leftrightarrow$ 5). \QED

\medskip
We are now able to give an alternative proof of Theorem~\ref{MiSchon}
in the particular case where $N = 2$ and $w_0 = \partial_1 (\uu_0)_2 
-\partial_2 (\uu_0)_1 \in L^2(4)$. Indeed, let $w(\xi,\tau)$ be the solution
of \reff{SV2} with initial data $w_0$. Since $\uu_0 \in L^1(\real^2)^2$, 
it follows from Corollary~\ref{velint2} that \reff{moment-cond} holds, 
hence $w(\xi,\tau)$ can be decomposed according to \reff{wvvdecomp}. 
Moreover, by Corollary~\ref{velint3}, the assumption $(1+|x|)\uu_0 \in 
L^1(\real^2)^2$ is equivalent to $\gamma_2(0) = \gamma_3(0) = 0$.

\smallskip\noindent{\bf a)}
Assume first that \reff{optidecay} holds, namely $w_0 \in W_s$. 
Then point 4) in Proposition~\ref{caractws} shows that $c_{12} = 0$ and 
$c_{11} = c_{22}$, hence the matrix $(c_{k\ell})$ is scalar. In addition,
$\gamma_1(0) = 0$, hence $b_{k\ell} = 0$ by Corollary~\ref{velint3}. 

\smallskip\noindent{\bf b)}
Conversely, assume that \reff{more-moments} holds. Then $\gamma_1(0) = 0$
by Corollary~\ref{velint3}, and since $\gamma_2(0) = 0 = c_{12}$, 
$\gamma_3(0) = 0 = c_{22}-c_{11}$ it follows from Proposition~\ref{caractws}
that $w_0 \in W_s$. This concludes the proof. \QED

\medskip
Note that as a consequence of these investigations, we see that
there are many other solutions $\uu(x,t)$ of the Navier-Stokes equation that
lie in the strong stable manifold (and hence satisfy the decay
estimate \reff{optidecay}), but which do not satisfy the 
moment condition $(1+|x|)\uu_0 \in L^1(\real^2)^2$.  Put another way, 
the conditions \reff{more-moments} are the necessary and sufficient 
conditions for \reff{optidecay} to hold within the class of solutions
satisfying $(1+|x|)\uu_0 \in L^1(\real^2)^2$, but if one looks at larger 
classes of initial conditions then one finds that \reff{more-moments}
is no longer necessary, and that a better characterization of solutions
satisfying \reff{optidecay} is that they lie in the strong stable
manifold $W_s$.  Indeed, as we observed above, the velocity
fields $\vv^{H_2}$ and $\vv^{H_3}$ corresponding to the vorticities 
$H_2$ and $H_3$ do not satisfy $(1+|\xi|) \vv^{H_j} \in L^1(\real^2)^2$. 
So if we choose as an initial condition any point in the strong stable 
manifold $W_s$ which does not lie in the hyperplane $E_s$ of
functions of the form \reff{wvvdecomp} with $\gamma_1(0) = \gamma_2(0) 
= \gamma_3(0)=0$, we obtain a solution of the Navier-Stokes equation
satisfying \reff{optidecay}, but not \reff{more-moments}.

Of course, this argument requires that there be some point in 
$W_s$ for which either $\gamma_2$ or $\gamma_3$ is nonzero. 
(Proposition~\ref{caractws} implies that all points in $\wsloc$
must have $\gamma_1 = 0$.) This could fail to occur only if $W_s$
coincided with the hyperplane $E_s$, namely if $E_s$ was invariant
under the semiflow defined by \reff{log1}. This in turn would happen 
only if $f_2(0,0,R) = f_3(0,0,R)=0$ for all $R \in E_s$. We now demonstrate 
explicitly that this is not the case. Choose 
$$
  K(\xi) = \partial_1 H_1(\xi) = \xi_1 (1-|\xi|^2/8)G(\xi)\ .
$$
Then $\cL K = -\frac{3}{2} K$, so that $K \in E_s$. The velocity field 
corresponding to $K$ is
$$
  \vv^K(\xi) = \partial_1 \vv^{H_1}(\xi) = \frac{G(\xi)}{4}\pmatrix{
  -\xi_1\xi_2 \cr \xi_1^2-2}\ ,
$$
and a direct calculation gives:
$$
  f_3(0,0,K) = \int_{\real^2} \bigl((v_1^K)^2 -(v_2^K)^2\bigr)\d\xi 
  = -\frac{1}{64\pi} \neq 0\ .
$$
This shows that $W_s \not\subset E_s$. On the other hand, it is easy to see 
that any $w \in E_s$ which is {\sl radially symmetric} actually lies
in $W_s$, because \reff{V2} reduces to the heat equation for radially
symmetric vorticities. Thus the intersection $W_s \cap E_s$ contains
at least an infinite-dimensional subspace of $L^2(m)$. 

\figurewithtex 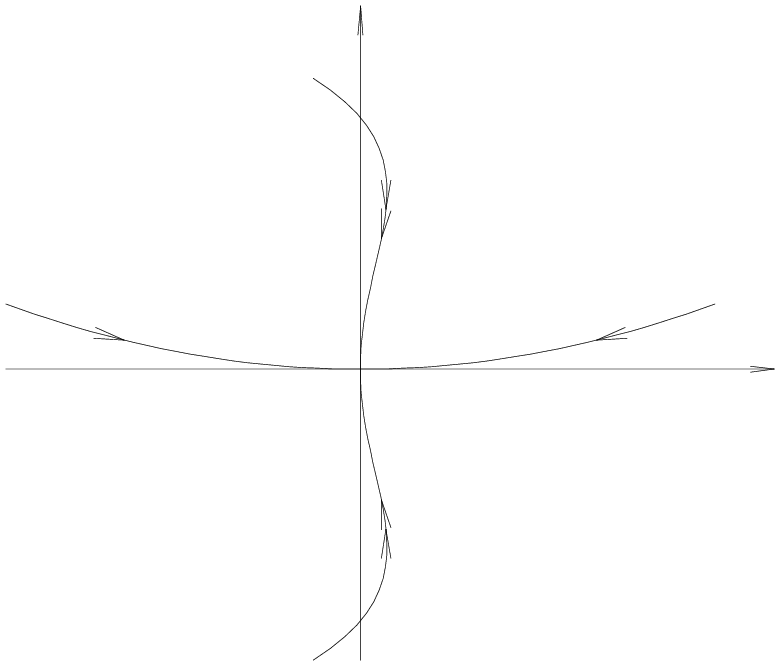 Fig3.tex 6.5 9.0 {\bf Fig. 3.} A schematic picture
of the dynamics defined by equations \reff{opt1} in a neighborhood of
the origin. All trajectories approach a three-dimensional invariant manifold 
$W_c$ at a rate $\cO(\e^{-\mu\tau})$ or faster. The strong-stable manifold
$W_s$ contains all solutions that converge to the origin faster than 
$\e^{-\tau}$. The intersection $W_s \cap E_s$, which is depicted here 
by two points only, contains in fact an infinite-dimensional subspace
of $L^2(m)$.\cr

Summing up, we see that for solutions of the Navier-Stokes equation
in $\real^2$, whose vorticity lies in $L^2(4)$, the stong-stable
manifold identifies exactly those which converge to zero faster
than $C/t$.  If the velocity field at $t=0$ satisfies
$(1+|x|)\uu \in L^1(\real^2)^2$, this gives a natural geometrical
interpretation of the moment conditions of \cite{miyakawa:2000}.
It also shows, however, that there are additional solutions which
decay with a rate faster than $C/t$ but which do not satisfy
the decay condition on the initial velocity field.

We conclude this section by describing a somewhat surprising property
of the strong-stable manifold.  Although it was constructed
to be invariant with respect to the dynamics of the vorticity
equation in the rescaled variables, it is also invariant
with respect to the dynamics
expressed in the unscaled variables.  Let $\Psi_t^m$ be the
semiflow on $L^2(m)$ defined by \reff{V2}.

\begin{proposition} The semiflow $\Psi_t^m$ leaves the manifold
$W_s$ invariant.
\end{proposition}

\begin{remark} This should be contrasted with the behavior
of $\wcloc$.  In order to obtain a manifold corresponding to
$\wcloc$ that is invariant with respect to $\Psi_t^m$ it
is necessary to go the extended phase space $L^2(m)\times \real_+$,
where $\real_+$ represents the time axis.  (See the discussion in 
Section~3 of \cite{wayne:1997}, in particular Theorem~3.6, where this 
question is discussed in a related context.)
\end{remark}

\proof The proposition follows from the characterization of
$W_s$ in point 3 of Proposition~\ref{caractws}.  Since $\uu(\cdot,t)$
is the velocity field corresponding to the vorticity field $\Psi_t^m(w_0)$, 
and since this condition is clearly invariant with respect to time 
translation, the invariance of $W_s$ with respect to $\Psi_t^m$ is 
immediate. The fact that $W_s$ is invariant with respect to both $\Psi_t^m$ 
and $\Phi_t^m$ is related to the fact that solutions in $W_s$ can be 
characterized by their asymptotics as $t \to +\infty$. To be more explicit, 
note that combining point~2 in Proposition~\ref{caractws}, with point~5 in 
the proof of that proposition, we see that $w_0$ lies in $W_s$ if and only 
if the solution $w(\cdot,\tau)$ of \reff{SV2} with initial data $w_0$ 
satisfies
\begin{equation}\label{wwdecay}
  \lim_{\tau\to+\infty} \e^\tau |w(\cdot,\tau)|_2 = 0\ .
\end{equation}
But letting $\omega(x,t) = \frac{1}{1+t} w(\frac{x}{\sqrt{1+t}},\log(1+t))$,
we see that \reff{wwdecay} is equivalent to
\begin{equation}\label{wdecay}
  \lim_{t\to+\infty} t^{3/2} |\omega(\cdot,t)|_2 = 0\ .
\end{equation}
Note also that $\omega(x,0) = w_0(x)$.  Thus, $W_s$ can also be 
characterized as the set of points in $L^2(4)$ satisfying \reff{moment-cond} 
and for which 
$$
  \lim_{t\to+\infty} t^{3/2} | \Psi_t^m(w_0)|_2 = 0\ .
$$
Since this characterization is again invariant with respect to
$\Psi_t^m(w_0)$, this shows that $W_s$ is an invariant manifold in the 
original variables as well as in the rescaled variables. \QED


\appendix
\section{Spectrum of the linear operator} \label{spectrum}

In this section, we assume that $N \in \intplus$, $N \ge 1$. We 
consider the linear operator $\cL$ given by
\begin{equation}\label{Lop}
  \cL = \Delta_\xi + \frac{1}{2}\xxi \cdot \nabla_\xi + \frac{N}{2}\ ,
  \quad \xi \in \real^N\ .
\end{equation}
As in Section~\ref{scaling}, we shall work in the weighted space 
$L^2(m)$ defined by
\begin{eqnarray}\label{L2m}
  L^2(m) &=& \bigl\{ f \in L^2(\real^N) \,|\, \|f\|_m < \infty \bigr\}\ ,\\
  \nonumber
  \|f\|_m &=& \left(\int_{\real^N}(1+|\xi|^2)^m |f(\xi)|^2 \d\xi 
  \right)^{1/2}\ .
\end{eqnarray}
Since its coefficients depend linearly on the space variable $\xi$, the 
operator $\cL$ becomes a first order differential operator when expressed
in the Fourier variable $p$. Our convention for Fourier transformation is
\begin{eqnarray}\label{Fourierdef}
  \hat f(p) &=& \int_{\real^N} f(\xi) \exp(-\I\pp\cdot\xxi)\d\xi\ ,\\
  f(\xi) &=& \frac{1}{(2\pi)^N} \int_{\real^N} \hat f(p) \exp(\I\pp\cdot\xxi)
  \d p\ .
\end{eqnarray} 
The expression of the operator $\cL$ in Fourier space is:
\begin{equation}\label{cLF}
  (\widehat{\cL f})(p) = -(|p|^2 + \half \pp\cdot\nabla_p)\hat f(p)\ .
\end{equation}

\medskip
The aim of this section is to prove the following result, which underlies
our approach for computing the long-time asymptotics of the vorticity
equation:

\begin{theorem}\label{Lspectrum}
Fix $m \ge 0$, and let $\cL$ be the linear operator \reff{Lop} in 
$L^2(m)$, defined on its maximal domain. Then the spectrum of 
$\cL$ is
$$
  \sigma(\cL) = \Bigl\{\lambda \in \complex\,\Big|\, \Re(\lambda) \le 
  \frac{N}{4} - \frac{m}{2}\Bigr\} \cup \Bigl\{-\frac{k}{2} \,\Big|\, 
  k \in \intplus\Bigr\}\ .
$$
Moreover, if $m > \frac{N}{2}$ and if $k \in \intplus$ satisfies 
$k + \frac{N}{2} < m$, then $\lambda_k = -\frac{k}{2}$ is an isolated 
eigenvalue of $\cL$, with multiplicity ${N+k-1 \choose k}$.
\end{theorem}

In the one-dimensional case, this result is proved in \cite{Gallay:1998}, 
Appendix~A. We give here a slightly different proof, which is valid 
for all $N \ge 1$. We begin with a few elementary observations.

\medskip
\noindent{\bf (1)} {\it The discrete spectrum of} $\cL$. Fix $k \in 
\intplus$, and take $\alpha = (\alpha_1,\dots,\alpha_N) \in \intplus^N$
such that $|\alpha| = \alpha_1 + \dots + \alpha_N = k$. Then the Hermite
function $\phi_\alpha \in \cS(\real^N)$ defined by
\begin{eqnarray}\nonumber
  \hat \phi_\alpha(p) &=& (\I p)^\alpha \e^{-|p|^2} \equiv \I^{|\alpha|}
  p_1^{\alpha_1}\dots p_N^{\alpha_N} \e^{-|p|^2}\ , \ \ \hbox{or}\\
  \label{phialpha}
  \phi_\alpha(\xi) &=& (\partial_\xi^\alpha \phi_0)(\xi)\ ,\ \
  \phi_0(\xi) = \frac{1}{(4\pi)^{N/2}}\,\e^{-|\xi|^2/4}\ ,
\end{eqnarray}
is an eigenfunction of $\cL$ with eigenvalue $-\frac{k}{2}$. Thus, for any 
$m \ge 0$, we have $\sigma(\cL) \supset \{-\frac{k}{2}\,|\,k \in \intplus\}$, 
and the multiplicity of the eigenvalue $\lambda_k = -\frac{k}{2}$ is 
greater or equal to ${N+k-1 \choose k} = \#\{\alpha \in \intplus^N\,|\,
|\alpha| = k\}$. 

\medskip
\noindent{\bf (2)} {\it The ``continuous'' spectrum of} $\cL$. Fix 
$\lambda \in \complex$ such that $\Re(\lambda) < N/4$ and $-\lambda 
\notin \intplus$. The function $\psi_\lambda : \real^N \to \real$ 
defined (in Fourier variables) by
$$
  \hat\psi_\lambda(p) = |p|^{-2\lambda}\,\e^{-|p|^2}\ , \quad 
  p \in \real^N\ ,
$$
is then an eigenfunction of $\cL$ with eigenvalue $\lambda$. It is
clear that $\psi_\lambda \in C^\infty(\real^N)$, and a standard 
calculation (see for instance \cite{guelfand:1962}, Section~2.3.3) shows 
that
$$
  \lim_{|\xi| \to \infty} |\xi|^{N-2\lambda} \psi_\lambda(\xi) =
  \frac{\Gamma(\frac{N}{2}-\lambda)}{2^{2\lambda}\pi^{\frac{N}{2}}
  \Gamma(\lambda)} \neq 0\ .
$$
In particular, $\psi_\lambda \in L^2(m)$ if and only if $\Re(\lambda) < 
\frac{N}{4}-\frac{m}{2}$. Since the spectrum of $\cL$ is closed, this shows 
that $\sigma(\lambda) \supset \{\lambda \in \complex\,|\, \Re(\lambda) \le 
\frac{N}{4}-\frac{m}{2}\}$.

\medskip
\noindent{\bf (3)} {\it The spectral projections.} For all $\alpha = 
(\alpha_1,\dots,\alpha_N) \in \intplus^N$, we define the Hermite 
polynomial $H_\alpha$ by
$$
  H_\alpha(\xi) = \frac{2^{|\alpha|}}{\alpha!} \,\e^{|\xi|^2/4}
  \partial^\alpha_\xi
  \left(\e^{-|\xi|^2/4}\right)\ ,\ \ \xi \in \real^N\ ,
$$
where $\alpha! = (\alpha_1!)\cdot\dots\cdot(\alpha_N!)$. It is
not difficult to verify that $H_\alpha$ is a polynomial of degree
$|\alpha|$ which satisfies $\cL^*H_\alpha = -\frac{|\alpha|}{2}H_\alpha$,
where $\cL^* = \Delta_\xi -\frac{1}{2}\xxi\cdot\nabla_\xi$ is the formal
adjoint of $\cL$. In addition, the following orthogonality relations
hold for all $\alpha,\beta \in \intplus^N$:
$$
  \int_{\real^N} H_\alpha(\xi)\phi_\beta(\xi)\d\xi = \delta_{\alpha\beta}\ .
$$
Finally, for any smooth function $\hat f(p)$, we have the relation
$$
  \left(H_\alpha(\I\nabla_p)\hat f\right)(0) = \frac{1}{\alpha!}
  \partial^\alpha_p \left(\hat f(p)\e^{|p|^2}\right)\Big|_{p=0}\ .
$$
Assume now that $m > \frac{N}{2}$. For all $n \in \allint$ such that 
$n + \frac{N}{2} < m$, we define a continuous projection $P_n : L^2(m) 
\to L^2(m)$ by the formula
$$
  (P_n f)(\xi) = \sum_{|\alpha| \le n} \left(\int_{\real^N}H_\alpha(\xi')
  f(\xi')\d\xi'\right)^{1/2}\phi_\alpha(\xi)\ .
$$
We also set $Q_n = \oone - P_n$. Remark that $P_n = 0$ and $Q_n = \oone$
if $n < 0$. If $n \ge 0$, it is clear from the definitions that 
$P_n$ is the spectral projection onto the $\sum_{k=0}^n
{N+k-1 \choose k}$--dimensional subspace spanned by the eigenfunctions
of $\cL$ corresponding to the eigenvalues $\{-\frac{k}{2}\,|\,k = 0,1,\dots,
n\}$. For later use, we note that the condition $P_n f = 0$ is 
equivalent to
$$
  \int_{\real^N} \xi^\alpha f(\xi)\d\xi = 0 \ \ \hbox{for all}\ 
  \alpha \in \intplus^N \ \hbox{with}\ |\alpha| \le n\ .
$$

\noindent{\bf (4)} {\it The semigroup} $\e^{\tau\cL}$. The 
operator $\cL$ is the generator of a linear semigroup $S(\tau) = 
\e^{\tau\cL}$ given by the following expressions:
\begin{eqnarray*}
  (\widehat{S(\tau)f})(p) &=& \e^{-a(\tau)|p|^2} \hat f(p\,\e^{-\tau/2})\ ,
  \quad \hbox{or}\\
  (S(\tau)f)(\xi) &=& \frac{\e^\frac{N\tau}{2}}{(4\pi a(\tau))^\frac{N}{2}} 
  \int_{\real^N} \exp\Bigl(-\frac{|\xi-\xi'|^2}{4a(\tau)}\Bigr) 
  f(\xi' \e^{\tau/2})\d\xi'\ , 
\end{eqnarray*}
where $a(\tau) = 1-\e^{-\tau}$. It follows from these formulas that
$S(\tau)$ is a strongly continuous (but not analytic) semigroup on $L^2(m)$ 
for any $m \ge 0$.

\medskip
The main technical result of this section is the following estimate
on the semigroup $S(\tau) = \e^{\tau\cL}$:

\begin{proposition}\label{Sgestim}
{\bf (a)} Fix $m \ge 0$, and take $n \in \allint$ such that
$n + \frac{N}{2} < m \le n+1 + \frac{N}{2}$. For all $\alpha \in \intplus^N$
and all $\epsilon > 0$, there exists $C > 0$ such that
\begin{equation}\label{Sg1}
  \|\partial^\alpha S(\tau)Q_n f\|_m \le \frac{C}{a(\tau)^{|\alpha|/2}}
  \,\e^{\frac{\tau}{2}(\frac{N}{2}-m+\epsilon)} \|f\|_m\ ,
\end{equation}
for all $f \in L^2(m)$ and all $\tau > 0$.

\smallskip\noindent
{\bf (b)} Fix $n \in \intplus \cup \{-1\}$, and take $m \in \real$ such that
$m > n+1+\frac{N}{2}$. For all $\alpha \in \intplus^N$
and all $\epsilon > 0$, there exists $C > 0$ such that
\begin{equation}\label{Sg2}
  \|\partial^\alpha S(\tau)Q_n f\|_m \le \frac{C}{a(\tau)^{|\alpha|/2}}
  \,\e^{-\frac{n+1}{2}\tau} \|f\|_m\ ,
\end{equation}
for all $f \in L^2(m)$ and all $\tau > 0$.
\end{proposition}

\begin{remark}
Estimate \reff{Sg2} has been obtained in \cite{wayne:1997}, 
\cite{eckmann:1997}, and \cite{eckmann:1998} under the assumption
that $m$ is sufficiently large, depending on $n$. As is clear from 
Theorem~\ref{Lspectrum}, the condition $m > n+1+\frac{N}{2}$ is 
optimal.
\end{remark}

Using Proposition~\ref{Sgestim}, it is easy to complete the proof of 
Theorem~\ref{Lspectrum}. Indeed, if $m$ and $n$ are as in part {\bf (a)}
of Proposition~\ref{Sgestim}, then $\sigma(\cL) = \sigma(\cL P_n) \cup
\sigma(\cL Q_n)$. By construction, $\sigma(\cL P_n) = \emptyset$ if
$n < 0$ and $\sigma(\cL P_n) = \{0;-\frac{1}{2};\dots;-\frac{n}{2}\}$ if
$n \in \intplus$. On the other hand, by the Hille-Yosida theorem 
(see for instance \cite{pazy:1983}), the bound \reff{Sg1} implies that 
$\sigma(\cL Q_n) \subset \{\lambda \in \complex\,|\, \Re(\lambda) \le 
\frac{N}{4}-\frac{m}{2}\}$. Thus
$$
  \sigma(\cL) \subset \Bigl\{\lambda \in \complex\,\Big|\, \Re(\lambda) \le 
  \frac{N}{4} - \frac{m}{2}\Bigr\} \cup \Bigl\{-\frac{k}{2} \,\Big|\, 
  k \in \intplus\Bigr\}\ ,
$$
and the reverse inclusion has already been established (see {\bf (1)} and
{\bf (2)} above). Finally, since $\sigma(\cL P_n) \cap \sigma(\cL Q_n) =
\emptyset$, we see that the multiplicity of the eigenvalue 
$\lambda_k = -\frac{k}{2}$ ($k = 0,\dots,n$) is exactly ${N+k-1 \choose
k}$. This concludes the proof of Theorem~\ref{Lspectrum}. \QED

\vskip 5pt
\noindent{\bf Proof of Proposition~\ref{Sgestim}.} We first show that
{\bf (a)} $\Rightarrow$ {\bf (b)}. Let $n \in \intplus \cup \{-1\}$ and
$m > n+1+\frac{N}{2}$. Taking $\bar n \in \intplus$ such that $\bar n
+\frac{N}{2} < m \le \bar n + 1 +\frac{N}{2}$, we can write, for any $f \in 
L^2(m)$, 
$$
  Q_n f = (\oone - P_n)f = (P_{\bar n} - P_n)f + Q_{\bar n} f\ .
$$
Since $m - \frac{N}{2} > \bar n \ge n+1$, we have by \reff{Sg1}
$$
  \|\partial^\alpha S(\tau)Q_{\bar n} f\|_m \le \frac{C}{a(\tau)^{|\alpha|/2}}
  \,\e^{-\frac{n+1}{2}\tau} \|f\|_m\ , \ \ \tau > 0\ .
$$
On the other hand, since $S(\tau)\phi_\alpha = \e^{-|\alpha|\tau/2}
\phi_\alpha$ and
$$
  (P_{\bar n} -P_n)f = \sum_{n < |\alpha| \le {\bar n}} 
  \left(\int_{\real^N}H_\alpha(\xi)f(\xi)\d\xi \right)^{1/2}\phi_\alpha\ ,
$$
we easily obtain $\|\partial^\alpha S(\tau)(P_{\bar n} -P_n)f\|_m \le
C\e^{-\frac{n+1}{2}\tau} \|f\|_m$. This proves \reff{Sg2}.

\medskip To prove {\bf (a)}, we first remark that, for all $\alpha \in 
\intplus^N$, 
\begin{equation}\label{DalphaS}
  (\partial^\alpha S(\tau)f)(\xi) = \frac{\e^\frac{N\tau}{2}}{a(\tau)^\frac{N
  +|\alpha|}{2}} \int_{\real^N} \phi_\alpha\Bigl(
  \frac{\xi-\xi'}{a(\tau)^{1/2}}\Bigr)f(\xi'\e^\frac{\tau}{2})\d\xi'\ ,
\end{equation}
where $\phi_\alpha$ is defined in \reff{phialpha}. Since $1+|\xi|^m \le
C(1+|\xi-\xi'|^m)(1+|\xi'|^m)$, we have
\begin{eqnarray*}
  (1+|\xi|^m)|(\partial^\alpha S(\tau)f)(\xi)| &\le& 
  C\frac{\e^\frac{N\tau}{2}}{a(\tau)^\frac{N+|\alpha|}{2}} \int_{\real^N} 
  (1+|\xi-\xi'|^m) \Big|\phi_\alpha\Bigl(\frac{\xi-\xi'}{a(\tau)^{1/2}}
  \Bigr)\Big|\\
  && \hspace{3cm} \times (1+|\xi'|^m) |f(\xi'\e^\frac{\tau}{2})|\d\xi'\ .
\end{eqnarray*}
Observing that
$$
  \frac{1}{a(\tau)^\frac{N}{2}} \int_{\real^N} (1+|\eta|^m) \Big|\phi_\alpha
  \Bigl(\frac{\eta}{a(\tau)^{1/2}}\Bigr)\Big| \d\eta \,\le\, C\ ,
$$
and applying Young's inequality $|g*h|_2 \le |g|_1 |h|_2$, we obtain
\begin{equation}\label{prebound}
  \|\partial^\alpha S(\tau)f\|_m \le C \frac{\e^\frac{N\tau}{4}}{
  a(\tau)^\frac{|\alpha|}{2}} \|f\|_m\ ,\ \ \tau > 0\ .
\end{equation}
In view of \reff{prebound}, it is sufficient to prove \reff{Sg1} for
$\tau \ge 1$. Since $1 - \e^{-1} \le a(\tau) \le 1$ when $\tau \ge 1$, 
we may drop for simplicity the factors $a(\tau)$ in \reff{DalphaS}.
Thus, all we really need to prove is:

\begin{lemma}\label{Sgaux}
Let $\phi \in \cS(\real^N)$, $m \ge 0$, and $n \in \allint$ such that
$n+\frac{N}{2} < m \le n+1+\frac{N}{2}$. For any $f \in L^2(m)$, 
define
$$
  (\bar S(\tau)f)(\xi) = \int_{\real^N} \phi(\xi-\eta)f(\eta 
  \e^\frac{\tau}{2}) \d\eta\ ,\ \ \tau \ge 0\ .  
$$
Then, for all $\epsilon > 0$, there exists $C > 0$ such that, for all
$f \in L^2(m)$, 
\begin{equation}\label{bdbar}
  \|\bar S(\tau)Q_n f\|_m \le C \e^{-\frac{\tau}{2}(m+\frac{N}{2}-\epsilon)}
  \|f\|_m\ ,\ \ \tau \ge 0.
\end{equation}
\end{lemma}

\proof Expanding $\phi$ in Taylor series, we obtain for all $\xi,\eta \in 
\real^N$:
$$
  \phi(\xi-\eta) = \sum_{|\alpha| \le n}\frac{(-1)^{|\alpha|}}{\alpha!}
  (\partial^\alpha \phi)(\xi)\eta^\alpha + \sum_{|\alpha| = n+1}
  \frac{(-1)^{n+1}}{\alpha!} \Phi_\alpha(\xi,\eta) \eta^\alpha\ ,
$$
where
$$
  \Phi_\alpha(\xi,\eta) = (n+1)\int_0^1 (1-s)^n (\partial^\alpha \phi)
  (\xi-s\eta)\d s\ .
$$
If $g = Q_n f$, then $P_n g = 0$, hence $\int_{\real^N}\eta^\alpha g(\eta)
\d\eta = 0$ for all $\alpha \in \intplus^N$ with $|\alpha| \le n$. 
Therefore
\begin{eqnarray*}
  (\bar S(\tau)g)(\xi) &=& \int_{\real^N} \biggl(\phi(\xi-\eta) - 
  \sum_{|\alpha| \le n}\frac{(-1)^{|\alpha|}}{\alpha!}(\partial^\alpha 
  \phi)(\xi)\eta^\alpha\biggr) g(\eta \e^\frac{\tau}{2})\d\eta \\
  &=& (\bar S_1(\tau)g)(\xi) +(\bar S_2(\tau)g)(\xi) +
  (\bar S_3(\tau)g)(\xi)\ , 
\end{eqnarray*}
where
\begin{eqnarray*}
  (\bar S_1(\tau)g)(\xi) &=& \int_{|\eta| \ge 1} \phi(\xi-\eta)g(\eta 
  \e^\frac{\tau}{2})\d\eta\ ,\\
  (\bar S_2(\tau)g)(\xi) &=& -\sum_{|\alpha| \le n}\frac{(-1)^{|\alpha|}}%
  {\alpha!}(\partial^\alpha \phi)(\xi)\int_{|\eta| \ge 1}\eta^\alpha
  g(\eta \e^\frac{\tau}{2})\d\eta\ ,\\
  (\bar S_3(\tau)g)(\xi) &=& \sum_{|\alpha| = n+1} \frac{(-1)^{n+1}}{\alpha!} 
  \int_{|\eta| \le 1}\Phi_\alpha(\xi,\eta)\eta^\alpha g(\eta 
  \e^\frac{\tau}{2})\d\eta\ .
\end{eqnarray*}
Since $1+|\xi|^m \le C(1+|\xi-\eta|^m)|\eta|^m$ when $|\eta| \ge 1$, we have
$$
  (1+|\xi|^m) |(\bar S_1(\tau)g)(\xi)| \le C \int_{|\eta|\ge 1} 
  (1+|\xi-\eta|^m) |\phi(\xi-\eta)| |\eta|^m |g(\eta \e^\frac{\tau}{2})|
  \d\eta\ .
$$
Applying Young's inequality, we conclude that 
\begin{equation}\label{bdfirst}
  \|(\bar S_1(\tau)g)\|_m \le C \e^{-\frac{\tau}{2}(m+\frac{N}{2})}
  \|g\|_m\ ,\ \ \tau \ge 0\ .
\end{equation}
On the other hand, if $|\alpha| \le n$, we have by H\"older's inequality
\begin{eqnarray*}
  \Big|\int_{|\eta|\ge 1} \eta^\alpha g(\eta \e^\frac{\tau}{2})\d\eta\Big|
  &\le& \int_{|\eta|\ge 1} |\eta|^{n-m} (|\eta|^m |g(\eta \e^\frac{\tau}{2})|)
  \d\eta \\
  &\le& \left(\int_{|\eta|\ge 1} |\eta|^{2(n-m)}\d\eta\right)^{1/2} 
  \e^{-\frac{\tau}{2}(m+\frac{N}{2})} \|g\|_m\ ,
\end{eqnarray*}
hence
\begin{equation}\label{bdsecond}
  \|(\bar S_2(\tau)g)\|_m \le C \e^{-\frac{\tau}{2}(m+\frac{N}{2})}
  \|g\|_m\ ,\ \ \tau \ge 0\ .
\end{equation}
To bound the last term $\bar S_3(\tau)g$, let $\alpha \in \intplus^N$ such 
that $|\alpha| = n+1$, and define $\Psi_\alpha(\xi,\eta) = (1+|\xi|^m) 
\Phi_\alpha(\xi,\eta)$. For any $r \ge 1$, it is easy to verify that 
$$
  \sup_{\xi \in \real^N} \left(\int_{|\eta|\le 1}|\Psi_\alpha(\xi,\eta)|^r
  \d\eta\right)^{1/r} + \sup_{|\eta| \le 1} \left(\int_{\real^N}|
  \Psi_\alpha(\xi,\eta)|^r \d\xi\right)^{1/r} < \infty\ .
$$
As is well-known, this implies that the linear operator defined by the
integral kernel $\Psi_\alpha(\xi,\eta)$ is bounded from 
$L^q(\{|\eta| \le 1\})$ into $L^p(\real^N)$ for all $q \le p$. We 
now distinguish two cases: 

\smallskip\noindent{\sl Case 1.} If $m \le n+1$ (which is possible only
if $N = 1$), we have
\begin{eqnarray}\label{S3first}
  \hspace{-0.6cm}
  (1+|\xi|^m)|(\bar S_3(\tau)g)(\xi)| &\le& \sum_{|\alpha|=n+1}\frac{1}
  {\alpha!} \int_{|\eta| \le 1}|\Psi_\alpha(\xi,\eta)||\eta|^{n+1}
  |g(\eta \e^\frac{\tau}{2})| \d\eta\\ \label{S3second}
  &\le& \sum_{|\alpha|=n+1}\frac{1}{\alpha!} 
  \int_{|\eta| \le 1}|\Psi_\alpha(\xi,\eta)||\eta|^m
  |g(\eta \e^\frac{\tau}{2})| \d\eta\ .
\end{eqnarray}
Using \reff{S3second} and the remark above with $p=q=2$, we immediately 
obtain
\begin{equation}\label{bdthird}
  \|(\bar S_3(\tau)g)\|_m \le C \e^{-\frac{\tau}{2}(m+\frac{N}{2})}
  \|g\|_m\ ,\ \ \tau \ge 0\ .
\end{equation}

\smallskip\noindent{\sl Case 2.} If $n+1 < m \le n+1+\frac{N}{2}$, 
we define $q_* \in [1,2)$ by the relation $N(\frac{1}{q_*}-\frac{1}{2}) = 
m - (n+1)$. Given any $\epsilon > 0$, we choose $q \in (q_*,2)$ such that
$N(\frac{1}{q_*}-\frac{1}{q}) \le \epsilon$. As is easily verified, 
there exists $C > 0$ such that $||\eta|^{n+1}g(\eta)|_q \le \|g\|_m$ 
for all $g \in L^2(m)$. Using \reff{S3first} and the remark above with
$p = 2$, we obtain
\begin{eqnarray}\nonumber
  \|(\bar S_3(\tau)g)\|_m &\le& C\left(\int_{|\eta| \le 1}\bigl(|\eta|^{n+1}
  |g(\eta \e^\frac{\tau}{2})|\bigr)^q \d\eta\right)^{1/q}\\ \label{bdfourth}
  &\le& C \e^{-\frac{\tau}{2}(n+1+\frac{N}{q})}||\eta|^{n+1}g(\eta)|_q
  \le C \e^{-\frac{\tau}{2}(m+\frac{N}{2}-\epsilon)} \|g\|_m\,
\end{eqnarray}
since $n+1+\frac{N}{q} = m+\frac{N}{2}-N(\frac{1}{q_*}-\frac{1}{q})
\ge m+\frac{N}{2}-\epsilon$. Combining \reff{bdfirst}, \reff{bdsecond}, 
and \reff{bdthird} or \reff{bdfourth}, we obtain \reff{bdbar}. This
concludes the proof of Lemma~\ref{Sgaux}, hence of Proposition~\ref{Sgestim}.
\QED

\medskip
The results of this section can easily be generalized to weighted
$L^p$ spaces with $p \neq 2$. For instance, the following bounds are 
used in Section~\ref{scaling} to control the nonlinearity in \reff{SV2}.

\begin{proposition}\label{LpLqestim} Let $1 \le q \le p \le \infty$,
$m \ge 0$ and $T > 0$. For all $\alpha \in \intplus^N$, there
exists $C > 0$ such that
\begin{equation}\label{LpLq}
  |b^m \partial^\alpha S(\tau)f|_p \le \frac{C}{a(\tau)^{\frac{N}{2}
  (\frac{1}{q}-\frac{1}{p})+\frac{|\alpha|}{2}}} \,|b^m f|_q\ ,
\end{equation}
where $b(\xi) = (1+|\xi|^2)^{1/2}$. 
\end{proposition}

\proof This estimate follows easily from \reff{DalphaS} and Young's 
inequality. \QED


\section{Bounds on the velocity field}
\label{velocity2d}

In this section, we study in more detail the relationship between the 
vorticity $w(\xi)$ and the velocity field $\vv(\xi)$ obtained from
$w(\xi)$ via the Biot-Savart law \reff{SBS2}. In particular, we obtain 
sharp estimates for the spatial decay of the velocity as $|\xi| \to \infty$. 
This information is systematically used in Section~\ref{applications} 
to derive properties of the solutions of the Navier-Stokes equation 
from the results we have on the vorticity. 

First, we give explicit formulas for the velocity fields corresponding 
to the first eigenfunctions of the linear operator $\cL$ defined in 
\reff{Ldef}. For our purposes in Section~\ref{applications}, it is 
sufficient to consider only the first three eigenvalues $\lambda_0 = 0$, 
$\lambda_1 = -1/2$, and $\lambda_2 = -1$. 

\smallskip\noindent
{\bf 1)} The first eigenvalue $\lambda_0 = 0$ is simple, with Gaussian 
eigenfunction
$$
   G(\xi) = \frac{1}{4\pi}\ \e^{-|\xi|^2/4}\ .
$$
The corresponding velocity field $\vv^G(\xi)$ is given by
\begin{equation}\label{vgfield}
  \vv^G(\xi) = \frac{1}{2\pi}\frac{\e^{-|\xi|^2/4}-1}{|\xi|^2}
  \pmatrix{\xi_2 \cr -\xi_1}\ .
\end{equation}
Remark that $|\vv^G(\xi)| \sim |\xi|^{-1}$ as $|\xi| \to \infty$, so that
$\vv^G \in L^q(\real^2)^2$ for all $q > 2$. 

\noindent
{\bf 2)} The second eigenvalue $\lambda_1 = -1/2$ is double, with 
eigenfunctions
$$
   F_i(\xi) = \partial_i G(\xi) = -\frac{\xi_i}{2} G(\xi)\ ,
   \quad i = 1,2\ .
$$
The corresponding velocity fields are 
\begin{equation}\label{vf1vf2}
  \vv^{F_1}(\xi) = \partial_1 \vv^G(\xi)\ , \quad 
  \vv^{F_2}(\xi) = \partial_2 \vv^G(\xi)\ .
\end{equation}
Note that $|\vv^{F_i}(\xi)| \sim |\xi|^{-2}$ as $|\xi| \to \infty$, so that
$\vv^{F_i} \in L^q(\real^2)^2$ for all $q > 1$. 

\smallskip\noindent 
{\bf 3)} The third eigenvalue $\lambda_1 = -1$ has multiplicity three. 
A convenient basis of eigenfunctions is
\begin{eqnarray*}
  H_1(\xi) &=& \Delta G(\xi) = (\textstyle{\frac{1}{4}}|\xi|^2 -1)G(\xi)\ ,\\
  H_2(\xi) &=& (\partial_1^2 - \partial_2^2)G(\xi) = \textstyle{\frac{1}{4}}
   (\xi_1^2 - \xi_2^2)G(\xi)\ ,\\
  H_3(\xi) &=& \partial_1 \partial_2 G(\xi) = \xi_1 \xi_2 G(\xi)\ ,  
\end{eqnarray*}
and the corresponding velocity fields read
\begin{eqnarray}\label{vhfield}
  \vv^{H_1}(\xi) & = & \Delta \vv^G(\xi) = \frac{1}{2}G(\xi)
    \pmatrix{\xi_2 \cr -\xi_1}\ ,\\ \nonumber
  \vv^{H_2}(\xi) & = & (\partial_1^2 - \partial_2^2) \vv^G(\xi)\ , \quad 
  \vv^{H_3}(\xi) \ =\  \partial_1 \partial_2 \vv^G(\xi)\ .
\end{eqnarray}
Remark that $\vv^{H_1}$ decreases rapidly as $|\xi| \to \infty$. 
In constrast, for $j = 2,3$, $|\vv^{H_j}(\xi)| \sim |\xi|^{-3}$ as 
$|\xi| \to \infty$, so that $\vv^{H_j} \in L^1(\real^2)^2$, but
$b \vv^{H_j} \notin L^1(\real^2)^2$. 

\medskip
Next, assume that $w \in L^2(m)$ for some $m > 3$. Then $w$ can be 
decomposed as follows:
\begin{equation}\label{wdecomp}
  w(\xi) = \alpha G(\xi) + \sum_{i=1}^2 \beta_i F_i(\xi) + 
  \sum_{j=1}^3 \gamma_j H_j(\xi) + R(\xi)\ ,
\end{equation}
where $\alpha, \beta_i$ are given by \reff{alphabetcoeff} and $\gamma_j$ 
by \reff{gammacoeff}. The corresponding velocity field has a similar 
decomposition
\begin{equation}\label{vdecomp}
  \vv(\xi) = \alpha \vv^G(\xi) + \sum_{i=1}^2 \beta_i \vv^{F_i}(\xi) + 
  \sum_{j=1}^3 \gamma_j \vv^{H_j}(\xi) + \vv^R(\xi)\ ,
\end{equation}
where $\vv^R$ is obtained from $R$ via the Biot-Savart law \reff{SBS2}.
By construction, the remainder term $R$ in \reff{wdecomp} satisfies
$\int_{\real^2} \xi_1^{n_1}\xi_2^{n_2} R(\xi)\d\xi = 0$ for all 
$n_1, n_2 \in \intplus$ such that $n_1 + n_2 \le 2$. Using this
information, we shall show that $|\vv^R(\xi)| \sim |\xi|^{-m}$ as 
$|\xi| \to \infty$. In view of \reff{wdecomp}, it will follow that 
$\vv \in L^2(\real^2)^2$ if and only if $\alpha = 0$, namely
\begin{equation}\label{CA}
  \int_{\real^2} w(\xi)\d\xi = 0\ .
\end{equation}
Moreover, if \reff{CA} holds, then $\vv \in L^1(\real^2)^2$ if and only 
if $\beta_1 = \beta_2 = 0$, namely
\begin{equation}\label{CB}
  \int_{\real^2} \xi_1 w(\xi)\d\xi = \int_{\real^2} \xi_2 w(\xi)\d\xi = 0\ .
\end{equation}
Finally, if \reff{CA} and \reff{CB} hold, then $b\vv \in L^1(\real^2)^2$
if and only if $\gamma_2 = \gamma_3 = 0$, namely
\begin{equation}\label{CC}
  \int_{\real^2} (\xi_1^2 -\xi_2^2) w(\xi)\d\xi = \int_{\real^2} \xi_1 
  \xi_2 w(\xi)\d\xi = 0\ .
\end{equation}
Remark that it is not necessary to assume here that $\gamma_1 = 0$, 
since the velocity field $\vv^{H_1}$ (unlike $\vv^{H_2}$ and $\vv^{H_3}$)
decreases rapidly at infinity. 

\begin{proposition}\label{velvort2}
Let $w \in L^2(m)$ for some $m > 0$, and denote by $\vv$ the velocity 
field obtained from $w$ via the Biot-Savart law \reff{SBS2}. Assume
that either\\
{\bf 1)} $0 < m \le 1$, or\\
{\bf 2)} $1 < m \le 2$ and \reff{CA} holds, or\\
{\bf 3)} $2 < m \le 3$ and \reff{CA}, \reff{CB} hold, or\\
{\bf 4)} $3 < m \le 4$ and \reff{CA}, \reff{CB}, \reff{CC} hold.\\
If $m \notin \intplus$, then for all $q \in (2,+\infty)$, there exists
$C > 0$ such that
\begin{equation}\label{velv1}
  |b^{m-2/q} \vv|_q \le C|b^m w|_2\ , 
\end{equation}
where $b(\xi) = (1+|\xi|^2)^{1/2}$. If $m \in \intplus$ and $b^m w \in 
L^p(\real^2) \cap L^r(\real^2)$ for some $p < 2$, $r > 2$, then
\begin{equation}\label{velv2}
  |b^m \vv|_\infty \le C (|b^m w|_p + |b^m w|_r)\ . 
\end{equation}
\end{proposition}

\medskip
The proof of Proposition~\ref{velvort2} relies on the following 
weighted Hardy-Littlewood-Sobolev inequality:

\begin{lemma}\label{HLSw}
If $0 < m < 1$ and 
$$
  u(\xi) = \int_{\real^2} \frac{\omega(\eta)}{|\xi-\eta|}\d\eta\ ,
  \quad \xi \in \real^2\ ,
$$
then, for all $q \in (2,+\infty)$, $|b^{m-2/q}u|_q \le C|b^m \omega|_2$. 
\end{lemma}

\proof We use the dyadic decomposition
$$
   \real^2 = \bigcup_{j=0}^\infty B_j\ ,
$$
where $B_0 = \{\xi \in \real^2\,|\,|\xi| \le 1\}$ and $B_j = \{\xi \in 
\real^2\,|\, 2^{j-1} < |\xi| \le 2^j\}$ for $j \in \intplus^*$. Let 
$u_i = u\oone_{B_i}$ and $\omega_i = \omega\oone_{B_i}$, $i \in \intplus$. 
Clearly $u_i = \sum_{j \in \intplus} \Delta_{ij}$, where
$$
   \Delta_{ij}(\xi) = \oone_{B_i}(\xi) \int_{B_j}\frac{\omega_j(\eta)}
   {|\xi-\eta|}\d\eta\ .
$$
Fix $q \in (2,+\infty)$, and define $p \in (1,2)$ by the relation 
$\frac{1}{q} = \frac{1}{p} -\frac{1}{2}$.

\smallskip\noindent
If $|i-j| \le 1$, it follows from \reff{HLS} that $|\Delta_{ij}|_q 
\le C|\omega_j|_p$. By H\"older's inequality, $|\omega_j|_p \le 
\mes(B_j)^{1/q}|\omega_j|_2 \le C 2^{2j/q}|\omega_j|_2$. Thus, 
$|\Delta_{ij}|_q \le C 2^{2j/q}|\omega_j|_2$.  

\smallskip\noindent
If $|i-j| \ge 2$, Young's inequality implies that $|\Delta_{ij}|_q \le 
M_1^{1-p/q} M_2^{p/q}|\omega_j|_2$, where
where
$$
   M_1 = \sup_{\xi\in B_i}\Bigl(\int_{B_j} \frac{1}{|\xi-\eta|^p}\d\eta
   \Bigr)^{1/p}\ , \quad 
   M_2 = \sup_{\eta\in B_j}\Bigl(\int_{B_i} \frac{1}{|\xi-\eta|^p}\d\xi
   \Bigr)^{1/p}\ .
$$
If $i \ge j+2$, then $|\xi-\eta| \ge |\xi|-|\eta| \ge 2^{i-1}-2^j \ge 
2^{i-2}$ for all $\xi \in B_i$, $\eta \in B_j$. Thus $M_1 \le 
C2^{-i}\mes(B_j)^{1/p} \le C2^{-i}2^{2j/p}$ and $M_2 \le C2^{-i}
\mes(B_i)^{1/p}\le C2^{-i}2^{2i/p}$, for some $C > 0$ independent of $i,j$.
It follows that $|\Delta_{ij}|_q \le C 2^{-(i-j)} 2^{2i/q}|\omega_j|_2$. 
If $j \ge i+2$, then $|\xi-\eta| \ge 2^{j-2}$ for all $\xi \in B_i$, 
$\eta \in B_j$, and a similar calculation shows that $|\Delta_{ij}|_q \le 
C 2^{2i/q}|\omega_j|_2$. Summarizing, we have shown that
$$
  2^{-2i/q} |u_i|_q \le C \sum_{j \in \intplus} K_{ij} |\omega_j|_2\ , 
  \quad i \in \intplus\ ,
$$
where $K_{ij} = 2^{-\frac{1}{2}|i-j|-\frac{1}{2}(i-j)}$. Now, by definition 
of the sets $B_i$, there exists $C \ge 1$ such that $C^{-1}2^{mi} \le 
b(\xi) \le C2^{mi}$ for all $\xi \in B_i$ and all $i \in \intplus$. 
It follows that 
$$
  |b^{m-2/q} u_i|_q \le C\sum_{j \in \intplus} K_{ij}^{(m)} |b^m 
  \omega_j|_2\ , \quad i \in \intplus\ ,
$$
where $K_{ij}^{(m)} = 2^{-\frac{1}{2}|i-j|+(m-\frac{1}{2})(i-j)}$. Since
$0 < m < 1$, $|K_{ij}^{(m)}| \le 2^{-\alpha|i-j|}$ for some $\alpha > 0$, 
hence $K^{(m)}$ defines a bounded linear operator from $\ell^2(\intplus)$ 
into $\ell^q(\intplus)$. This concludes the proof. \QED

\medskip\noindent{\bf Proof of Proposition~\ref{velvort2}.}
The proof is naturally divided into four steps, according to the 
values of $m$.\\
{\bf 1)} If $\vv$ and $w$ are related via the Biot-Savart law \reff{SBS2}, 
then
\begin{equation}\label{newBS}
  |\vv(\xi)| \le \frac{1}{2\pi} \int_{\real^2}
  \frac{1}{|\xi - \eta |} |w(\eta)| \d\eta\ , \quad \xi \in \real^2\ ,
\end{equation}
and \reff{velv1} follows immediately from Lemma~\ref{HLSw}. To prove
\reff{velv2}, we remark that $b(\xi) \le 1+|\xi| \le 1+|\xi-\eta|+|\eta|
\le |\xi-\eta| + 2b(\eta)$ for all $\xi,\eta \in \real^2$, hence
\begin{equation}\label{dummy}
  b(\xi)|\vv(\xi)| \le \frac{1}{2\pi} \int_{\real^2} |w(\eta)|\d\eta + 
   \frac{1}{\pi} \int_{\real^2} \frac{1}{|\xi-\eta|} b(\eta)
   w(\eta)\d\eta\ .
\end{equation}
In view of \reff{interpol}, the second integral in \reff{dummy} 
is uniformly bounded by $C(|b w|_p + |b w|_r)$ if $p < 2$ and $r > 2$.
On the other hand, it is clear that $|w|_1 \le C|b w|_p$. This concludes
the proof of case {\bf 1}. 

\smallskip\noindent{\bf 2)} For all $\xi,\eta \in \real^2$ with $\xi \neq 0$
and $\xi \neq \eta$, we have the identity
$$
  \frac{\xi_1-\eta_1}{|\xi-\eta|^2} - \frac{\xi_1}{|\xi|^2} =
  \frac{1}{|\xi|^2|\xi-\eta|^2}\Bigl((\xi_1-\eta_1)(\xxi\cdot\eeta) +
  (\xi_2-\eta_2)(\xxi\wedge\eeta)\Bigr)\ ,
$$
where $\xxi\cdot\eeta = \xi_1\eta_1 + \xi_2\eta_2$ and $\xxi\wedge\eeta = 
\xi_1\eta_2 - \xi_2\eta_1$. Therefore, for all $w \in L^2(m)$ satisfying
\reff{CA}, we find
\begin{eqnarray*}
 |v_2(\xi)| & = & \frac{1}{2\pi} \left| \int_{\real^2} 
 \left(\frac{\xi_1-\eta_1}{|\xi-\eta|^2} - \frac{\xi_1}{|\xi|^2}\right)
 w(\eta)\d\eta \right| \\
 & \le & C \int_{\real^2} \frac{1}{|\xi||\xi-\eta|} |\eta||w(\eta)|\d\eta\ .
\end{eqnarray*}
A similar bound holds for the first component $v_1(\xi)$ of the 
velocity field. Combining these estimates with \reff{newBS}, we obtain
$$
  |b(\xi)\vv(\xi)| \le C \int_{\real^2} \frac{1}{|\xi - \eta|} 
  |b(\eta)w(\eta)| \d\eta\ , \quad \xi \in \real^2\ .
$$
This inequality has exactly the same form as \reff{newBS}, with $\vv,w$
replaced by $b\vv,bw$. Therefore, applying the preceding result to 
$b\vv,bw$ instead of $\vv,w$, we obtain \reff{velv1} and \reff{velv2} in 
case {\bf 2}. 

\smallskip\noindent{\bf 3)} For all $\xi,\eta \in \real^2$ with $\xi \neq 0$
and $\xi \neq \eta$, we have the identity
\begin{eqnarray*}
  &&\frac{\xi_1-\eta_1}{|\xi-\eta|^2} - \frac{\xi_1-\eta_1}{|\xi|^2} 
  -\frac{2\xi_1(\xxi\cdot\eeta)}{|\xi|^4}\\ 
  &&\qquad\qquad
  = \frac{1}{|\xi|^4|\xi-\eta|^2}\Bigl((\xi_1-\eta_1)(2(\xxi\cdot\eeta)^2
  - |\xi|^2|\eta|^2) + 2(\xi_2-\eta_2)(\xxi\cdot\eeta)(\xxi\wedge\eeta)
\Bigr)\ .
\end{eqnarray*}
Therefore, for all $w \in L^2(m)$ satisfying \reff{CA} and \reff{CB}, we find
\begin{eqnarray*}
 |v_2(\xi)| & = & \frac{1}{2\pi} \left| \int_{\real^2} 
 \left(\frac{\xi_1-\eta_1}{|\xi-\eta|^2} - \frac{\xi_1-\eta_1}{|\xi|^2} 
  -\frac{2\xi_1(\xxi\cdot\eeta)}{|\xi|^4}\right)w(\eta)\d\eta \right| \\
 & \le & C \int_{\real^2} \frac{1}{|\xi|^2 |\xi-\eta|}|\eta|^2 
  |w(\eta)|\d\eta\ .
\end{eqnarray*}
A similar bound holds for the first component $v_1(\xi)$, hence
$$
  |b(\xi)^2\vv(\xi)| \le C \int_{\real^2} \frac{1}{|\xi - \eta|} 
  |b(\eta)^2w(\eta)| \d\eta\ , \quad \xi \in \real^2\ .
$$
This inequality has again the same form as \reff{newBS}, with $\vv,w$
replaced by $b^2\vv,b^2w$. Therefore, proceeding as above, we obtain 
\reff{velv1} and \reff{velv2} in case {\bf 3}.

\smallskip\noindent{\bf 4)} For all $\xi,\eta \in \real^2$ with $\xi \neq 0$
and $\xi \neq \eta$, we have the identity
\begin{eqnarray}\nonumber
  &&\frac{\xi_1-\eta_1}{|\xi-\eta|^2} - \frac{\xi_1-\eta_1}{|\xi|^2} 
  -\frac{2(\xi_1-\eta_1)(\xxi\cdot\eeta)}{|\xi|^4} + 
  \frac{\xi_1 |\eta|^2}{|\xi|^4} -\frac{4\xi_1(\xxi\cdot\eeta)^2}{|\xi|^6}
  \\ \label{ident}
  &&\hspace{30pt} = \frac{1}{|\xi|^6|\xi-\eta|^2}
  \Bigl((\xi_1-\eta_1)(\xxi\cdot\eeta)(4(\xxi\cdot\eeta)^2-3|\xi|^2|\eta|^2)
  \\ \nonumber
  &&\hspace{100pt} +\,(\xi_2-\eta_2)(\xxi\wedge\eeta)
  (4(\xxi\cdot\eeta)^2-|\xi|^2|\eta|^2)\Bigr)\ .
\end{eqnarray}
The left-hand side of \reff{ident} has the form $(\xi_1{-}\eta_1)/
|\xi{-}\eta|^2 + \Lambda(\xi,\eta)$, where $\Lambda$ is a polynomial of
degree $2$ in $\eta$ with $\xi$-dependent coefficients. As is easily 
verified, the terms in $\Lambda$ that are quadratic in $\eta$ can be
written as
$$
  \frac{1}{|\xi|^6}\Bigl(\xi_1(\xi_1^2-3\xi_2^2)(\eta_2^2-\eta_1^2)
  + 2\xi_2(\xi_2^2-3\xi_1^2)\eta_1\eta_2\Bigr)\ .
$$
For all $w \in L^2(m)$ satisfying \reff{CA}, \reff{CB} and
\reff{CC}, we thus have $\int_{\real^2} \Lambda(\xi,\eta)w(\eta)\d\eta 
\equiv 0$. Therefore, using \reff{ident}, we obtain
$$
  |v_2(\xi)| \le C \int_{\real^2} \frac{1}{|\xi|^3 |\xi-\eta|}|\eta|^3 
  |w(\eta)|\d\eta\ .
$$
A similar bound holds for the first component $v_1(\xi)$, hence
$$
  |b(\xi)^3 \vv(\xi)| \le C \int_{\real^2} \frac{1}{|\xi - \eta|} 
  |b(\eta)^3 w(\eta)| \d\eta\ , \quad \xi \in \real^2\ .
$$
This inequality has again the same form as \reff{newBS}, with $\vv,w$
replaced by $b^3\vv,b^3w$. Therefore, proceeding as above, we obtain 
\reff{velv1} and \reff{velv2} in case {\bf 4}. This concludes the proof 
of Proposition~\ref{velvort2}. \QED

\begin{corollary}\label{velint1}
Assume that $w \in L^2(m)$ for some $m > 1$, and denote by $\vv$ the 
velocity field obtained from $w$ via the Biot-Savart law \reff{SBS2}.
Then $\vv \in L^2(\real^2)^2$ if and only if \reff{CA} holds. 
\end{corollary}

\proof Without loss of generality, we assume that $1 < m < 2$. 
For any $w \in L^2(m)$, we have the decompositions $w = \alpha G +\bar w$, 
$\vv = \alpha \vv^G + \bar \vv$, where $\alpha = \int_{\real^2}w(\xi)\d\xi$
and $\bar w$ satisfies the assumptions of point~2 in 
Proposition~\ref{velvort2}. Setting $q = 2/(m{-}1)$ in \reff{velv1}, we 
obtain $b \bar \vv \in L^q(\real^2)^2$, which implies $\bar \vv \in 
L^2(\real^2)^2$. Since $\vv^G \notin L^2(\real^2)^2$, it is clear that
$\vv \in  L^2(\real^2)^2$ if and only if $\alpha = 0$. \QED

\begin{corollary}\label{velint2}
Assume that $w \in L^2(m)$ for some $m > 2$, and denote by $\vv$ the 
velocity field obtained from $w$ via the Biot-Savart law \reff{SBS2}.
Then $\vv \in L^1(\real^2)^2$ if and only if \reff{CA} and \reff{CB} 
hold. In that case, $\int_{\real^2}v_i(\xi)\d\xi = 0$ for $i = 1,2$. 
\end{corollary}

\proof Without loss of generality, we assume that $2 < m < 3$. 
For any $w \in L^2(m)$, we have the decompositions
$$
   w = \alpha G + \beta_1 F_1 + \beta_2 F_2 + \tilde w\ , \quad
   \vv = \alpha \vv^G + \beta_1 \vv^{F_1} + \beta_2 \vv^{F_2} 
     + \tilde \vv\ ,
$$
where $\alpha,\beta_1,\beta_2$ are given by \reff{alphabetcoeff} and 
$\tilde w$ satisfies the assumptions of point~3 in 
Proposition~\ref{velvort2}. Thus $b^2 \tilde\vv \in L^q(\real^2)^2$ 
with $q = 2/(m{-}2)$, hence $\tilde\vv \in L^1(\real^2)^2$. On the
other hand, it is easy to verify that $\vv - \tilde\vv \in L^1(\real^2)^2$
if and only if $\alpha = \beta_1 = \beta_2 = 0$. Thus, $\vv \in 
L^1(\real^2)^2$ if and only if \reff{CA} and \reff{CB} hold. 
In that case, since $\div \vv = 0$, we must have $\int_{\real^2}v_i(\xi)
\d\xi = 0$ for $i = 1,2$. \QED

\begin{corollary}\label{velint3}
Assume that $w \in L^2(m)$ for some $m >  3$, and denote by $\vv$ the 
velocity field obtained from $w$ via the Biot-Savart law \reff{SBS2}.
Then $b\vv \in L^1(\real^2)^2$ if and only if \reff{CA}, \reff{CB} and 
\reff{CC} hold. In that case, 
\begin{eqnarray*}
  \int_{\real^2}\xi_1 v_1(\xi)\d\xi & = & \int_{\real^2}\xi_2 v_2(\xi)
    \d\xi \ = \ 0\ ,\\
  \int_{\real^2}\xi_2 v_1(\xi)\d\xi & = & -\int_{\real^2}\xi_1 v_2(\xi)
    \d\xi \ = \ \gamma_1\ .\\
\end{eqnarray*}
\end{corollary}

\proof Without loss of generality, we assume that $3 < m < 4$. Let 
$w \in L^2(m)$, and consider the decomposition \reff{wdecomp}, 
\reff{vdecomp}. Then $R$ satisfies the assumptions of point~4 in 
Proposition~\ref{velvort2}. In particular, $b^3 \vv^R \in L^q(\real^2)^2$
for $q = 2/(m{-}3)$, hence $b\vv^R \in L^1(\real^2)^2$. On the other
hand, it is easy to verify that $b(\vv - \vv^R) \in L^1(\real^2)^2$
if and only if $\alpha = \beta_1 = \beta_2 = \gamma_2 = \gamma_3 = 0$. 
Thus, $b\vv \in L^1(\real^2)^2$ if and only if \reff{CA}, \reff{CB} and
\reff{CC} hold. 

Assume now that $b\vv \in L^1(\real^2)$, and consider the vector 
field $\AA(\xi) = \xi_1^2 \vv(\xi)$. Then $\AA \in L^2(\real^2)^2$ and 
$\div \AA = 2\xi_1v_1 \in L^1(\real^2)^2$, hence
$$
   0 = \int_{\real^2}\div \AA(\xi)\d\xi = \int_{\real^2} 2\xi_1 v_1(\xi)
   \d\xi\ .
$$
Similarly, using the vector fields $\xi_2^2 \vv$ and $\xi_1 \xi_2 \vv$, 
one obtains $\int_{\real^2}\xi_2 v_2\d\xi = 0$ and $\int_{\real^2}
(\xi_1v_2 + \xi_2v_1)\d\xi = 0$. Finally, if $\BB(\xi) = |\xi|^2 
\vv^\perp = |\xi|^2(-v_2,v_1)$, then $\BB \in L^2(\real^2)^2$ and 
$\div \BB = -|\xi|^2 w + 2(\xi_2 v_1 -\xi_1v_2) \in L^1(\real^2)^2$, 
hence $\int_{\real^2} \div \BB\d\xi = 0$. Using \reff{gammacoeff}, 
we thus find
$$
   4 \gamma_1 = \int_{\real^2} |\xi|^2 w(\xi) \d\xi = 
   2\int_{\real^2} (\xi_2 v_1(\xi) - \xi_1 v_2(\xi))\d\xi = 
   4 \int_{\real^2} \xi_2 v_1(\xi)\d \xi\ .
$$
This concludes the proof of Corollary~\ref{velint3}. \QED


\section{Proof of Proposition~\ref{invman2}} \label{Gapp}
The existence and attractivity of the manifold follow from the
results of \cite{chen:1997} exactly as in the proof of Theorem~\ref{invman}
and Corollary~\ref{locinvman}.  By construction, the function
$\cG$ whose graph gives the invariant manifold satisfies the
integral equation (see, for instance, \cite{gallay:1993}):
\begin{equation}\label{Gint}
  \cG(\bbeta,\ggamma) = \int_0^\infty \e^{\tau \cL}
  \cN_0(\e^{\tau/2} \bbeta,\phi^{-\tau}(\ggamma;\bbeta),
  \cG(\e^{\tau/2} \bbeta,\phi^{-\tau}(\ggamma;\bbeta)) )\d\tau\ ,
\end{equation}
where $(\bbeta,\ggamma) \mapsto (\e^{-\tau/2}\bbeta, 
\phi^\tau(\ggamma;\bbeta))$ is (the projection onto $E_c$ of) the flow 
defined by equations \reff{log2} on $\wcloc$, and $\cN_0$ is the nonlinear
term in the equation for $\rho$ in \reff{log2}, cutoff outside a neighborhood
of the origin as in \reff{SV2c}.  Note that in computing 
$\phi^\tau(\ggamma;\bbeta)$, one can first solve the equations for 
$\bbeta$ explicitly and then insert these solutions into the equations for 
$\gamma_j$, thus treating them as inhomogeneous terms.

\begin{remark}\label{gflowdef} Note that the nonlinear terms $\tf_j$
in the equations for $\gamma_j$ are also cutoff so that they vanish outside 
a neighborhood of the origin of size $r_0$. It is an easy exercise to show 
that for any $\epsilon>0$, there exists $r_0 > 0$ sufficiently small so 
that the resulting flow $\phi^\tau(\ggamma;\bbeta)$ satisfies
\begin{equation}\label{gflowest}
  |\phi^{-\tau}(\ggamma;\bbeta) | \le C \e^{(1+\epsilon)\tau}(|\ggamma|
  +|\bbeta|^2)\ ,
\end{equation}
for all $\tau\ge 0$.
\end{remark}
Define a mapping
\begin{equation} \label{cFdef}
  \cF(\cG)(\bbeta,\ggamma) = \int_0^{\infty} \e^{\tau\cL} 
  \cN_0(\e^{\tau/2} \bbeta,\phi^{-\tau}(\ggamma;\bbeta),
  \cG(\e^{\tau/2} \bbeta,\phi^{-\tau}(\ggamma;\bbeta)) ) \d\tau\ .
\end{equation}
Equation \reff{Gint} implies that $\cF$ has a fixed point in the space of 
$C^1$ functions with globally bounded Lipshitz constant. Moreover, if one 
fixes the cutoff function $\chi_{r_0}$ in \reff{SV2c}, then (for $r_0 > 0$
sufficiently small) $\cF$ is a contraction in this function space, so that
$\cG$ is the unique fixed point of $\cF$. Thus, the proof of the proposition 
is completed if we can show that for some $C_\delta > 0$, $\cF$ maps the 
set of functions satisfying \reff{cGest} into itself.

Let $X$ be the Banach space $\real^2 \times \real^3 \times E_s$,
equipped with the norm $\|(\bbeta,\ggamma,\rho)\|_X = |\bbeta|+|\ggamma|
+\|\rho\|_4$, and let $\chi_{r_0}$ be a cutoff function like that introduced
just prior to \reff{SV2c}, but defined with respect to the norm $\|\cdot\|_X$ 
rather than $\| \cdot \|_m$. Then the nonlinear term in \reff{Gint} can be 
written as
\begin{equation}\label{ncdef}  
  \cN_0(\bbeta,\ggamma,\rho) =
  Q_s \{ \nabla \cdot (\chi_{r_0}(\bbeta,\ggamma,\rho) 
  \FF(\bbeta,\ggamma,\rho)) - \vv^\rho \cdot \nabla (\chi_{r_0}
  (\bbeta,\ggamma,\rho) \rho) \}\ ,
\end{equation}
where $Q_s$ is the projection onto $E_s$. Commuting the derivatives
in the nonlinear term through the factor of $\e^{\tau\cL}$ as
we did in \reff{winteq}, the right-hand side of \reff{cFdef}
can be rewritten as
\begin{eqnarray*}
&& \int_0^1 Q_s \nabla\cdot\e^{\tau(\cL-\frac{1}{2})}
   [(\chi_{r_0}\FF)(\bbeta(\tau),\ggamma(\tau),\cG(\bbeta(\tau),
   \ggamma(\tau)))\\ 
&& \qquad - \vv^{\cG(\beta(\tau),\gamma(\tau))} \chi_{r_0}(\bbeta(\tau),
   \ggamma(\tau),\cG(\bbeta(\tau),\ggamma(\tau)))\cG(\bbeta(\tau),
   \ggamma(\tau))] \d\tau \\ 
&& +\int_1^{\infty} \e^{(\tau-1)\cL} Q_s \nabla \cdot \e^{\cL-\frac{1}{2}}
   [(\chi_{r_0}\FF)(\bbeta(\tau),\ggamma(\tau),\cG(\bbeta(\tau),
   \ggamma(\tau))) \\ 
&& \qquad -  \vv^{\cG(\beta(\tau),\gamma(\tau))} (\chi_{r_0}(\bbeta(\tau),
   \ggamma(\tau),\cG(\bbeta(\tau),\ggamma(\tau)))\cG(\bbeta(\tau),
   \ggamma(\tau) ) ] \d\tau\ ,
\end{eqnarray*}
where in keeping with our usual notational convention, $\vv^{\cG}$ denotes 
the velocity field associated with the vorticity field $\cG$.  Also, to 
save space, we have used the short-hand notation $\bbeta(\tau)$ for 
$\e^{\tau/2}\bbeta$ and $\ggamma(\tau)$ for $\phi^{-\tau}(\ggamma;\bbeta)$.
Fixing $\epsilon > 0$, we can take the $\|\cdot\|_4$ norm of this expression 
and use the estimates on the semigroup from Appendix \ref{spectrum} to 
obtain
\begin{eqnarray}\nonumber
&& \|\cF(\cG)(\bbeta,\ggamma)\|_4 \le  C\int_0^{\infty}
  \e^{-(\frac{3}{2} -\epsilon)(\tau-1)} [(\chi_{r_0} \|\FF\|_4)
  (\bbeta(\tau),\ggamma(\tau),\cG(\bbeta(\tau),\ggamma(\tau))) \|_4 
  \\ \label{cFbound1}
&& \qquad \qquad +\, \chi_{r_0}(\bbeta(\tau),\ggamma(\tau),\cG(\bbeta(\tau),
  \ggamma(\tau)))\|\cG(\bbeta(\tau),\ggamma(\tau))\|_4^2] \d\tau\ ,
\end{eqnarray}
where we estimated the term involving $\vv^{\cG} (\chi_{r_0}(\bbeta,\ggamma,
\cG(\bbeta,\ggamma))\cG)$ in the same way as we estimated 
the term involving $\vv(s) w(s)$ in Lemma~\ref{Restim}. Using the estimate 
on $\FF$ in \reff{ffbounds}, we see that
$$
  \|\chi_{r_0}(\bbeta,\ggamma,\cG) 
  \FF(\bbeta,\ggamma,\cG) \|_4 \le C ((|\bbeta| +|\ggamma|)
  (|\ggamma|+\|\cG\|_4) +\|\cG\|_4^2)\chi_{r_0}(\bbeta,\ggamma,\cG) \ .
$$
If we now assume that $\|\cG(\bbeta,\ggamma)\|_4 \le C_{\delta} 
(|\bbeta|^{3-\delta} +|\ggamma|^{\frac{3}{2}-\delta})$, then, using the 
fact that $\chi_{r_0}(\bbeta,\ggamma,\rho)=0$ if either $|\bbeta|$ or
$|\ggamma|$ is bigger than $r_0$, it is easy to verify that 
\begin{equation}
  \|\chi_{r_0}(\bbeta,\ggamma,\cG)\FF(\bbeta,\ggamma,\cG)\|_4 
  \le C r_0^\delta (|\bbeta|^{3-\delta} + |\ggamma|^{\frac{3}{2}-\delta})\ .
\end{equation}
In a similar fashion we can bound
$$
  \chi_{r_0}(\bbeta,\ggamma,\cG(\bbeta,\ggamma))\|\cG(\bbeta,\ggamma))\|_4^2
  \le C r_0^\delta (|\bbeta|^{3-\delta} + |\ggamma|^{\frac{3}{2}-\delta})\ .
$$
Inserting these estimates into \reff{cFbound1} and using 
\reff{gflowest} to bound $\gamma(\tau) = \phi^{-\tau}(\ggamma;\bbeta)$, 
we find
\begin{eqnarray}\nonumber \label{cFbound2}
 \| \cF(\cG)(\bbeta,\ggamma) \|_4 &\le& C r_0^{\delta} \int_0^{\infty} 
  \e^{-(\frac{3}{2}-\epsilon)(\tau-1)} [\,\e^{(3-\delta)(\tau/2)} 
  |\bbeta|^{3-\delta} +\e^{(\frac{3}{2}-\delta)(1+\epsilon)\tau} 
  |\ggamma|^{\frac{3}{2}-\delta} ]\d\tau\\ 
 &\le& C(\epsilon,\delta) r_0^{\delta} (|\bbeta|^{3-\delta}
  + |\ggamma|^{\frac{3}{2}-\delta} )\ ,
\end{eqnarray}
provided that $\epsilon<2\delta/5$.  For fixed $\epsilon$ and $\delta$
(satisfying $\epsilon<2\delta/5$) we can always pick $r_0$
sufficiently small so that $C(\epsilon,\delta) r_0^\delta < C_\delta$.
Thus $\cF$ maps the set of functions satisfying \reff{cGest}
into itself (for $r_0$ sufficiently small) and hence the fixed
point of $\cF$ whose graph gives the invariant manifold must satisfy
\reff{cGest}.  This completes the proof of Proposition~\ref{invman2}.
\QED

\bibliographystyle{plain}
\bibliography{ref}

\end{document}